\documentclass[10pt]{article}
\usepackage{amscd}
\usepackage{amsmath}
\usepackage{latexsym}
\usepackage{amsfonts}
\usepackage{amssymb}
\usepackage{amsthm}
\usepackage{graphicx}

\newtheorem{theorem}{Theorem}[section]
\newtheorem{lemma}{Lemma}[section]

\newtheorem{proposition}{Proposition}[section]
\newtheorem{definition}{Definition}[section]

\newtheorem{corollary}{Corollary}[section]
\newtheorem{remark}{Remark}[section]

\makeatletter
\newcommand{\Extend}[5]{\ext@arrow0099{\arrowfill@#1#2#3}{#4}{#5}}
\makeatother

\begin{document}

 \title{ Global well-posedness and scattering for the energy-critical, defocusing Hartree equation in $\mathbb{R}^{1+n}$}
 \author{{Changxing Miao$^{\dagger}$,\ \ Guixiang Xu$^{\dagger}$,\ \ and \ Lifeng Zhao$^{\ddagger}$ }\\
         {\small $^{\dagger}$Institute of Applied Physics and Computational Mathematics}\\
         {\small P. O. Box 8009,\ Beijing,\ China,\ 100088}\\
         {\small $^{\ddagger}$Department of Mathematics, University of Science and Technology of
         China}\\
         {\small Hefei, China, 230026}\\
         {\small (miao\_changxing@iapcm.ac.cn, \ xu\_guixiang@iapcm.ac.cn, zhao\_lifeng@iapcm.ac.cn ) }\\
         \date{}
        }
\maketitle

\begin{abstract}
Using the same induction on energy argument in both frequency
space and spatial space simultaneously as in \cite{CKSTT07},
\cite{RyV05} and \cite{Vi05}, we obtain global well-posedness and
scattering of energy solutions of defocusing energy-critical
nonlinear Hartree equation in $\mathbb{R}\times
\mathbb{R}^n$($n\geq 5$), which removes the radial assumption on
the data in \cite{MiXZ07a}. The new ingredients are that we use a
modified long time perturbation theory to obtain the frequency
localization (Proposition \ref{freqdelocaimplystbound} and
Corollary \ref{frequencylocalization}) of the minimal energy blow
up solutions, which can not be obtained from the classical long
time perturbation and bilinear estimate and that we obtain the
spatial concentration of minimal energy blow up solution after
proving that $L^{\frac{2n}{n-2}}_x$-norm of minimal energy blow up
solutions is bounded from below, the $L^{\frac{2n}{n-2}}_x$-norm
is larger than the potential energy.
\end{abstract}

 \begin{center}
 \begin{minipage}{140mm}
 \end{minipage}
 \end{center}


{\center
\section{Introduction}
} \setcounter{section}{1}\setcounter{equation}{0}

In this paper, we consider the following initial value problem
\begin{equation} \label{equ1}
\left\{ \aligned
    iu_t +  \Delta u  & = f(u), \quad  \text{in}\  \mathbb{R}^n \times \mathbb{R}, \quad n\geq 5,\\
     u(0)&=u_0(x), \ \ \text{in} \ \mathbb{R}^n.
\endaligned
\right.
\end{equation}
where $u(t,x)$ is a complex-valued function in spacetime $\mathbb{R}
\times \mathbb{R}^n$ and $\Delta$ is the Laplacian in
$\mathbb{R}^{n}$, $f(u)=\big( |x|^{-4}* |u|^2 \big)
u=(|\nabla|^{-(n-4)}|u|^2 ) u$. It is introduced as a classical
model in \cite{website}. In practice, we use the integral
formulation of $(\ref{equ1})$
\begin{equation}\label{intequa1}
u(t)=U(t)u_0(x) -i \int^{t}_{0} U(t-s)f(u(s))ds,
\end{equation}
where $U(t)=e^{it\Delta}$.

This equation has the Hamiltonian
\begin{equation}\label{energy}
\aligned E(u(t))=\frac12 & \big\|\nabla
u(t)\big\|^2_{L^2}+\frac{1}{4} \iint \frac{1}{|x-y|^{4}} |u(t,x)|^2
|u(t,y)|^2\ dxdy.
\endaligned
\end{equation}
Since $(\ref{energy})$ is preserved by the flow corresponding to
$(\ref{equ1})$ we shall refer to it as the energy and often write
$E(u)$ for $E(u(t))$.

We are primarily interested in $(\ref{equ1})$ since it is critical
with respect to the energy norm. That is, the scaling $u\mapsto
u_{\lambda}$ where
\begin{equation}\label{scaling}
u_{\lambda}(t,x)=\lambda^{\frac{n-2}{2}}u(\lambda^2 t, \lambda x), \
\lambda>0
\end{equation}
leaves the energy invariant, in other words, the energy $E(u)$ is a
dimensionless quantity.

As is well-known that if the initial data $u_0(x)$ has finite
energy, then $(\ref{equ1})$ is locally well-posed (see, for instance
\cite{Ca03}, \cite{MiXZ06}). That is, there exists a unique
local-in-time solution that lies in $C^0_t\dot{H}^1_x \cap
L^6_tL^{\frac{6n}{3n-8}}_x$ and the map from the initial data to the
solution is locally Lipschitz in these norms. If the energy is
small, it is known that the solution exists globally in time and
scattering occurs; However, for initial data with large energy, the
local well-posedness argument do not extend to give global
well-posedness, only with the conservation of the energy
$(\ref{energy})$, because the lifespan of existence given by the
local theory depends on the profile of the data as well as on the
energy.

A large amount of works have been devoted to the theory of
scattering for Hartree equation, see \cite{ChHKY08},
\cite{GiO93}-\cite{HaT87}, \cite{LiMZ08}, \cite{Mi97} and
\cite{MiXZ07a}-\cite{NaO92}. In particular, global well-posedness
and scattering of $(\ref{equ1})$ with radial data in $\dot{H}^1_x$
was obtained in \cite{MiXZ07a} by taking advantage of the term
$\displaystyle - \int_{I}\int_{|x|\leq A|I|^{1/2}}|u|^{2}\Delta
\Big(\frac{1}{|x|}\Big)dxdt$ in the localized Morawetz identity to
rule out the possibility of energy concentration at origin. In this
paper, we continue this investigation. In order to prevent
concentration at any location in spacetime, we should take advantage
of the interaction Morawetz estimate achieved in \cite{CKSTT04}, or
the frequency-localized interaction Morawetz estimate in
\cite{CKSTT07}, see also \cite{RyV05}, \cite{Vi05}, \cite{Vi06},
etc.

Here, we give their brief differences in the case of the defocusing
Schr\"{o}dinger equation. After proving the negative regularity of
soliton solutions and double low-to-high frequency cascade solutions
(some kinds of minimal energy blow up solutions or almost periodic
solutions modulo symmetries), we can utilize the interaction
Morawetz estimate to prevent the concentration of them at any
location. While the negative regularity of soliton solutions and
double low-to-high frequency cascade solutions can be obtained under
the additional assumption of spatial dimension $n\geq 5$ due to the
fact that the Schr\"{o}dinger dispersion is not strong enough to
perform the double Duhamel trick for the low dimensions $n=3,4$. But
we can utilize the frequency localized interaction Morawetz estimate
to prevent the concentration of them at any location in low
dimensions as well as in high dimensions. See details in
\cite{CKSTT07}, \cite{KiV08}, \cite{RyV05}, \cite{Vi05},
\cite{Vi06}, etc.

Together with the frequency-localized interaction Morawetz estimate,
we will use the same induction on energy argument in both frequency
space and spatial space simultaneously as in \cite{CKSTT07} to
obtain global well-posedness and scattering for general large data,
which removes the radial assumption in \cite{MiXZ07a}. As for
induction on energy argument, we can also refer to \cite{B99a},
\cite{RyV05}, \cite{Vi05} and \cite{Vi06}. Induction on energy
argument is quantitative. In contrast with this method, D. Li, C.
Miao and X. Zhang \cite{LiMZ08} recently use concentration
compactness principle to obtain the similar result, that method is
qualitative and firstly introduced by Kenig and Merle \cite{KeM06}
to deal with the global well-posedness and scattering for focusing
energy-critical NLS. There are many applications in this direction,
for example \cite{KeM07}, \cite{KiTV}, \cite{KiV08}, \cite{KiVZ08},
etc.

However, the stability theory for the equation $(\ref{equ1})$ is an
essential tool for induction on energy argument. In the frame work
of the classical long time perturbation, we inevitably demand to
control the non-local interaction between the low and high
frequencies
\begin{equation*}
\aligned \big\| (|\nabla|^{-(n-4)}|u_{lo}|^2 ) u_{hi}
\big\|_{L^{\frac32}_t\dot{H}_x^{1, \frac{6n}{3n+4}}(I\times
\mathbb{R}^n)}
\endaligned
\end{equation*}
with
\begin{equation}\label{contol} \aligned
 \big\|u_{lo}\big\|_{\dot{S}^{1+s}(I \times \mathbb{R}^n)}
\leq C(\eta)\epsilon^{s}, \quad  \big\|u_{hi}\big\|_{\dot{S}^{1-k}(I
\times \mathbb{R}^n)} \leq C(\eta)\epsilon^{k}, \quad \forall \
0\leq s, k \leq 1.
\endaligned
\end{equation}
where the definition of norm $\big\|\cdot\big\|_{\dot{S}^1}$ refers
to $(\ref{defs})$. Because $|\nabla|^{-(n-4)}$ destroys the direct
interaction between $u_{lo}$ and $u_{hi}$ and we cannot use bilinear
estimate to obtain any decay even though we have the estimate
$(\ref{contol})$. It is different from the local interaction case of
the Schr\"{o}digner equation \cite{CKSTT07}, \cite{RyV05} and
\cite{Vi05}. We only have
\begin{equation*}
\aligned \big\| (|\nabla|^{-(n-4)}|u_{lo}|^2 ) u_{hi}
\big\|_{L^{\frac32}_t\dot{H}_x^{1, \frac{6n}{3n+4}}(I\times
\mathbb{R}^n)} \leq C(\eta).
\endaligned
\end{equation*}
No decay dues to one derivative on spatial variable in the spacetime
space $L^{\frac32}_t\dot{H}_x^{1, \frac{6n}{3n+4}}(I\times
\mathbb{R}^n)$. In deed, because we can not use the bilinear
estimate when one derivative falls on $u_{hi}$, there is no any
decay by $(\ref{contol})$ . However, when we would like to transfer
some part of derivative to integral on spatial variable, we can
obtain the small interaction
\begin{equation*}
\aligned \big\| (|\nabla|^{-(n-4)}|u_{lo}|^2 ) u_{hi}
\big\|_{L^{\frac32}_t(I; \dot{H}^{1-\frac{2}{n},
\frac{6n^2}{3n^2+4n-12}}_x)}  \leq C(\eta)\epsilon^{\frac{4}{n} }
\endaligned
\end{equation*}
according to $(\ref{contol})$. Inspired by this fact together with
the inhomogeneous Strichartz estimate \cite{Fos05}, \cite{Ka94} and
\cite{Vile07}, we set down a modified long time perturbation, which
replaces the role of the classical long time perturbation and the
bilinear estimate in some senses and is important to establish the
frequency localization (Proposition \ref{freqdelocaimplystbound} and
Corollary \ref{frequencylocalization}) of minimal energy blow up
solutions. See details in Section 4.

In addition, we obtain the spatial concentration of minimal energy
blow up solution after we prove that $L^{\frac{2n}{n-2}}_x$-norm of
minimal energy blow up solutions is bounded from below, which is
stronger than the statement that the potential energy of minimal
energy blow up solution is bounded from below.

Now, we give the main result of this paper.
\begin{theorem}\label{main}
Let $n\geq 5$. For any $u_0$ with finite energy, $E(u_0)<\infty$,
there exists a unique global solution $u\in C^0_t(\dot{H}^1_x) \cap
L^6_t(L^{\frac{6n}{3n-8}}_x)$ to $(\ref{equ1})$ such that
\begin{equation}\label{bound}
\aligned \big\| u\big\|_{L^6_tL^{\frac{6n}{3n-8}}_x(\mathbb{R}\times
\mathbb{R}^n)} \leq C(E(u_0))
\endaligned
\end{equation}
for some constant $C(E(u_0))$ that depends only on the energy.
\end{theorem}

As is well-known, the $L^{6}_tL^{\frac{6n}{3n-8}}_x$ bound above
also gives scattering, asymptotic completeness, and uniform
regularity.

\begin{corollary}
Let $u_0$ have finite energy. Then there exist finite energy
solutions $u_{\pm}(t,x)$ to the free Schr\"{o}dinger equation
$(i\partial_t + \Delta ) v  =0$ such that
\begin{equation*}
\aligned \big\|u(t)-u_{\pm}(t)\big\|_{\dot{H}^1} \rightarrow 0 \quad
\text{as} \quad t\rightarrow\pm \infty.
\endaligned
\end{equation*}
Furthermore, the maps $u_0 \mapsto u_{\pm}(0)$ are homeomorphisms
from $\dot{H}^1$ to $\dot{H}^1$. Finally, if $u_0 \in H^s$ for some
$s>1$, then $u(t)\in H^s$ for all time $t$, and one has the uniform
bounds
\begin{equation*}
\aligned \sup_{t\in \mathbb{R}} \big\|u(t)\big\|_{H^s} \leq
C(E(u_0), s)\big\|u_0\big\|_{H^s}.
\endaligned
\end{equation*}
\end{corollary}

Next, we introduce some notations.  If $X, Y$ are nonnegative
quantities, we use $X\lesssim Y $ or $X=O(Y)$ to denote the estimate
$X\leq CY$ for some $C$ which may depend on the critical energy
$E_{crit}$ (see Section 3) but not on any parameter such as $\eta$,
and $X \thicksim Y$ to denote the estimate $X\lesssim Y\lesssim X$.
We use $X\ll Y$ to mean $X \leq c Y$ for some small constant $c$
which is again allowed to depend on $E_{crit}$. We use $C\gg1$ to
denote various large finite constants, and $0< c \ll 1$ to denote
various small constants.

We use $L^q_tL^r_x$ to denote the spacetime norm
\begin{equation*}
\aligned \big\|u\big\|_{L^q_tL^r_x(\mathbb{R} \times
\mathbb{R}^n)}:=\Big( \int_{\mathbb{R}}\Big(
\int_{\mathbb{R}^n}\big|u(t,x)\big|^rdx\Big)^{q/r}dt\Big)^{1/q}
\endaligned
\end{equation*}
with the usual modifications when $q$ or $r$ is infinity, or when
the domain $\mathbb{R} \times \mathbb{R}^n$ is replaced by some
smaller spacetime region. When $q=r$, we abbreviate $L^q_tL^r_x$ by
$L^q_{t,x}$.

When $n\geq 5$, we say that a pair $(q, r)$ is sharp admissible if
\begin{equation*}
\frac{2}{q} = n\Big(\frac{1}{2}-\frac{1}{r}\Big),\ 2 \leq r \leq
 \frac{2n}{n-2}.
\end{equation*}
We say that the pair $(q,r)$ is acceptable if
\begin{equation*}
\aligned 1\leq q,r \leq \infty, \quad \frac{1}{q} < n
(\frac12-\frac1r), \quad \text{or} \ \ (q,r)=(\infty, 2).
\endaligned
\end{equation*}
For a spacetime slab $I\times \mathbb{R}^n$, we define the {\it
Strichartz} norm $\dot{S}^0(I)$ by
\begin{equation*}
\big\|u\big\|_{\dot{S}^0(I)}:= \sup_{(q, r)\ \text{sharp
admissible}} \Big(\sum_{N}\big\|P_{N}u\big\|^2_{L^q_tL^r_x(I\times
\mathbb{R}^n)}\Big)^{1/2}.
\end{equation*}
and for $k>0$, we define $\dot{S}^k(I)$ by
\begin{equation}\label{defs}
\big\|u\big\|_{\dot{S}^k(I)}:= \big\||\nabla|^k
u\big\|_{\dot{S}^0(I)}.
\end{equation}
From the Littlewood-Paley inequality, Sobolev embedding and
Minkowski's inequality, we have
\begin{equation}\label{lemstri2}
\aligned \big\|\nabla u\big\|_{L^{\infty}_tL^2_x} + \big\|\nabla u
\big\|_{L^6_tL^{L^{\frac{6n}{3n-2}}}_x}  + \big\| \nabla u
\big\|_{L^{\frac{6(n-2)}{n}}_tL^{\frac{6(n-2)}{3n-8}}_x} +
\big\|\nabla u\big\|_{L^3_tL^{\frac{6n}{3n-4}}_x} + &\big\|\nabla u\big\|_{L^2_tL^{\frac{2n}{n-2}}_x}  \\
+ \big\| u\big\|_{L^{\infty}_tL^{\frac{2n}{n-2}}_x} + \big\|
u\big\|_{L^6_tL^{\frac{6n}{3n-8}}_x} + \big\|
u\big\|_{L^4_tL^{\frac{2n}{n-3}}_x} +\big\| u
\big\|_{L^3L^{\frac{6n}{3n-10}}}+ \big\|
u\big\|_{L^2_tL^{\frac{2n}{n-4}}_x}& \lesssim \big\|
u\big\|_{\dot{S}^1},
\endaligned
\end{equation}
where all spacetime norms are taken on $I\times \mathbb{R}^n$.

The Fourier transform on $\mathbb{R}^n$ is defined by
\begin{equation*}
\aligned \widehat{f}(\xi):= \big( 2\pi
\big)^{-\frac{n}{2}}\int_{\mathbb{R}^n}e^{- ix\cdot \xi}f(x)dx ,
\endaligned
\end{equation*}
giving rise to the fractional differentiation operators
$|\nabla|^{s}$,  defined by
\begin{equation*}
\aligned \widehat{|\nabla|^sf}(\xi):=|\xi|^s\widehat{f}(\xi).
\endaligned
\end{equation*}
These define the homogeneous Sobolev norms
\begin{equation*}
\big\|f\big\|_{\dot{H}^s_x}:= \big\| |\nabla|^sf
\big\|_{L^2_x(\mathbb{R}^n)}.
\end{equation*}

Let $e^{it\Delta }$ be the free Schr\"{o}dinger propagator. This
propagator preserves the above Sobolev norms and obeys the
dispersive estimate
\begin{equation}\label{dispersivee}
\aligned \big\|e^{it\Delta} f\big\|_{L^{\infty}_x(\mathbb{R}^n)}
\lesssim |t|^{-\frac{n}{2}}\big\|f\big\|_{L^1_x(\mathbb{R}^n)}
\endaligned
\end{equation}
for all times $t\not=0$. We also recall {\it Duhamel's} formula
\begin{equation}\label{duhamelfor}
\aligned
u(t)=e^{i(t-t_0)\Delta}u(t_0)-i\int^t_{t_0}e^{i(t-s)\Delta}f(u)(s)ds.
\endaligned
\end{equation}

We will occasionally use subscripts to denote spatial derivatives
and will use the summation convention over repeated indices.

We will also need the Littlewood-Paley projection operators.
Specifically, let $\varphi(\xi)$ be a smooth bump function adapted
to the ball $|\xi|\leq 2$ which equals 1 on the ball $|\xi|\leq 1$.
For each dyadic number $N\in 2^{\mathbb{Z}}$, we define the
Littlewood-Paley operators
\begin{equation*}
\aligned \widehat{P_{\leq N}f}(\xi)& :=
\varphi\Big(\frac{\xi}{N}\Big)\widehat{f}(\xi), \quad \widehat{P_{>
N}f}(\xi) :=
\Big(1-\varphi\Big(\frac{\xi}{N}\Big)\Big)\widehat{f}(\xi), \\
\widehat{P_{N}f}(\xi)& :=
\Big(\varphi\Big(\frac{\xi}{N}\Big)-\varphi\Big(\frac{2\xi}{N}\Big)\Big)\widehat{f}(\xi).
\endaligned
\end{equation*}
Similarly we can define $P_{<N}$, $P_{\geq N}$, and $P_{M<\cdot\leq
N}=P_{\leq N}-P_{\leq M}$, whenever $M$ and $N$ are dyadic numbers.
We will frequently write $f_{\leq N}$ for $P_{\leq N}f$ and
similarly for the other operators.

The Littlewood-Paley operators commute with derivative operators,
the free propagator, and the conjugation operation. They are
self-adjoint and bounded on every $L^p_x$ and $\dot{H}^s_x$ space
for $1\leq p\leq \infty$ and $s\geq 0$. They also obey the following
Sobolev and Bernstein estimates
\begin{equation*}\label{bernstein}
\aligned
 \big\| P_{\geq N} f \big\|_{L^p} & \lesssim   N^{-s} \big\|
|\nabla|^{s}P_{\geq N} f \big\|_{L^p}, \\
\big\||\nabla|^s P_{\leq N} f \big\|_{L^p} & \lesssim    N^{s}
\big\| P_{\leq N} f \big\|_{L^p}, \quad \qquad \big\| P_{\leq N} f
\big\|_{L^q} \lesssim   N^{\frac{n}{p}-\frac{n}{q}} \big\|
P_{\leq N} f \big\|_{L^p},  \\
\big\||\nabla|^{\pm s} P_{N} f \big\|_{L^p} & \thicksim   N^{\pm s}
\big\| P_{N} f \big\|_{L^p}, \quad \ \ \qquad \big\| P_{ N} f
\big\|_{L^q} \lesssim   N^{\frac{n}{p}-\frac{n}{q}} \big\|P_{ N} f
\big\|_{L^p},
\endaligned
\end{equation*} whenever $s\geq 0$ and $1\leq p\leq q \leq \infty$.
Note that the kernel of the operator $P_{\leq N}$ is not positive.
To overcome this problem, we use the operator $P'_{\leq N}$ in
\cite{Vi05}, etc. More precisely, if $K_{\leq N}$ is the kernel
associated to $P_{\leq N}$, we let $P'_{\leq N}$ be the operator
associated to $N^{-n}\big( K_{\leq N} \big)^2$. The kernel of $
P'_{\leq N}$ is bounded in $L^1_x$ independently of $N$. Therefore,
the operator  $P'_{\leq N}$ is bounded on every $L^p_x$ for $1\leq
p\leq \infty$. Furthermore, for $s\geq 0$ and $1\leq p\leq q \leq
\infty$, we have
\begin{equation*}
\aligned \big\| | \nabla \big|^s P'_{\leq N}f\big\|_{L^p_x} &
\lesssim N^s \big\|P'_{\leq N} f\big\|_{L^p_x},\quad \big\|P'_{\leq
N} f\big\|_{L^q_x} & \lesssim N^{\frac{n}{q}-\frac{n}{p}}
\big\|P'_{\leq N} f\big\|_{L^p_x}.
\endaligned
\end{equation*}

Last, the paper is organized as follows. In Section 2, we introduce
Strichartz estimates and perturbation theory in $\mathbb{R}^{1+n}$;
In Section 3, we overview the proof of main theorem; In Section 4,
we show that the frequency delocalization at one time implies
spacetime bound, which means that the frequency localization of
minimal energy blow up solutions; In Section $5$, we show that
$L^{2n/(n-2)}_x$-norm of minimal energy blow up solutions is bounded
from below, which means that the spatial concentration of minimal
energy blow up solutions; In Section 6, we establish the
frequency-localized interaction Morawetz estimate of minimal energy
blow up solutions, which is used to eliminate soliton-like
solutions; Finally, we prevent energy evacuation of minimal energy
blow up solutions in Section 7, which is used to exclude the finite
time blow-up solutions and double low-to-high frequency cascade
solutions.

{\center
\section{Strichartz estimates and perturbation theory in
$\mathbb{R}^{1+n}$}}\label{strandpert} \setcounter{equation}{0} In
this section, we recall Strichartz estimates and the classical long
time perturbation in $\mathbb{R} \times \mathbb{R}^n$ for $n\geq 5$,
and give a modified long time perturbation, which is important to
establish the frequency localization of minimal energy blow up
solutions.

The following {\it Strichartz } inequalities are tied up with the
local well-posedness theory.

\begin{lemma}[\cite{Ca03}, \cite{CKSTT07}, \cite{KeT98}, \cite{Stri77}]\label{lemstri1}
Let $I$ be a compact time interval, and Let $u: I \times
\mathbb{R}^n \rightarrow \mathbb{C}$ be a Schwartz solution to the
forced  Schr\"{o}dinger equation
\begin{equation*}
\aligned iu_t + \Delta u = \sum^M_{m=1}F_m
\endaligned
\end{equation*}
for some Schwartz functions $F_1, \cdots, F_M$. Then
\begin{equation*}
\big\|u\big\|_{\dot{S}^k(I\times \mathbb{R}^n)} \lesssim
\big\|u(t_0)\big\|_{\dot{H}^k_x} + C \sum^M_{m=1} \big\|\nabla^k
F_m(u)\big\|_{L^{q'_m}_tL^{r'_m}_x(I\times \mathbb{R}^n)}
\end{equation*}
for any  $k\geq 0$ and $t_0 \in I$, and any sharp admissible pairs
$(q_1, r_1)$, $\cdots$, $(q_M, r_M)$, where we use $p'$ to denote
the dual exponent to $p$, i.e.  $1/p'+1/p=1$.
\end{lemma}
On the inhomogeneous Strichartz estimate, we also have
\begin{lemma}[Inhomogenous Strichartz estimate, \cite{Fos05}, \cite{Ka94}, \cite{Vile07}]\label{foschistrichartz}
If $v$ is the solution of
\begin{equation*}
\aligned iv_t+\Delta v = F(t,x)
\endaligned
\end{equation*}
with zero data and inhomogeneous term  $F$ supported on
$\mathbb{R}\times \mathbb{R}^d$, then we have the estimate
\begin{equation*}
\aligned \big\|v\big\|_{L^q_tL^r_x} \lesssim
\big\|F\big\|_{L^{\widetilde{q}'}_tL^{\widetilde{r}'}_x}
\endaligned
\end{equation*}
whenever $(q,r)$, $(\widetilde{q},\widetilde{r})$  are acceptable,
verify the scaling condition
\begin{equation*}
\aligned \frac{1}{q}+\frac{1}{\widetilde{q}}=\frac{n}{2}\big(
1-\frac{1}{r}-\frac{1}{\widetilde{r}} \big),
\endaligned
\end{equation*} and  either the conditions
\begin{equation*}
\aligned \frac{1}{q}+\frac{1}{\widetilde{q}}<1,\quad
\frac{n-2}{n}\leq \frac{r}{\widetilde{r}}\leq \frac{n}{n-2},
\endaligned
\end{equation*}
or the conditions
\begin{equation*}
\aligned \frac{1}{q}+\frac{1}{\widetilde{q}}=1,\quad
\frac{n-2}{n}\leq \frac{r}{\widetilde{r}} \leq \frac{n}{n-2}, \quad
\frac{1}{r} \leq \frac{1}{q}, \ \ \frac{1}{\widetilde{r}} \leq
\frac{1}{\widetilde{q}}.
\endaligned
\end{equation*}

\end{lemma}

Now, similar as in \cite{CKSTT07}, \cite{RyV05}, \cite{TaV05} and
\cite{Vi05}, we first have

\begin{lemma}[Classical long time perturbation]\label{longperturbations}
Let I be a compact interval, and let $\widetilde{u}$ be a function
on $I\times \mathbb{R}^n$ which obeys the bounds
\begin{equation*}
\big\| \widetilde{u} \big\|_{L^6_tL_x^{\frac{6n}{3n-8}}(I\times
\mathbb{R}^n)}   \leq M
\end{equation*}
and
\begin{equation*}
\big\| \widetilde{u} \big\|_{L^{\infty}_t\dot{H}^1_x(I\times
\mathbb{R}^n)} \leq E
\end{equation*}
for some $M, E>0$. Suppose also that $\widetilde{u}$ is a
near-solution to $(\ref{equ1})$ in the sense that it solves
\begin{equation}\label{near-solution}
\aligned (i\partial_t + \Delta ) \widetilde{u} =
(|x|^{-4}*|\widetilde{u}|^2)\widetilde{u}+e
\endaligned
\end{equation} for some function $e$. Let $t_0\in I$, and
let $u(t_0)$ be close to $\widetilde{u}(t_0)$ in the sense that
\begin{equation*}
\big\| u(t_0)-\widetilde{u}(t_0) \big\|_{\dot{H}^1_x(\mathbb{R}^n)}
\leq E'
\end{equation*}
for some $E'>0$. Assume also that we have the smallness conditions
\begin{eqnarray}
\Big( \sum_{N}\big\| P_N \nabla e^{i(t-t_0)\Delta}
 \big(u(t_0)-\widetilde{u}(t_0)\big)
\big\|^2_{L^6_tL_x^{\frac{6n}{3n-2}}(I\times \mathbb{R}^n)}\Big)^{\frac12} \qquad \qquad \qquad & &\nonumber \\
 + \Big( \sum_{N}\big\| P_N \nabla e^{i(t-t_0)\Delta}
\big(u(t_0)-\widetilde{u}(t_0)\big) \big\|^2_{L^3_tL^{
\frac{6n}{3n-4}}_x(I\times \mathbb{R}^n)} \Big)^{\frac12} & \leq&
\epsilon , \label{longcond4}\\
 \big\| e \big\|_{L^{\frac32}_t\dot{H}_x^{1, \frac{6n}{3n+4}}(I\times \mathbb{R}^n)}
& \leq & \epsilon \label{error1}
\end{eqnarray}
for some $0< \epsilon< \epsilon_1$, where $\epsilon_1$ is some
constant $\epsilon_1=\epsilon_1(E, E', M)>0$.

We conclude that there exists a solution $u$ to $(\ref{equ1})$ on
$I\times \mathbb{R}^n$ with the specified initial data $u(t_0)$ at
$t_0$, and furthermore
\begin{equation*}\aligned
\big\| u- \widetilde{u} \big\|_{\dot{S}^1(I\times \mathbb{R}^n)}  &
\leq C(M, E, E'), \\
\big\| u \big\|_{\dot{S}^1(I\times \mathbb{R}^n)}  & \leq C(M, E, E'), \\
\big\| u-\widetilde{u} \big\|_{L^6_tL_x^{\frac{6n}{3n-8}}}+ \big\|
u- \widetilde{u} \big\|_{L^3_t\dot{H}_x^{1, \frac{6n}{3n-4}}}  &
\leq C(M, E, E') \epsilon.
\endaligned\end{equation*}
\end{lemma}

\begin{remark}\label{ltpl}
Note that $u(t_0)-\widetilde{u}(t_0)$ is allowed to have large
energy, albeit at the cost of forcing $\epsilon$ to be smaller, and
worsening the bounds in $\big\| u- \widetilde{u}
\big\|_{\dot{S}^1(I\times \mathbb{R}^n)} $. From the Strichartz
estimate and Plancherel's theorem, we have
\begin{equation*}
\aligned \text{L.H.S of } (\ref{longcond4}) &  \lesssim \Big(
\sum_{N} \big\|P_N \nabla \big( u(t_0) -u(t_0)\big) \big\|^2_{L^2_x}
\Big)^{1/2}\\
 & \lesssim \big\|\nabla \big( u(t_0) -u(t_0)\big) \big\|_{L^2_x} \\
 & \lesssim E'.
\endaligned
\end{equation*}
Hence, the hypotheses  $(\ref{longcond4})$ are redundant if one is
willing to take $E'=O(\epsilon)$.
\end{remark}

Based on Lemma \ref{foschistrichartz}, we can also obtain the
following long time perturbation,

\begin{lemma}[Modified long time
perturbation]\label{ltp} Let I be a compact interval, and let
$\widetilde{u}$ be a function on $I\times \mathbb{R}^d$ which obeys
the bounds
\begin{equation}
\big\| \widetilde{u} \big\|_{L^6_tL_x^{\frac{6n}{3n-8}}(I\times
\mathbb{R}^n)}   \leq M
\end{equation}
and
\begin{equation}
\big\| \widetilde{u} \big\|_{L^{\infty}_t( I;\dot{H}^1_x)} \leq E
\end{equation}
for some $M, E>0$. Suppose also that $\widetilde{u}$ is a
near-solution to $(\ref{equ1})$ in the sense that it solves
$(\ref{near-solution})$ for some function $e$. Let $t_0\in I$, and
let $u(t_0)$ be close to $\widetilde{u}(t_0)$ in the sense that
\begin{equation*}
\big\| u(t_0)-\widetilde{u}(t_0) \big\|_{\dot{H}^1_x} \leq E'
\end{equation*}
for some $E'>0$. Assume also that we have the smallness conditions
\begin{eqnarray}
\big\| e^{i(t-t_0)\Delta}
 \big(u(t_0)-\widetilde{u}(t_0)\big)
\big\|_{L^3_t(I; \dot{H}^{1-\frac{2}{n},
\frac{6n^2}{3n^2-4n-12}}_x)} & \leq &
\epsilon , \label{longcondition4}\\
 \big\| e \big\|_{L^{\frac32}_t(I; \dot{H}^{1-\frac{2}{n},
\frac{6n^2}{3n^2+4n-12}}_x)} & \leq  &\epsilon \label{econdition}
\end{eqnarray}
for some $0< \epsilon< \epsilon_1$, where $\epsilon_1$ is some
constant $\epsilon_1=\epsilon_1(E, E', M)>0$.

We conclude that there exists a solution $u$ to $(\ref{equ1})$ on
$I\times \mathbb{R}^d$ with the specified initial data $u(t_0)$ at
$t_0$, and
\begin{equation*}\aligned
\big\| u \big\|_{L^6_tL_x^{\frac{6n}{3n-8}}(I\times \mathbb{R}^n)}
\leq C(M, E, E'), \quad \big\| u \big\|_{\dot{S}^1(I)}&   \leq C(M,
E,
E').\\
\big\| u-\widetilde{u} \big\|_{L^6_tL_x^{\frac{6n}{3n-8}}}+ \big\|
u- \widetilde{u} \big\|_{L^3_t(I; \dot{H}^{1-\frac{2}{n},
\frac{6n^2}{3n^2-4n-12}}_x)}  & \leq C(M, E, E') \epsilon.
\endaligned\end{equation*}
\end{lemma}

\begin{remark} As discussions in induction, checking condition $(\ref{econdition})$ is more convinient than checking condition
$(\ref{error1})$ as one deal with the interaction between the low
and high frequency. It plays an essential role to deal with the
nonlocal interaction  like Hartree equation, etc. Refer to details
in Section $4$.
\end{remark}

We end this section with a few related results. First, if a solution
cannot be continued strongly beyond a time $T_*$, then the
$L^6_tL_x^{\frac{6n}{3n-8}}$ norm must blow up near that time.

\begin{lemma}[Standard blow-up criterion,
\cite{MiXZ07a}]\label{blowupcri}
Let $u_0 \in \dot{H}^1$, and let $u$ be a strong solution to
$(\ref{equ1})$ on the slab $[t_0, T_0)\times \mathbb{R}^n$ such that
\begin{equation*}
\big\| u \big\|_{L^6_tL_x^{\frac{6n}{3n-8}}([t_0,T_0)\times
\mathbb{R}^n)}< \infty.
\end{equation*}
Then there exists $\delta>0$ such that the solution $u$ extends to a
strong solution to $(\ref{equ1})$ on the slab $[t_0,
T_0+\delta]\times \mathbb{R}^n$.
\end{lemma}

Last,  once we have $L^6_tL_x^{\frac{6n}{3n-8}}$ control of a finite
energy solution, we can control all Strichartz norms as well by the
standard argument (partion the time interval), .

\begin{lemma}[Persistence of regularity]\label{regularitypersis} Let $ s \geq 0$, I be a
compact time interval, and let $u$ be a finite energy solution to
$(\ref{equ1})$ on $I\times \mathbb{R}^n$ obeying the bounds
\begin{equation*}
\big\| u \big\|_{L^6_tL_x^{\frac{6n}{3n-8}}(I\times
\mathbb{R}^n)}\leq M.
\end{equation*}
Then, if $t_0 \in I$ and $u(t_0)\in \dot{H}^s$, we have
\begin{equation*}
\big\|u\big\|_{\dot{S}^s\big(I\times \mathbb{R}^n\big)} \leq C(M,
E(u))\big\|u(t_0)\big\|_{\dot{H}^s}.
\end{equation*}
\end{lemma}

{\center
\section{Overview of proof of global spacetime bounds}}
\label{overview} \setcounter{equation}{0}
 We now outline the
proof of Theorem \ref{main}, breaking it down into a number of
smaller Propositions, which are the same  as in \cite{CKSTT07},
\cite{Vi05}, see also \cite{RyV05} and \cite{Vi06}. On one hand,
note that the non-local interaction of Hartree equation, we have to
use the modified long time perturbation  to establish the frequency
localization of minimal energy blow up solutions, instead of the
classical long time perturbation and bilinear estimate. On the other
hand, we obtain the spatial concentration of minimal energy blow up
solution after we prove that $L^{\frac{2n}{n-2}}_x$-norm of minimal
energy blow up solutions is bounded from below, which is stronger
than the statement that the potential energy of minimal energy blow
up solution is bounded from below

\subsection{Zeroth stage: Induction on energy}
We say that a solution $u$ to $(\ref{equ1})$ is \emph{Schwartz} on a
slab $I\times \mathbb{R}^n$ if $u(t)$ is a Schwartz function for all
$t \in I$. Note that such solutions are then also smooth in time as
well as space, thanks to $(\ref{equ1})$.

The first observation is that  it suffices to do so for Schwartz
solutions in order to prove Theorem $\ref{main}$. For every energy
$E\geq 0$, we define the quantity $0\leq S(E)\leq +\infty$ by
\begin{equation*}
S(E):=\sup\big\{
\big\|u\big\|_{L^6_tL_x^{\frac{6n}{3n-8}}(I_{*}\times \mathbb{R}^n)}
\big\}
\end{equation*}
where the supreme is taken over all compact interval $I_{*}\subset
\mathbb{R}$ , and   over all Schwartz solution $u$  to
$(\ref{equ1})$ on $I_{*}\times \mathbb{R}^n$ with $E(u) \leq E$. We
shall adopt the convention that $S(E)=0$ for $E<0$.

From the local well-posedness theory, we know that $(\ref{equ1})$ is
locally wellposedness in $\dot{H}^1$. Moreover, from the global
well-posedness theory for small initial data, we see that $S(E)$ is
finite for small energy $E$. Our task is to show that
\begin{equation*}
S(E)<\infty, \ \text{for all} \ E>0
\end{equation*}
Assume that $S(E)$ is not always finite. From Lemma
$\ref{longperturbations}$, we see that the set $\{ E: S(E)< \infty
\}$ is open. Clearly it is also connected and contains 0. By our
contradiction hypothesis, there must therefore exist a critical
energy $0<E_{crit}< \infty$ such that $S(E_{crit}) =+\infty$, but
$S(E) <\infty$ for all $E< E_{crit}$. One can think of $E_{crit}$ as
the minimal energy required to create a blowup solution. From the
definition of $E_{crit}$, the local well-posedness theory, and Lemma
$\ref{regularitypersis}$, we have

\begin{lemma}[Induction on energy hypothesis]\label{inductionenergy}
Let $t_0 \in \mathbb{R}$, and let $v(t_0)$ be a Schwartz function
such that $E(v(t_0)) \leq E_{crit}- \eta$ for some $\eta> 0$. Then
there exists a Schwartz global solution $v:\mathbb{R}\times
\mathbb{R}^n \rightarrow \mathbb{C} $ to $(\ref{equ1})$ with initial
data $v(t_0)$ at time $t_0$ such that
\begin{equation*}
\big\|v\big\|_{L^6_tL_x^{\frac{6n}{3n-8}}(\mathbb{R} \times
\mathbb{R}^n)}\leq S (E_{crit}- \eta)=C(\eta).
\end{equation*} Furthermore we have
$\big\|v\big\|_{\dot{S}^1(\mathbb{R}\times \mathbb{R}^n)}\leq
C(\eta)$.
\end{lemma}

For the contradiction argument, we will use six such parameters
\begin{equation*}
1\gg \eta_0 \gg \eta_1 \gg \eta_2 \gg \eta_3 \gg \eta_4 \gg \eta_5
>0
\end{equation*}
Specifically, we will need a small parameter $0<
\eta_0=\eta_0(E_{crit})\ll 1$ depending on $E_{crit}$. Then we need
a smaller quantity $0< \eta_1=\eta_1(\eta_0, E_{crit})\ll 1$ assumed
sufficiently small depending on $E_{crit} $ and $\eta_0$. We
continue in this fashion, choosing each $0<\eta_j\ll 1$ to be
sufficiently small depending on all previous quantities $\eta_0,
\cdots, \eta_{j-1}$ and the energy $E_{crit}$, all the way down to
$\eta_5$ which is extremely small, much smaller than any quantity
depending on $E_{crit}, \eta_0, \cdots, \eta_4$ that will appear in
our argument. We will always assume implicitly that each $\eta_j$
has been chosen to be sufficiently small depending on the previous
parameters. We will often display the dependence of constants on a
parameter, e. g. $C(\eta)$ denotes a large constant depending on
$\eta$, and $c(\eta)$ will denote a small constant depending upon
$\eta$. When $\eta_1 \gg \eta_2$, we will understand $c(\eta_1) \gg
c(\eta_2)$ and $C(\eta_1)\ll C(\eta_2)$.

Since $S(E_{crit})$ is infinite, it is in particular larger than
$\frac{1}{\eta_5}$. By definition of $S$, this means that we may
find a compact interval $I_{*}\subset \mathbb{R}$ and a smooth
solution $u:I_{*} \times \mathbb{R}^n \rightarrow \mathbb{C}$ to
$(\ref{equ1})$ with $\frac{E_{crit}}{2} \leq E(u) \leq E_{crit}$ so
that $u$ is ridiculously large in the sense that
\begin{equation}\label{contra1}
\big\| u\big\|_{L^6_tL_x^{\frac{6n}{3n-8}}(I_{*}\times
\mathbb{R}^n)} \geq \frac{1}{\eta_5}.
\end{equation}
We will show that this leads to a contradiction. Although $u$ does
not actually blow up, it is still convenient to thinks of $u$ as
almost blowing up in $L^6_tL_x^{\frac{6n}{3n-8}}$ in the sense of
$(\ref{contra1})$.

\begin{definition}[Definition of the minimal energy blowup
solution]\label{definitionminienergy}
 A minimal energy blowup solution of $(\ref{equ1})$ is a Schwartz
 solution on a time interval $I_{*}$ with energy
\begin{equation}\label{minienergy}
\aligned \frac{1}{2} E_{crit} \leq E(u(t))  \leq E_{crit}
\endaligned\end{equation}
and
$L^6_tL_x^{\frac{6n}{3n-8}}$ norm enormous in sense of
$(\ref{contra1})$.
\end{definition}

We remark that both conditions $(\ref{contra1})$ and
$(\ref{minienergy})$ are invariant under the scaling
$(\ref{scaling})$. Thus applying the scaling $(\ref{scaling})$ to a
minimal energy blowup solution produces another minimal energy
blowup solution. Some proofs of the sub-proposition below will
revolve around a specific frequency $N$. Henceforth we will not
mention the $E_{crit}$ dependence of our constants explicitly, even
though all our constants will depend on $E_{crit}$. We shall need
however to keep careful track of the dependence of our argument on
$\eta_0, \cdots, \eta_5$. Broadly speaking , we will start with the
largest $\eta$, namely $\eta_0$, and slowly ``retreat'' to
increasingly smaller values of $\eta$ as the argument progresses.
However we will only retreat as far as $\eta_4$, not $\eta_5$, so
that $(\ref{contra1})$ will eventually lead to a contradiction when
we show that
\begin{equation*}
\big\| u\big\|_{L^6_tL_x^{\frac{6n}{3n-8}}(I_{*}\times
\mathbb{R}^n)} \leq C(\eta_0, \cdots, \eta_4).
\end{equation*}

Together with our assumption that we are considering a minimal
energy blowup solution $u$ as in Definition
$\ref{definitionminienergy}$, the Hardy-Littlewood-Sobolev
inequality implies the bounds on kinetic energy
\begin{equation}\label{kinetic}
\big\|u\big\|_{L^{\infty}_t\dot{H}^1_x(I_{*}\times \mathbb{R}^n)}
\sim 1,
\end{equation}
and potential energy
\begin{equation}\label{potential}
\sup_{t\in I_{*}}\iint  \frac{1}{|x-y|^{4}} |u(t,x)|^2 |u(t,y)|^2\
dxdy \lesssim  1.
\end{equation}

Having displayed our preliminary bounds on the kinetic and potential
energy, we briefly discuss the mass
\begin{equation*}
\int |u(t,x)|^2 dx,
\end{equation*}
which is another conserved quantity. From $(\ref{kinetic})$ and the
Bernstein inequality, we know that the high frequencies of $u$ have
small mass:
\begin{equation}\label{hifremassbound}
\big\|P_{>M}u\big\|_{L^2}\lesssim \frac1M \ \text{for all}\  M \in
2^{\mathbb{Z}}.
\end{equation}
Thus we will still be able to use the concept of mass in our
estimates as long as we restrict our attention to sufficiently high
frequencies.

\subsection{First stage: Frequency localization and spatial concentration}

We aim to show that a minimal energy blowup solution as Definition
$\ref{definitionminienergy}$ does not exist. Intuitively, it seems
reasonable to expect that a minimal-energy blow up solution should
be ``irreducible'' in the sense that it cannot be decoupled into two
or more components of strictly smaller energy that essentially do
not interact with each other, since one of the components must then
also blow up, contradicting the minimal-energy hypothesis. In
particular, we expect at every time that such a solution should be
localized in frequency and have the spatial concentration result,
which are inspired by those in \cite{CKSTT07}, \cite{RyV05} and
\cite{Vi05}.

\begin{proposition}[Frequency delocalization implies spacetime
bound]\label{freqdelocaimplystbound} Let $\eta > 0$, and suppose
there exists a dyadic frequency $N_{lo}>0$ and a time $t_0\in I_{*}$
such that we have the energy separation conditions
\begin{equation}\label{lowfreqbound}
\big\|P_{\leq N_{lo}} u(t_0)\big\|_{\dot{H}^1(\mathbb{R}^n)} \geq
\eta
\end{equation}
and
\begin{equation}\label{hifrebound}
\big\|P_{\geq K(\eta)N_{lo}} u(t_0)\big\|_{\dot{H}^1(\mathbb{R}^n)}
\geq \eta.
\end{equation}
If $K(\eta)$ is sufficiently large depending on $\eta$, i.e.
\begin{equation*}
K(\eta)\gg C(\eta).
\end{equation*}
Then we have
\begin{equation}\label{stbound1}
\big\| u\big\|_{L^6_tL_x^{\frac{6n}{3n-8}}(I_{*}\times
\mathbb{R}^n)} \leq C(\eta).
\end{equation}
\end{proposition}
 The basic idea is as above discussion,
the main tool we need is the modified long time perturbation, which
replace the role of the classical long time perturbation and
bilinear estimate. See details in Section \ref{fdistb}.

Clearly the conclusion of Proposition $\ref{freqdelocaimplystbound}$
is in conflict with the hypothesis $(\ref{contra1})$, and so we
should expect the solution to be localized in frequency for every
time $t$. This is indeed the case:

\begin{corollary}[Frequency localization of energy at each time]\label{frequencylocalization}
A minimal energy blowup solution of $(\ref{equ1})$ satisfies: For
every time $t\in I_{*}$, there exists a dyadic frequency $N(t) \in
2^{\mathbb{Z}}$ such that for every $\eta_4 \leq \eta \leq \eta_0$,
we have small energy at frequencies $\ll N(t)$,
\begin{equation}\label{lowerenergy}
\big\|P_{\leq c(\eta)N(t)} u(t)\big\|_{\dot{H}^1} \leq \eta,
\end{equation}
small energy at frequencies $\gg N(t)$,
\begin{equation}\label{hienergy}
\big\|P_{\geq C(\eta)N(t)} u(t)\big\|_{\dot{H}^1} \leq \eta.
\end{equation}
and large energy at frequencies $\sim N(t)$,
\begin{equation}\label{mediumenergy}
\big\|P_{c(\eta)N(t)<\cdot<C(\eta)N(t)} u(t)\big\|_{\dot{H}^1} \sim
1.
\end{equation}
Here $0< c(\eta) \ll1 \ll C(\eta) < \infty$ are quantities depending
on $\eta$.
\end{corollary}

{\bf Proof: } See Corollary $4.4$ in \cite{CKSTT07},  \cite{RyV05} ,
\cite{Vi05} and \cite{Vi06}.

Having shown that a minimal energy blowup solution must be localized
in frequency, we turn our attention to space. In physical space, we
will not need the full strength of a localization result
(Proposition $4.7$ in \cite{CKSTT07}). We will settle instead for a
weaker property concerning the spatial concentration of a minimal
energy blowup solution. To derive it, we use an idea of \cite{B99a},
\cite{CKSTT07}, \cite{RyV05} and \cite{Vi05}, and restrict our
analysis to a subinterval $I_0 \subset I_*$. We need to use both the
frequency localization result and the fact that the
$L^{\frac{2n}{n-2}}_x$-norm of a minimal energy blowup solution is
bounded away from zero in order to prove spatial concentration.

Since $u$ is Schwartz, we may divide the interval $I_*$ into three
consecutive pieces $I_*=I_{-}\cup I_{0} \cup I_{I}$ where each of
the three intervals contains a third of the
$L^6_tL^{\frac{6n}{3n-8}}_x$ norm:

\begin{equation*}
\aligned \int_I \Big| \int_{\mathbb{R}^n} \big|u(t,x)
\big|^{\frac{6n}{3n-8}} dx\Big|^{\frac{3n-8}{n}}dt =\frac13
\int_{I_*} \Big| \int_{\mathbb{R}^n} \big|u(t,x)
\big|^{\frac{6n}{3n-8}} dx\Big|^{\frac{3n-8}{n}}dt \quad \text{for}\
I=I_{-}, I_{0}, I_{+}.
\endaligned
\end{equation*}

In particular from $(\ref{contra1})$ we have
\begin{equation}\label{contra2}
\big\| u\big\|_{L^6_tL_x^{\frac{6n}{3n-8}}(I \times \mathbb{R}^n)}
\gtrsim \frac{1}{\eta_5} \ \text{for}\ I=I_{-}, I_{0}, I_{+}.
\end{equation}
Thus to contradict $(\ref{contra1})$, it suffices to obtain
$L^6_tL^{\frac{6n}{3n-8}}_x$ bounds on one of the three intervals
$I_{-}, I_{0}, I_{+}$.

It is in the middle interval $I_0$ that we can obtain physical space
concentration; this shall be done in two stages. The first step is
to ensure that the norm
$\big\|u(t)\big\|_{L^{\frac{2n}{n-2}}_{x}(\mathbb{R}^n)}$ is bounded
from below.

\begin{proposition}[$L^{\frac{2n}{n-2}}_x$-norm bounded from
below]\label{stbound}
For any minimal energy blowup solution to $(\ref{equ1})$ and all
$t\in I_0$, we have
\begin{equation}\label{potentialenergybound}
\aligned \big\| u(t)\big\|_{L^{\frac{2n}{n-2}}_{x}}\geq \eta_1.
\endaligned
\end{equation}
\end{proposition}

This is proved in Section \ref{potentialenergyboundbelow}. Using
$(\ref{potentialenergybound})$ and some simple Fourier analysis as
in \cite{CKSTT07}, \cite{RyV05} and \cite{Vi05}, we can thus
establish the following concentration result:

\begin{proposition}[Spatial concentration of energy at each
time]\label{sce} For any minimal energy blowup solution to
$(\ref{equ1})$ and for each $t\in I_0$, there exists $x(t)\in
\mathbb{R}^n$ such that
\begin{equation} \label{kc}
\aligned
 \int_{|x-x(t)|\leq C(\eta_1)/N(t)} \big| \nabla u(t,x)
\big|^2dx \gtrsim c(\eta_1)
\endaligned
\end{equation} and
\begin{equation} \label{pc}
\aligned
 \int_{|x-x(t)|\leq C(\eta_1)/N(t)} \big| u(t,x)
\big|^pdx  \gtrsim c(\eta_1) N(t)^{\frac{n-2}{2}p-n}
\endaligned
\end{equation} for all $1<p<\infty$, where the implicit
constants depend on $p$.
\end{proposition}

Similar result was obtained in \cite{MiXZ07a} in the radial case; To
summarize, the statements above tell us that any minimal energy
blowup solution to the equation $(\ref{equ1})$ must be localized in
frequency space at every time and have spatial concentration result
at every time. We are still far from done: we have not yet precluded
blowup solutions in finite time, nor have we eliminated soliton or
soliton-like solutions and double low to high frequency cascade
solutions. To achieve this we need spacetime integrability bounds on
$u$. Our main tool for this is a frequency localized interaction
Morawetz estimate.

\subsection{Second stage: Frequency localized Morawetz estimate}

From Bernstein estimate, we have
\begin{equation}
\big\|P_{c(\eta_0)N(t)< \cdot < C(\eta_0)N(t)}u(t)\big\|_{\dot{H}^1}
\leq C(\eta_0)N(t)\big\|u\big\|_{L^{\infty}L^2_x}.
\end{equation}
Comparing this with $(\ref{mediumenergy})$, we obtain the lower
bound
\begin{equation*}
N(t)\geq c(\eta_0)\big\| u\big\|^{-1}_{L^{\infty}L^2_x} \quad
\text{for}\ t\in I_0.
\end{equation*}
Similar analysis as in \cite{Vi05}, we know that the quantity
\begin{equation*}
N_{min}:=\inf_{t\in I_0} N(t)
\end{equation*}
is strictly positive.

From $(\ref{lowerenergy})$ we see that the low frequency portion of
the solution has small energy; one can use Strichartz estimates to
obtain some spacetime control on those low frequencies. However, we
do not yet have much control on the high frequencies, apart from the
energy bounds $(\ref{kinetic})$.

Our initial spacetime bound in the high frequencies is provided by
the following interaction Morawetz estimate.

\begin{proposition}[Frequency-localized interaction Morawetz
estimate]\label{flim} Assuming $u$ is a minimal energy blowup
solution of $(\ref{equ1})$, and $N_{*}< c(\eta_2)N_{min}$. Then we
have
\begin{equation}\label{flime}
\aligned
  \iiint_{I_0\times  \mathbb{R}^n \times \mathbb{R}^n} \frac{\big| u_{\geq N_*}(t,y) \big|^2\big| u_{\geq N_*}(t,x) \big|^2 }{|x-y|^3} \ dxdydt
\lesssim \eta_1N^{-3}_*.
\endaligned
\end{equation}
\end{proposition}

This proposition is proven in Section \ref{imi}. It is based on the
interaction Morawetz inequality developed in \cite{CKSTT04},
\cite{CKSTT07}, \cite{RyV05} and \cite{Vi05}. The key thing about
this estimate is that the right-hand side  does not depend on $I_0$,
thus it is useful in eliminating soliton or pseudosolitons, at least
for frequencies close to $N_{min}$.

Moreover, we also obtain Proposition $\ref{imsc}$ during the proof
of Proposition $\ref{flim}$, which is Strichartz control on low and
high frequencies of the minimal energy blowup solution. By meaning
of scaling and Proposition $\ref{imsc}$, we obtain the following:
\begin{corollary}\label{himora}
Let $n\geq 5$, $u$ a minimal energy blowup solution to
$(\ref{equ1})$, and $N_{*}< c(\eta_2)N_{min}$. Then, we have
\begin{equation}\label{hfc}
\aligned \big\|P_{\geq N_*} u\big\|_{L^3L^{\frac{6n}{3n-4}}(I_0
\times \mathbb{R}^n)} & \lesssim\eta^{\frac13}_1N^{-1}_{*}.
\endaligned
\end{equation}
\end{corollary}

{\bf Proof: } The claim follows interpolating between
\begin{equation*}
\nabla P_{\geq N_*} u \in L^{\infty}_t L^2_x(I_0 \times
\mathbb{R}^n)
\end{equation*}
and
\begin{equation*}
 P_{\geq N_*} u  \in L^2_tL^{\frac{2n}{n-2}}_x (I_0 \times
\mathbb{R}^n)
\end{equation*}
which comes from Proposition $\ref{imsc}$.

Combining $(\ref{hfc})$ with Proposition $\ref{sce}$, we obtain the
following integral bound on $N(t)$.

\begin{corollary}\label{intercas}
For any minimal energy blowup solutions of $(\ref{equ1})$, we have
\begin{equation*}
\aligned \int_{I_0} \big(N(t)\big)^{-1} dt \lesssim C(\eta_1,
\eta_2) N^{-3}_{min}.
\endaligned
\end{equation*}
\end{corollary}

{\bf Proof: } By $(\ref{hfc})$, we have
\begin{equation*}
\aligned \int_{I_0} \Big( \int_{\mathbb{R}^n}\big|P_{\geq N_*}
u\big|^{\frac{6n}{3n-4}} dx\Big)^{\frac{3n-4}{2n}} dt
\lesssim\eta_1N^{-3}_{*}
\endaligned
\end{equation*}
for all $N_* \leq c(\eta_2) N_{min}$. Let  $N_* = c(\eta_2)
N_{min}$, then
\begin{equation}\label{stb}
\aligned \int_{I_0} \Big( \int_{\mathbb{R}^n}\big|P_{\geq N_*}
u\big|^{\frac{6n}{3n-4}} dx\Big)^{\frac{3n-4}{2n}} dt \lesssim
C\big( \eta_1, \eta_2\big) N^{-3}_{min}.
\endaligned
\end{equation}

On the other hand, by the Bernstein estimate and the conservation of
energy, we have
\begin{equation}\label{lowfb}
\aligned \int_{|x-x(t)| \leq \frac{C(\eta_1)}{N(t)}}
\big|P_{<N_*}u(t) \big|^{\frac{6n}{3n-4}} dx & \lesssim C(\eta_1)
N(t)^{-n} \big\|  P_{<N_*}u(t) \big\|^{\frac{6n}{3n-4}}_{L^{\infty}_x} \\
& \lesssim  C(\eta_1)
N(t)^{-n} N(t)^{\frac{3n(n-2)}{3n-4}} c(\eta_2) \big\| P_{<N_*}u(t) \big\|^{\frac{6n}{3n-4}}_{L^{\frac{2n}{n-2}}_x} \\
& \lesssim  c(\eta_2) N(t)^{-\frac{2n}{3n-4}} \big\| P_{<N_*}u(t)
\big\|^{\frac{6n}{3n-4}}_{\dot{H}^1_x}   \lesssim c(\eta_2)
N(t)^{-\frac{2n}{3n-4}}.
\endaligned
\end{equation}
By Proposition $\ref{sce}$, we also have
\begin{equation}\label{alfb}
\aligned \int_{|x-x(t)| \leq \frac{C(\eta_1)}{N(t)}} \big| u(t)
\big|^{\frac{6n}{3n-4}} dx & \gtrsim c(\eta_1)
N(t)^{-\frac{2n}{3n-4}}.
\endaligned
\end{equation}
Combining $(\ref{lowfb})$, $(\ref{alfb})$ and using the triangle
inequality, we have
\begin{equation*}
 \int_{|x-x(t)| \leq
\frac{C(\eta_1)}{N(t)}} \big| P_{\geq N_*} u(t)
\big|^{\frac{6n}{3n-4}} dx \gtrsim c(\eta_1)
N(t)^{-\frac{2n}{3n-4}}.
\end{equation*}
Integrating over $I_0$ and comparing with $(\ref{stb})$, we get the
desired result.

This corollary allows us to obtain some useful bounds in the case
when $N(t)$ is bounded from above.

\begin{corollary}[Nonconcentration implies spacetime
bound]\label{nonconcen} Let $I \subseteq I_0$, and suppose there
exists an $N_{max}>0$ such that $N(t)\leq N_{max}$ for all $t\in I$.
Then for any minimal energy blowup solution of $(\ref{equ1})$, we
have
\begin{equation*}
\big\| u \big\|_{L^6_tL_x^{\frac{6n}{3n-8}}(I\times \mathbb{R}^n)}
\lesssim C(\eta_0, \eta_1, \eta_2, N_{max}/N_{min}),
\end{equation*}
and furthermore
\begin{equation*}
\big\| u \big\|_{\dot{S}^1(I\times \mathbb{R}^n)} \lesssim C(\eta_0,
\eta_1, \eta_2, N_{max}/N_{min}).
\end{equation*}
\end{corollary}

{\bf Proof: } We will prove it by stability theory. First we may use
scale invariance $(\ref{scaling})$ to rescale $N_{min}=1$. From
Corollary $\ref{intercas}$, we obtain the useful bound
\begin{equation*}
|I_0| \lesssim C(\eta_1, \eta_2, N_{max}).
\end{equation*}

Let $\delta=\delta(\eta_0, N_{max})>0$ be a small number to be
chosen later. Partition $I_0$ into $O(\frac{|I_0|}{\delta})$
subintervals $I_1, \cdots, I_J$ with $|I_j|\leq \delta$. Let $t_j\in
I_j$. Since $N(t_j)\leq N_{max}$, Corollary
$\ref{frequencylocalization}$ yields
\begin{equation*}
\big\|P_{\geq C(\eta_0)N_{max}} u(t_j)\big\|_{\dot{H}^1_x} \leq
\eta_0.
\end{equation*}

Let $\widetilde{u}(t) = e^{i(t-t_j)\Delta}P_{<C(\eta_0)N_{max}}
u(t_j)$ be the free evolution of the low and medium frequencies of
$u(t_j)$. Then we have
\begin{equation*}
\big\|u(t_j)-\widetilde{u}(t_j)\big\|_{\dot{H}^1_x} \leq \eta_0.
\end{equation*}
Moreover, by Remark $\ref{ltpl}$, we have
\begin{equation*}\aligned
\Big( \sum_{N}\big\| P_N \nabla e^{i(t-t_j)\Delta}
 \big(u(t_j)-\widetilde{u}(t_j)\big)
\big\|^2_{L^6_tL_x^{\frac{6n}{3n-2}}}\Big)^{\frac12}
 + \Big( \sum_{N}\big\| P_N \nabla e^{i(t-t_j)\Delta}
\big(u(t_j)-\widetilde{u}(t_j)\big) \big\|^2_{L^3_tL^{
\frac{6n}{3n-4}}_x} \Big)^{\frac12} & \lesssim \eta_0 .
\endaligned
\end{equation*}
By the Bernstein estimate, Sobolev embedding, and conservation of
energy, we obtain
\begin{equation*}
\aligned \big\|\widetilde{u}(t)\big\|_{L^{\frac{6n}{3n-8}}_x} &
\lesssim C(\eta_0, N_{max}) \big\|
\widetilde{u}(t)\big\|_{L^{\frac{2n}{n-2}}_x} \lesssim C(\eta_0,
N_{max}) \big\|
\widetilde{u}(t)\big\|_{\dot{H}^1_x} \\
& \lesssim C(\eta_0, N_{max}) \big\| u(t_j)\big\|_{\dot{H}^1_x}
\lesssim C(\eta_0, N_{max})
\endaligned
\end{equation*}
for all $t\in I_j$, so
\begin{equation*}
\aligned
\big\|\widetilde{u}(t)\big\|_{L^6_tL^{\frac{6n}{3n-8}}_x(I_j \times
\mathbb{R}^n)} & \lesssim C(\eta_0, N_{max})\delta^{\frac16}.
\endaligned
\end{equation*}

Similarly, we have
\begin{equation*}
\aligned \big\| \nabla \big( |\nabla|^{-(n-4)}|\widetilde{u}(t)|^2
\widetilde{u}(t) \big) \big\|_{L^{\frac{6n}{3n+4}}_x} & \lesssim
\big\|\nabla \widetilde{u} (t)\big\|_{L^{\frac{6n}{3n-4}}_x} \big\|
\widetilde{u} (t) \big\|^2_{L^{\frac{6n}{3n-8}}_x} \\
& \lesssim C(\eta_0, N_{max}) \big\|\nabla \widetilde{u}(t)
\big\|_{L^2_x} \big\| \widetilde{u} (t) \big\|^2_{L^{\frac{6n}{3n-8}}_x} \\
& \lesssim C(\eta_0, N_{max}) \big\|\widetilde{u}(t)
\big\|^3_{\dot{H}^1_x}   \lesssim C(\eta_0, N_{max}), \\
\endaligned
\end{equation*}
which shows that

\begin{equation*}
\aligned \big\| \nabla \big( |\nabla|^{-(n-4)}|\widetilde{u}(t)|^2
\widetilde{u}(t) \big)
\big\|_{L^{\frac32}L^{\frac{6n}{3n+4}}_x(I_j\times \mathbb{R}^n)} &
\lesssim C(\eta_0, N_{max})\delta^{\frac23}, \\
\endaligned
\end{equation*}

Therefore, Lemma $\ref{longperturbations}$ with
$e=-(|\nabla|^{-(n-4)} |\widetilde{u}|^2) \widetilde{u}$ implies
that
\begin{equation*}
\big\| u\big\|_{L^6_tL^{\frac{6n}{3n-8}}_x(I_j \times \mathbb{R}^n)}
\lesssim 1,
\end{equation*}
provided $\delta $ and $\eta_0$ are chosen small enough. Summing
these bounds in $j$, we obtain
\begin{equation*}
\big\| u\big\|_{L^6_tL^{\frac{6n}{3n-8}}_x(I \times \mathbb{R}^n)}
\lesssim \frac{|I|}{\delta}  \lesssim \frac{|I_0|}{\delta} \lesssim
C(\eta_0, \eta_1, \eta_2, N_{max}).
\end{equation*}
The $\dot{S}^1$ bound then follows from Lemma
$\ref{regularitypersis}$.

This above corollary gives the desired contradiction to
$(\ref{contra2})$ when $ N_{max}/N_{min}$ is bounded, i.e., $N(t)$
stays in a bounded range.

\subsection{Third stage: Nonconcentration of energy}
Now, we will make use of almost conservation law of frequency
localized mass to show that any minimal energy blowup solution
cannot concentrate energy to very high frequencies. Instead the
solution always leaves a nontrival amount of mass and energy behind
at medium frequencies. This ``littering'' of the solution will serve
to keep $N(t)$ from escaping to infinity, which is inspired by the
ideas in \cite{CKSTT07}, \cite{RyV05} and \cite{Vi05}. We will prove
in Section \ref{preventenergyevacuation}.

\begin{proposition}[Energy cannot evacuate from low
frequencies]\label{ecneflf} For any minimal energy blowup solution
of $(\ref{equ1})$, we have
\begin{equation*}
N(t) \leq C(\eta_4)N_{min}
\end{equation*}
for all $t\in I_0$
\end{proposition}

By combining Proposition $\ref{ecneflf}$ with Corollary
$\ref{nonconcen}$, we encounter a contradiction to $(\ref{contra2})$
which completes the proof of Theorem $\ref{main}$.

\section{Frequency delocalization at one time implies spacetime bound
}\label{fdistb}  \setcounter{equation}{0}Using the modified long
time perturbation,  we now prove Proposition
$\ref{freqdelocaimplystbound}$ as in \cite{CKSTT07}, \cite{RyV05}
and \cite{Vi05}. Let $0<\epsilon=\epsilon(\eta)\ll 1$ be a small
number to be chosen later. If $K(\eta)$ is sufficiently large
depending on $\epsilon$, then one can find $\epsilon^{-2}$ disjoint
intervals $[\epsilon^2N_j, \epsilon^{-2}N_j]$ contained in $[N_{lo},
K(\eta)N_{lo}]$. By $(\ref{kinetic})$ and the pigeonhole principle,
we may find an $N_j$ such that the localization of $u(t_0)$ to the
interval $[\epsilon^2N_j, \epsilon^{-2}N_j]$ has very little energy:
\begin{equation}\label{kineticlittlenear1}
\big\|P_{\epsilon^2N_j\leq \cdot \leq  \epsilon^{-2}N_j}
u(t_0)\big\|_{\dot{H}^1_x} \lesssim \epsilon.
\end{equation}
Since both the statement and conclusion of the proposition are
invariant under the scaling $(\ref{scaling})$, we normalize $N_j=1$.

Define
\begin{equation*}
u_{lo}(t_0)=P_{\leq \epsilon}u(t_0), \qquad u_{hi}(t_0)=P_{\geq
1/\epsilon}u(t_0).
\end{equation*}
We claim
that $u_{hi}$ and $u_{lo}$ have smaller energy than $u$.

\begin{lemma}
If $\epsilon$ is sufficiently small depending on $\eta$, we have
\begin{equation*}
E(u_{lo}(t_0)), E(u_{hi}(t_0))\leq E_{crit}-c\eta^C.
\end{equation*}
\end{lemma}

{\bf Proof: } Without loss of generality, we will prove this for
$u_{lo}$, the proof for $u_{hi}$ is similar. Define
\begin{equation*}
u_{hi'}(t_0)=P_{>\epsilon}u(t_0),
\end{equation*}
so that $u(t_0)= u_{lo}(t_0)+ u_{hi'}(t_0)$ and consider the
quantity
\begin{equation*}
\big| E(u(t_0)) - E(u_{lo}(t_0)) - E(u_{hi'}(t_0))\big|.
\end{equation*}
By the definition of energy, we can bound this by
\begin{equation}\label{energydiff}
\aligned
  \Big| \langle \nabla u_{lo}(t_0), & \nabla u_{hi'}(t_0) \rangle \Big|   +
\Big| \int \Big(  |\nabla|^{-(n-4)} |u(t_0)|^2 |u(t_0)|^2 \\
 &  -
|\nabla|^{-(n-4)} |u_{lo}(t_0)|^2 |u_{lo}(t_0)|^2 -|\nabla|^{-(n-4)}
|u_{hi'}(t_0)|^2 |u_{hi'}(t_0)|^2 \Big) dx \Big|.
\endaligned
\end{equation}

We first deal with the kinetic energy. By the Bernstein estimate,
$(\ref{kinetic})$ and $(\ref{kineticlittlenear1})$, We have
\begin{eqnarray}
\big\| u_{hi'}(t_0)\big\|_{L^2_x} & \lesssim & \sum_{N>\epsilon}
\big\| P_N u(t_0)\big\|_{L^2_x}   \lesssim    \sum_{\epsilon< N \leq
\epsilon^{-2}} N^{-1} \epsilon + \sum_{N>\epsilon^{-2}} N^{-1}
\lesssim   1,\ \label{ineh}
\end{eqnarray}
therefore,
\begin{eqnarray}
\Big| \langle \nabla u_{lo}(t_0), \nabla u_{hi'}(t_0) \rangle \Big|
& \lesssim & \Big| \langle \nabla P_{>\epsilon} P_{\leq \epsilon}
u(t_0), \nabla u(t_0) \rangle \Big|  \lesssim  \big\| \nabla
P_{>\epsilon} P_{\leq \epsilon} u(t_0) \big\|_{L^2_x} \big\| \nabla
u(t_0)  \big\|_{L^2_x} \nonumber\\
&  \lesssim & \big\| \xi \varphi(\xi/\epsilon)\big(1-
\varphi(\xi/\epsilon)\big) \widehat{u(t_0)}(\xi) \big\|_{L^2_x}
\lesssim  \epsilon \big\| u_{hi'}(t_0)\big\|_{L^2_x}  \lesssim \
\epsilon. \label{ki}
\end{eqnarray}

Next we deal with the potential energy part of $(\ref{energydiff})$.
By the Bernstein inequality, $(\ref{kinetic})$ and
$(\ref{kineticlittlenear1})$, We have
\begin{equation*}
\aligned
 \big\|u_{hi'}(t_0)\big\|_{L^{\frac{3n}{2n-4}}_x}
 & \lesssim \sum_{N\geq
 \epsilon}  \big\|P_Nu\big\|_{L^{\frac{3n}{2n-4}}_x}  \lesssim \sum_{N\geq
 \epsilon} N^{-1} N^{\frac{n}{2}-\frac{2n-4}{3}} \big\|\nabla P_Nu\big\|_{L^2_x} \\
 & \lesssim \sum_{
 \epsilon \leq N \leq \epsilon^{-2}} N^{-1} N^{\frac{n}{2}-\frac{2n-4}{3}} \epsilon+ \sum_{N\geq
 \epsilon^{-2}} N^{-1} N^{\frac{n}{2}-\frac{2n-4}{3}}   \lesssim \epsilon^{\frac{8-n}{6}},\\
    \big\|u_{hi'}(t_0)\big\|_{L^{\frac{n}{n-2}}_x}
 & \lesssim \sum_{N\geq
 \epsilon}  \big\|P_Nu\big\|_{L^{\frac{n}{n-2}}_x}  \lesssim \sum_{N\geq
 \epsilon} N^{-1} N^{\frac{n}{2}-n+2} \big\|\nabla P_Nu\big\|_{L^2_x} \\
 & \lesssim \sum_{
 \epsilon \leq N \leq \epsilon^{-2}} N^{-1} N^{\frac{n}{2}-n+2} \epsilon+ \sum_{N\geq
 \epsilon^{-2}} N^{-1} N^{\frac{n}{2}-n+2} \lesssim \epsilon^{2-\frac{n}2},\\
\big\|u_{lo}(t_0)\big\|_{L^{\infty}_x}  & \lesssim
\epsilon^{\frac{n-2}{2}}\big\|
u_{lo}(t_0)\big\|_{L^{\frac{2n}{n-2}}_x} \lesssim
\epsilon^{\frac{n-2}{2}} \big\| u_{lo}(t_0)\big\|_{\dot{H}^1_x}
\lesssim \epsilon^{\frac{n-2}{2}}
\\
 \big\|u_{lo}(t_0)\big\|_{L^{\frac{6n}{3n-8}}_x}
 & \lesssim \epsilon^{\frac{n-2}{2}
 -\frac{3n-8}{6}}  \big\|u_{lo}(t_0)\big\|_{L^{\frac{2n}{n-2}}_x}  \lesssim \epsilon^{\frac13},
 \endaligned
 \end{equation*}
combining the above estimates with $(\ref{ineh})$, we obtain
\begin{equation}\label{po}
\aligned  \int & \Big||\nabla|^{-(n-4)} |u(t_0)|^2 |u(t_0)|^2 -
|\nabla|^{-(n-4)} |u_{lo}(t_0)|^2 |u_{lo}(t_0)|^2 -|\nabla|^{-(n-4)}
|u_{hi'}(t_0)|^2 |u_{hi'}(t_0)|^2\Big| dx  \\
\lesssim &   \int |\nabla|^{-(n-4)} |u_{lo}(t_0)|^2
\big(|u_{lo}(t_0)|+ |u_{hi'}(t_0)|\big)|u_{hi'}(t_0)| dx  \\
 &+  \int |\nabla|^{-(n-4)}  |u_{lo}(t_0)u_{hi'}(t_0)| \big(
|u_{lo}(t_0)|^2+ |u_{hi'}(t_0)|^2 \big) dx  \\
& +  \int |\nabla|^{-(n-4)}  |u_{hi'}(t_0)|^2 \big(
|u_{lo}(t_0)|+ |u_{hi'}(t_0)| \big)|u_{lo}(t_0)| dx  \\
\lesssim &  \big\|u_{hi'}\big\|^3_{L^{\frac{3n}{2n-4}}_x}
\big\|u_{lo}\big\|_{L^{\infty}_x}
+\big\|u_{hi'}\big\|^2_{L^{\frac{n}{n-2}}_x}
\big\|u_{lo}\big\|^2_{L^{\infty}_x} + \big\|u_{hi'}\big\|_{L^2_x}
\big\|u_{lo}\big\|^3_{L^{\frac{6n}{3n-8}}_x}  \\
\lesssim & \Big( \epsilon^{\frac{8-n}{6}}\Big)^3
\epsilon^{\frac{n-2}{2}} + \Big( \epsilon^{2-\frac{n}{2}}\Big)^2
\Big( \epsilon^{\frac{n-2}{2}} \Big)^2 + \Big(
\epsilon^{\frac{1}{3}}\Big)^3
 \lesssim  \epsilon.
\endaligned\end{equation} Therefore, we have from $(\ref{ki})$ and
$(\ref{po})$
\begin{equation*}
\big| E(u(t_0)) - E(u_{lo}(t_0)) - E(u_{hi'}(t_0))\big|\lesssim
\epsilon.
\end{equation*}

Since
\begin{equation*}
E(u)\leq E_{crit},
\end{equation*}
and by hypothesis $(\ref{hifrebound})$, we have
\begin{equation*}
E(u_{hi'}(t_0)) \gtrsim
\big\|u_{hi'}(t_0)\big\|^2_{\dot{H}^1_x}\gtrsim \eta^2,
\end{equation*}
the triangle inequality implies
\begin{equation*}
E(u_{lo}(t_0)) \leq E_{crit}-c\eta^C,
\end{equation*}
provided we choose $\epsilon$ sufficiently small. Similarly, one
proves
\begin{equation*}
E(u_{hi}(t_0))\leq E_{crit}-c\eta^C.
\end{equation*}

Now, since
\begin{equation*}
\aligned E(u_{lo}(t_0)) \leq E_{crit}-c\eta^C,\quad
 E(u_{hi}(t_0))  \leq
E_{crit}-c\eta^C.
\endaligned
\end{equation*}
We can apply Lemma $\ref{inductionenergy}$, we know that there exist
Schwartz solutions $u_{lo}(t), u_{hi}(t)$ to $(\ref{equ1})$ on the
slab $I_{*}\times \mathbb{R}^n$  with initial data $u_{lo}(t_0),
u_{hi}(t_0)$ at time $t_0$, and furthermore
\begin{equation*}
\big\|u_{lo}\big\|_{\dot{S}^1(I_{*}\times \mathbb{R}^n)} \leq
C(\eta),\quad  \big\|u_{hi}\big\|_{\dot{S}^1(I_{*}\times
\mathbb{R}^n)} \leq C(\eta) .
\end{equation*}
From Lemma $\ref{regularitypersis}$, we also have
\begin{equation*}
\aligned
 \big\|u_{lo}\big\|_{\dot{S}^{1+s}(I_{*}\times \mathbb{R}^n)} &
\leq
C(\eta)\big\|u_{lo}(t_0)\big\|_{\dot{H}^{1+s}_x} \leq C(\eta)\epsilon^{s}, \ \forall \ 0\leq s \leq 1,\\
 \big\|u_{hi}\big\|_{\dot{S}^{1-k}(I_{*}\times \mathbb{R}^n)}
& \leq C(\eta)\big\|u_{hi}(t_0)\big\|_{\dot{H}^{1-k}_x} \leq
C(\eta)\epsilon^{k},\ \forall \ 0\leq k \leq 1.
\endaligned
\end{equation*}

Define
\begin{equation*}
\widetilde{u}(t):=u_{lo}(t)+u_{hi}(t).
\end{equation*}
We claim that $\widetilde{u}(t)$ is a near-solution to
$(\ref{equ1})$.

\begin{lemma}
We have
\begin{equation*}
i\widetilde{u}_t+\Delta \widetilde{u} =(
|\nabla|^{-(n-4)}|\widetilde{u}|^2 )\widetilde{u} -e
\end{equation*}
where the error $e$ obeys the bound
\begin{equation}\label{condition3}
\big\|  e \big\|_{L^{\frac32}_t(I; \dot{H}^{1-\frac{2}{n},
\frac{6n^2}{3n^2+4n-12}}_x)}  \leq C(\eta)\epsilon^{\frac{4}{n} }.
\end{equation}
\end{lemma}

{\bf Proof: }  In order to estimate the error term
\begin{equation*}
\aligned e= & |\nabla|^{-(n-4)}|\widetilde{u}|^2 \widetilde{u} -
|\nabla|^{-(n-4)}|u_{lo}|^2 u_{lo}-|\nabla|^{-(n-4)}|u_{hi}|^2
u_{hi} \\
= & |\nabla|^{-(n-4)}|u_{lo}|^2 u_{hi} + 2
|\nabla|^{-(n-4)}\text{Re}(u_{lo}\overline{u}_{hi}) (u_{lo}+u_{hi})
+|\nabla|^{-(n-4)}|u_{hi}|^2 u_{lo},
\endaligned
\end{equation*}
we obtain by the Leibniz rule
\begin{equation*}
\aligned  \big\| e \big\|&_{L^{\frac32}_t(I; \dot{H}^{1-\frac{2}{n},
\frac{6n^2}{3n^2+4n-12}}_x)} \\
 \lesssim &  \quad \ \big\||\nabla|^{-(n-4)}  ( |\nabla|^{1-\frac{2}{n}} u_{lo}
\overline{u}_{lo} ) u_{hi} \big\|_{L^{\frac32}_t(I; L^{
\frac{6n^2}{3n^2+4n-12}}_x)} + \big\| |\nabla|^{-(n-4)}|u_{lo}|^2
|\nabla|^{1-\frac{2}{n}} u_{hi} \big\|_{L^{\frac32}_t(I; L^{
\frac{6n^2}{3n^2+4n-12}}_x)} \\
& + \big\|  |\nabla|^{-(n-4)}( |\nabla|^{1-\frac{2}{n}}
u_{lo}\overline{u}_{hi}) u_{lo} \big\|_{L^{\frac32}_t(I; L^{
\frac{6n^2}{3n^2+4n-12}}_x)} + \big\| |\nabla|^{-(n-4)}( u_{lo}
|\nabla|^{1-\frac{2}{n}} \overline{u}_{hi}) u_{lo}
\big\|_{L^{\frac32}_t(I; L^{
\frac{6n^2}{3n^2+4n-12}}_x)} \\
& +  \big\| |\nabla|^{-(n-4)}( u_{lo} \overline{u}_{hi})
|\nabla|^{1-\frac{2}{n}} u_{lo} \big\|_{L^{\frac32}_t(I; L^{
\frac{6n^2}{3n^2+4n-12}}_x)} + \big\| |\nabla|^{-(n-4)}
(|\nabla|^{1-\frac{2}{n}} u_{lo}\overline{u}_{hi})u_{hi}
\big\|_{L^{\frac32}_t(I; L^{ \frac{6n^2}{3n^2+4n-12}}_x)}\\
&  + \big\| |\nabla|^{-(n-4)} (u_{lo}|\nabla|^{1-\frac{2}{n}}
\overline{u}_{hi})u_{hi}\big\|_{L^{\frac32}_t(I; L^{
\frac{6n^2}{3n^2+4n-12}}_x)}+ \big\| |\nabla|^{-(n-4)}
(u_{lo}\overline{u}_{hi})|\nabla|^{1-\frac{2}{n}}
u_{hi}\big\|_{L^{\frac32}_t(I; L^{
\frac{6n^2}{3n^2+4n-12}}_x)}   \\
& + \big\| |\nabla|^{-(n-4)}|u_{hi}|^2 |\nabla|^{1-\frac{2}{n}}
u_{lo} \big\|_{L^{\frac32}_t(I; L^{ \frac{6n^2}{3n^2+4n-12}}_x)} \ \
\ +  \big\| |\nabla|^{-(n-4)} ( |\nabla|^{1-\frac{2}{n}} u_{hi}
\overline{u}_{hi} ) u_{lo} \big\|_{L^{\frac32}_t(I; L^{
\frac{6n^2}{3n^2+4n-12}}_x)} \\
=:&  \sum^{10}_{i=1}I_i.
\endaligned
\end{equation*}

We first deal with the terms which contains
$|\nabla|^{1-\frac{2}{n}} u_{lo}$. By the Hardy-Littlewood-Sobolev
inequality and Sobobev inequality, we have
\begin{equation*}
\aligned   I_1+I_3+I_5
& \lesssim  \big\| \nabla u_{lo} \big\|_{L^6_tL^{\frac{6n}{3n-8}}_x} \big\|  u_{lo} \big\|_{L^6_tL^{\frac{6n}{3n-8}}_x} \big\|  u_{hi} \big\|_{L^3_tL^{\frac{6n}{3n-4}}_x}\\
& \lesssim  \big\|  u_{lo} \big\|_{\dot{S}^2} \big\|  u_{lo}
\big\|_{\dot{S}^1} \big\|  u_{hi} \big\|_{\dot{S}^0} \lesssim C(\eta) \epsilon^2, \\
I_6+I_9 & \lesssim  \big\|  u_{lo} \big\|_{\dot{S}^2} \big\|  u_{hi} \big\|_{\dot{S}^1} \big\|  u_{hi} \big\|_{\dot{S}^0}\\
& \lesssim C(\eta) \epsilon^2, \\
\endaligned
\end{equation*}

Next, we deal with the terms which contains
$|\nabla|^{1-\frac{2}{n}} u_{hi}$. Using the
Hardy-Littlewood-Sobolev inequality, we get
\begin{equation}\label{importine1}
\aligned I_2 + I_4 \lesssim & \big\||\nabla|^{1-\frac{2}{n}}
u_{hi}\big\|_{L^6_tL^{\frac{6n}{3n-8}}_x}
\big\|u_{lo}\big\|^2_{L^3_tL^{\frac{6n^2}{3n^2-10n-6}}_x}\\
\lesssim & \big\| u_{hi}\big\|_{\dot{S}^{1-\frac{2}{n}}}
\big\|u_{lo}\big\|^2_{\dot{S}^{1+\frac{1}{n}} }   \lesssim C(\eta)
\epsilon^{\frac{4}{n}}.
\endaligned
\end{equation}
It is the place where we must use Lemma \ref{ltp} instead of Lemma
\ref{longperturbations}, otherwise we can only obtain the
boundedness of $I_2, I_4$, no any decay !

The last three estimates follow from the Hardy-Littlewood-Sobolev
inequality and Sobobev inequality
\begin{equation*}
\aligned I_7+I_8+I_{10}& \lesssim \big\| \nabla
u_{hi} \big\|_{L^2_tL^{\frac{2n}{n-2}}_x} \big\| u_{hi} \big\|_{L^6_tL^{\frac{6n}{3n-2}}_x}  \big\| u_{lo} \big\|_{L^{\infty}_tL^{\frac{2n}{n-4}}_x}  \\
& \lesssim \big\| u_{hi} \big\|_{\dot{S}^{1}} \big\| u_{hi}
\big\|_{\dot{S}^{0}}  \big\| u_{lo} \big\|_{\dot{S}^2} \lesssim
C(\eta) \epsilon^2.
\endaligned
\end{equation*}

Next, we derive estimates on $u$ from those on $\widetilde{u}$ via
perturbation theory. More precisely, we know from
$(\ref{kineticlittlenear1})$ that
\begin{equation}
\big\| u(t_0) -\widetilde{u}(t_0) \big\|_{\dot{H}^1_x} \lesssim
\epsilon,
\end{equation}
and hence, we have
\begin{equation*}\aligned
\big\| e^{i(t-t_0)\Delta}
 \big(u(t_0)-\widetilde{u}(t_0)\big)
\big\|_{L^3_t(I; \dot{H}^{1-\frac{2}{n},
\frac{6n^2}{3n^2-4n-12}}_x)}  \lesssim  \epsilon
\endaligned\end{equation*}
By the Strichartz estimate, we also have that
\begin{equation*}
\big\|\widetilde{u}\big\|_{L^{\infty}\dot{H}^1(I_* \times
\mathbb{R}^n)}\lesssim \big\|\widetilde{u}\big\|_{\dot{S}^1(I_*
\times \mathbb{R}^n)} \lesssim\big\|u_{lo}\big\|_{\dot{S}^1(I_*
\times \mathbb{R}^n)} + \big\|u_{hi}\big\|_{\dot{S}^1(I_* \times
\mathbb{R}^n)}  \lesssim C(\eta),
\end{equation*}
and hence,
\begin{equation*}
\big\| \widetilde{u} \big\|_{L^6_tL_x^{\frac{6n}{3n-8}}(I_* \times
\mathbb{R}^n)} \lesssim \big\|\widetilde{u}\big\|_{\dot{S}^1(I_*
\times \mathbb{R}^n)} \lesssim C(\eta).
\end{equation*}
So in view of $(\ref{condition3})$, if $\epsilon$ is sufficiently
small depending on $\eta$, we can apply Lemma $\ref{ltp}$ and obtain
the desired bound $(\ref{stbound1})$.  This completes the proof of
Proposition $\ref{freqdelocaimplystbound}$.

\section{Small $L^{\frac{2n}{n-2}}_x$ norm implies spacetime bound
} \label{potentialenergyboundbelow}
\setcounter{section}{5}\setcounter{equation}{0}

 We now prove Proposition $\ref{stbound}$. Here $L^{\frac{2n}{n-2}}_x$ norm is not the potential
energy of the solution. We will argue by contradiction just as in
\cite{CKSTT07}, see also \cite{RyV05} and \cite{Vi05}.  Suppose
there exists some time $t_0 \in I_0$ such that
\begin{equation}\label{stboundcontradicition}
\big\|u(t_0)\big\|_{L^{\frac{2n}{n-2}}_x} < \eta_1.
\end{equation}

Using $(\ref{scaling})$, we scale $N(t_0)=1$. If the linear
evolution $e^{i(t-t_0)\Delta }u(t_0)$ had small
$L^6_tL_x^{\frac{6n}{3n-8}}$-norm, then by perturbation theory, the
nonlinear solution would have small
$L^6_tL_x^{\frac{6n}{3n-8}}$-norm as well. Hence, we may assume
\begin{equation*}
\big\| e^{i(t-t_0)\Delta}u(t_0)
\big\|_{L^6_tL_x^{\frac{6n}{3n-8}}(\mathbb{R}\times \mathbb{R}^n)}
\gtrsim 1.
\end{equation*}

On the other hand, Corollary $\ref{frequencylocalization}$ implies
that
\begin{equation*}
\big\|P_{lo}u(t_0)\big\|_{\dot{H}^1_x} +
\big\|P_{hi}u(t_0)\big\|_{\dot{H}^1_x} \lesssim \eta_0,
\end{equation*}
where $P_{lo}=P_{<c(\eta_0)}$ and $P_{hi}=P_{>C(\eta_0)}$. The
Strichartz estimates yield
\begin{equation*}
\big\| e^{i(t-t_0)\Delta}P_{lo}u(t_0)
\big\|_{L^6_tL_x^{\frac{6n}{3n-8}}(\mathbb{R}\times \mathbb{R}^n)} +
\big\| e^{i(t-t_0)\Delta}P_{hi}u(t_0)
\big\|_{L^6_tL_x^{\frac{6n}{3n-8}}(\mathbb{R}\times \mathbb{R}^n)}
\lesssim \eta_0.
\end{equation*}
Thus,
\begin{equation*}
\big\| e^{i(t-t_0)\Delta}P_{med}u(t_0)
\big\|_{L^6_tL_x^{\frac{6n}{3n-8}}(\mathbb{R}\times \mathbb{R}^n)}
\thickapprox 1.
\end{equation*}
where $P_{med}=1-P_{lo}-P_{hi}$. However, $P_{med}u(t_0)$ has
bounded energy and Fourier support in $c(\eta_0) \lesssim |\xi|
\lesssim C(\eta_0)$, an application of the Strichartz and Bernstein
estimate yields
\begin{equation*}
\big\| e^{i(t-t_0)\Delta}P_{med}u(t_0)
\big\|_{L^{\frac{6(n-2)}{n}}_tL^{\frac{6(n-2)}{3n-8}}_x(\mathbb{R}\times
\mathbb{R}^n)} \lesssim \big\|P_{med}u(t_0)\big\|_{L^2_x}\lesssim
C(\eta_0).
\end{equation*}
Combining these estimates with the H\"{o}lder inequality, we get
\begin{equation*}
\big\| e^{i(t-t_0)\Delta}P_{med}u(t_0)
\big\|_{L^{\infty}_{t,x}(\mathbb{R}\times \mathbb{R}^n)} \gtrsim
c(\eta_0).
\end{equation*}
In particular, there exist a time $t_1\in \mathbb{R}$ and a point
$x_1 \in \mathbb{R}^n$ so that
\begin{equation}\label{concen1}
\big| e^{i(t_1-t_0)\Delta}\big(P_{med}u(t_0)\big)(x_1) \big|\gtrsim
c(\eta_0).
\end{equation}
We may perturb $t_1$ such that $t_1 \not= t_0$ and, by time reversal
symmetry, we may take $t_1 <t_0$. Let $\delta_{x_1}$ be the Dirac
mass at $x_1$. Define $f(t_1):=P_{med}\delta_{x_1}$ and for $t>t_1$
define $f(t):=e^{i(t-t_1)\Delta}f(t_1)$. We first recall a property
about $f(t)$ as in \cite{CKSTT07}.

\begin{lemma}\label{stboundlemma}
For any $t\in \mathbb{R}$ and any $1\leq p \leq \infty$, we have
\begin{equation*}
\big\|f(t)\big\|_{L^p_x} \leq C(\eta_0) \langle t-t_1
\rangle^{\frac{n}{p}-\frac{n}{2}}.
\end{equation*}
\end{lemma}

From $(\ref{stboundcontradicition})$ and the H\"{o}lder inequality,
we have
\begin{equation*}
\big| \langle f(t_0), u(t_0)\rangle \big| \lesssim
\big\|f(t_0)\big\|_{L^{\frac{2n}{n+2}}_x}
\big\|u(t_0)\big\|_{L^{\frac{2n}{n-2}}_x} \lesssim \eta_1 C(\eta_0)
\langle t_0-t_1 \rangle.
\end{equation*}
On the other hand, by $(\ref{concen1})$, we get
\begin{equation*}
\big| \langle f(t_0), u(t_0)\rangle \big| = \big|
e^{i(t_1-t_0)\Delta}\big(P_{med}u(t_0)\big)(x_1) \big|\gtrsim
c(\eta_0).
\end{equation*}
So $\langle t_0-t_1 \rangle \gtrsim \frac{c(\eta_0)}{\eta_1}$, i.e.,
$t_1$ is far from $t_0$. In particular, the time $t_1$ of
concentration must be far from $t_0$ where the
$L^{\frac{2n}{n-2}}_x$-norm is small.

Also, from Lemma $\ref{stboundlemma}$, we have
\begin{eqnarray}
\big\| f\big\|_{L^{6}_tL^{\frac{6n}{3n-2}}_x([t_0, \infty)\times \mathbb{R}^n)} &  \lesssim & C(\eta_0)\big\|\langle \cdot -t_1\rangle^{-\frac13}\big\|_{L^6([t_0, \infty))} \lesssim C(\eta_0) \eta^{\frac16}_1, \label{stle1}\\
\big\| f\big\|_{L^3_tL^{\frac{6n}{3n-4}}_x([t_0, \infty)\times
\mathbb{R}^n)} & \lesssim & C(\eta_0)\big\|\langle \cdot
-t_1\rangle^{-\frac23}\big\|_{L^3([t_0, \infty))}  \lesssim
C(\eta_0) \eta^{\frac13}_1. \label{stle2}
\end{eqnarray}

Now we use the induction hypothesis. Split $u(t_0)=v(t_0)+w(t_0)$
where $w(t_0)=\delta e^{i\theta} \Delta^{-1} f(t_0)$ for some small
$\delta=\delta(\eta_0)>0$ and phase $\theta$ to be chosen later. One
should think of $w(t_0)$ as the contribution coming from the point
$(t_1, x_1)$ where the solution concentrates. We will show that for
an appropriate choice of $\delta$ and $\theta$, $v(t_0)$ has
slightly smaller energy than $u$. By the definition of $f$ and an
integration by parts, we have
\begin{equation*}
\aligned \frac12 \int_{\mathbb{R}^n} \big| \nabla v(t_0) \big|^2 dx
\leq E_{crit} + \delta\text{Re}\ e^{-i\theta}\langle u(t_0),
f(t_0)\rangle + O(\delta^2C(\eta_0)).
\endaligned
\end{equation*}
Choosing  $\delta$ and $\theta$ appropriately, we get
\begin{equation*}
\frac12 \int_{\mathbb{R}^n} \big| \nabla v(t_0) \big|^2 dx \leq
E_{crit}-c(\eta_0).
\end{equation*}
Also, by Lemma $\ref{stboundlemma}$, we have
\begin{equation*}
\aligned \big\|w(t_0)\big\|_{L^{\frac{2n}{n-2}}_x} \lesssim
C(\eta_0) \eta_1.
\endaligned
\end{equation*}
So, by $(\ref{stboundcontradicition})$, the Hardy-Littlewood-Sobolev
inequality, and the triangle inequality we obtain
\begin{equation*}
\aligned  \iint_{\mathbb{R}^{n}\times \mathbb{R}^{n}}
\frac{|v(t_0,x)|^2 |v(t_0,y)|^2}{|x-y|^{4}}  dxdy & \lesssim \big\|
v(t_0) \big\|^4_{L^{\frac{2n}{n-2}}_x} \lesssim \big\| u(t_0)
\big\|^4_{L^{\frac{2n}{n-2}}_x} + \big\| w(t_0)
\big\|^4_{L^{\frac{2n}{n-2}}_x} \lesssim C(\eta_0) \eta^4_1.
\endaligned
\end{equation*}
Combining the above two energy estimates and taking $\eta_1$
sufficiently small depending on $\eta_0$, we obtain
\begin{equation*}
E(v(t_0)) \leq E_{crit}-c(\eta_0).
\end{equation*}
Lemma $\ref{inductionenergy}$ implies that there exists a global
solution $v$ to $(\ref{equ1})$ with initial data $v(t_0)$ at time
$t_0$ satisfying
\begin{equation*}
\big\| v \big\|_{\dot{S}^1(\mathbb{R}\times \mathbb{R}^n)} \lesssim
C(\eta_0).
\end{equation*}
In particular,
\begin{equation*}
\big\| v \big\|_{L^{\infty}_t\dot{H}^1_x([t_0, \infty)\times
\mathbb{R}^n)} + \big\| v \big\|_{L^6_tL^{\frac{6n}{3n-8}}_x([t_0,
\infty)\times \mathbb{R}^n)} \lesssim C(\eta_0).
\end{equation*}

Moreover, by the Bernstein estimate,
\begin{equation*}
\big\| w(t_0)\big\|_{\dot{H}^1_x} \lesssim C(\eta_0).
\end{equation*}
By $(\ref{stle1})$, $(\ref{stle2})$ and the frequency localization,
we estimate
\begin{equation*}\aligned
 \sum_{N}\big\| P_N \nabla e^{i(t-t_0)\Delta}
 w(t_0)
\big\|^2_{L^6_tL_x^{\frac{6n}{3n-2}}([t_0, \infty)\times
\mathbb{R}^n)}  \lesssim  & \quad \   \sum_{N \leq C(\eta_0)}\big\|
P_N \nabla e^{i(t-t_0)\Delta}
 w(t_0)
\big\|^2_{L^6_tL_x^{\frac{6n}{3n-2}}([t_0, \infty)\times
\mathbb{R}^n)} \\
& + \sum_{N> C(\eta_0)}\big\| P_N \nabla e^{i(t-t_0)\Delta}
 w(t_0)
\big\|^2_{L^6_tL_x^{\frac{6n}{3n-2}}([t_0, \infty)\times
\mathbb{R}^n)} \\
\lesssim &   \quad \  \delta(\eta_0)\sum_{N \leq C(\eta_0)} N^2 C(\eta_0)  \big\| f\big\|^2_{L^{6}_tL^{\frac{6n}{3n-2}}_x([t_0, \infty)\times \mathbb{R}^n)} \\
& + \delta(\eta_0)  \sum_{N> C(\eta_0)} N^{-2} \big\| f\big\|^2_{L^{6}_tL^{\frac{6n}{3n-2}}_x([t_0, \infty)\times \mathbb{R}^n)} \\
\lesssim & C(\eta_0) \eta^{\frac13}_1, \endaligned
\end{equation*}
\begin{equation*}\aligned
 \sum_{N}\big\| P_N  \nabla e^{i(t-t_0)\Delta}
w(t_0) \big\|^2_{L^3_tL^{ \frac{6n}{3n-4}}_x([t_0, \infty)\times
\mathbb{R}^n)} \lesssim & C(\eta_0) \eta^{\frac23}_1.
\endaligned\end{equation*}
So, if $\eta_1$ is sufficiently small depending on $\eta_0$, we can
apply Lemma $\ref{longperturbations}$ with $\widetilde{u}=v$ and
$e=0$ to conclude that $u$ extends to all of $[t_0, \infty)$ and
obeys
\begin{equation*}
\big\|u\big\|_{L^6_tL^{\frac{6n}{3n-8}}_x([t_0, \infty)\times
\mathbb{R}^n)} \lesssim C(\eta_0, \eta_1).
\end{equation*}
Since $[t_0, \infty)$ contains $I_{+}$, the above estimate
contradicts $(\ref{contra2})$ if $\eta_5$ is chosen sufficiently
small. This concludes the proof of Proposition $\ref{stbound}$.

\section{Interaction Morawetz inequality }\label{imi}
\setcounter{section}{6}\setcounter{equation}{0} The goal of this
section is to prove Proposition $\ref{flim}$, which is used to
eliminate the soliton-like solutions.

\subsection{Interaction Morawetz: Generalities} We start by recalling the standard Morawetz action centered at a
point. Let $a$ be a function on the slab $I \times \mathbb{R}^n$ and
$\phi$ satisfying
\begin{equation}\label{gee}
i\partial_t \phi + \Delta \phi =\mathcal{N}
\end{equation}
on $I \times \mathbb{R}^n$. We define the Morawetz action centered
at zero to be
\begin{equation*}
\aligned M^0_a(t)=2\int_{\mathbb{R}^n} a_j(x)\text{Im}\big(
\overline{\phi}(x)\phi_j(x) \big) dx
\endaligned
\end{equation*}
where repeated indices are implicitly summed. A simple calculation
yields
\begin{lemma}
\begin{equation*}
\aligned
\partial_t M^0_a = \int_{\mathbb{R}^n} (-\Delta \Delta a) |\phi|^2 + 4
\int_{\mathbb{R}^n}a_{jk} \text{Re} \big(\overline{\phi}_j
\phi_k\big) + 2 \int_{\mathbb{R}^n} a_j \big\{ \mathcal{N},
\phi\big\}^j_p ,
\endaligned
\end{equation*}
where we define the momentum bracket to be $\big\{ f, g \big\}_{p} =
\text{Re}( f\nabla \overline{g} -g \nabla \overline{f})$.
\end{lemma}

Now let $a(x)=|x|$, easy computations show that in dimension $n\geq
5$ we have the following identities:
\begin{equation*}
\aligned a_j(x) & = \frac{x_j}{|x|}, \quad  a_{jk}(x)   =
\frac{\delta_{jk}}{|x|}- \frac{x_j x_k}{|x|^3},\\
\Delta a(x) & = \frac{n-1}{|x|}, \quad -\Delta \Delta a(x)  =
\frac{(n-1)(n-3)}{|x|^3},
\endaligned
\end{equation*}
and hence,
\begin{equation*}
\aligned
\partial_t M^0_a & = (n-1)(n-3)\int_{\mathbb{R}^n}
\frac{|\phi(x)|^2}{|x|^3}dx + 4 \int_{\mathbb{R}^n}\Big(
\frac{\delta_{jk}}{|x|}- \frac{x_j x_k}{|x|^3}
\Big)\text{Re}\big(\overline{\phi}_j \phi_k\big)(x)dx\\
& \qquad \qquad \qquad \qquad \qquad \qquad \quad \ \  + 2
\int_{\mathbb{R}^n} \frac{x_j}{|x|} \big\{ \mathcal{N},
\phi\big\}^j_p(x)dx \\
& =(n-1)(n-3)\int_{\mathbb{R}^n} \frac{|\phi(x)|^2}{|x|^3}dx + 4
\int_{\mathbb{R}^n} \frac{1}{|x|} \big| \nabla_0 \phi(x) \big|^2 dx\\
& \qquad \qquad \qquad \qquad \qquad \qquad \quad \ \  + 2
\int_{\mathbb{R}^n} \frac{x_j}{|x|} \big\{ \mathcal{N},
\phi\big\}^j_p(x)dx \\
\endaligned
\end{equation*}
where we use $ \nabla_0$ to denote the complement of the radial
portion of the gradient.

We may center the above argument at any other point $y\in
\mathbb{R}^n$. Choosing $a(x)=|x-y|$, we define the Morawetz action
centered at $y$ to be
\begin{equation*}
\aligned M^y_a(t)=2\int_{\mathbb{R}^n} \frac{x-y}{|x-y|}
\text{Im}\big( \overline{\phi}(x) \nabla \phi(x) \big) dx.
\endaligned
\end{equation*}
The same calculations now yield that
\begin{equation*}
\aligned
\partial_t M^y_a & =(n-1)(n-3)\int_{\mathbb{R}^n} \frac{|\phi(x)|^2}{|x-y|^3}dx + 4
\int_{\mathbb{R}^n} \frac{1}{|x-y|} \big| \nabla_y \phi(x) \big|^2
dx + 2 \int_{\mathbb{R}^n} \frac{x-y}{|x-y|} \big\{ \mathcal{N},
\phi\big\}_p(x)dx .
\endaligned
\end{equation*}

We are now ready to define the interaction Morawetz potential:
\begin{equation*}
\aligned M^{interact}(t) & = \int_{\mathbb{R}^n} \big| \phi(t,y)
\big|^2 M^y_a(t)dy = 2 \text{Im} \iint_{\mathbb{R}^{n} \times
\mathbb{R}^{n}} \big|\phi(t,y) \big|^2
\frac{x-y}{|x-y|}\overline{\phi}(t,x) \nabla \phi(t,x) dxdy.
\endaligned
\end{equation*}
One gets immediately the estimate
\begin{equation*}
\aligned  \big|M^{interact}(t)\big|\leq 2
\big\|\phi(t)\big\|^3_{L^2_x} \big\|\phi(t)\big\|_{\dot{H}^1_x}.
\endaligned
\end{equation*}

Calculating the time derivative of the interaction Morawetz
potential, we get the following virial-type identity:
\begin{eqnarray}
\partial_t M^{interact} &=&
\iint_{\mathbb{R}^{n}\times \mathbb{R}^{n}} \Big((n-1)(n-3)
\frac{\big| \phi(y) \big|^2 \big|\phi(x) \big|^2}{|x-y|^3} +
2\big|\phi(y) \big|^2
\frac{x-y}{|x-y|} \big\{ \mathcal{N}, \phi \big\}_p(x) \Big) dxdy \nonumber\\
&& +\ 4 \iint_{\mathbb{R}^n\times \mathbb{R}^n}
\frac{\big|\phi(y)\big|^2 \big| \nabla_y \phi(x)\big|^2}{|x-y|} dxdy
+\ 2 \int_{\mathbb{R}^n} \partial_{y^k} \text{Im}\big( \phi \overline{\phi}_k \big)(y)M^y_a dy\label{imi1}  \\
&& +\ 4 \text{Im}\iint_{\mathbb{R}^n\times \mathbb{R}^n}
\big\{\mathcal{ N} , \phi \big\}_m(y)\frac{x-y}{|x-y|} \nabla
\phi(x) \overline{\phi}(x)dxdy. \nonumber
\end{eqnarray}
where the mass bracket is defined to be $\big\{f,g \big\}_m =
\text{Im}\big(f\overline{g} \big)$.

Similar proof as  Proposition $2.5$ in \cite{CKSTT04}, see also
Proposition $10.3$ in \cite{CKSTT07}, Lemma 5.3 in \cite{RyV05} or
Proposition 5.5 in \cite{Vi05},  we have
\begin{lemma}
$ (\ref{imi1})  \geq 0.$
\end{lemma}

Thus, integrating the virial-type identity over the compact interval
$I_0$, we get

\begin{proposition}[Interaction Morawetz inequality]\label{IMG}
Let $\phi$ be a (Schwzrtz) solution to the equation
\begin{equation*}
i\partial_t \phi + \Delta \phi =\mathcal{N}
\end{equation*}
on a spacetime slab $I_0\times \mathbb{R}^n$. Then we have
\begin{equation*}
\aligned &  \iiint_{I_0\times \mathbb{R}^{n}\times \mathbb{R}^{n}}  (n-1)(n-3)\frac{\big| \phi(t,y) \big|^2\big| \phi(t,x) \big|^2 }{|x-y|^3}  +2 \big| \phi(t,y) \big|^2 \frac{x-y }{|x-y|} \{\mathcal{N}, \phi \}_{p}(t,x) \ dxdydt \\
& \leq 4\big\|\phi \big\|^3_{L^{\infty}_tL^2_x(I_0\times
\mathbb{R}^n)} \big\|\phi \big\|_{L^{\infty}_t\dot{H}^1_x(I_0\times
\mathbb{R}^n)} + 4 \iiint_{I_0\times \mathbb{R}^{n}\times
\mathbb{R}^{n}} \big| \{ \mathcal{N} , \phi \}_{m}(t, y) \big|\
\big| \nabla\phi(t,x) \big| \ \big| \phi(t,x) \big| \ dxdydt.
\endaligned
\end{equation*}
\end{proposition}

\subsection{Interaction Morawetz: The setup}
We are ready to start the proof of Proposition $\ref{flim}$. By
scaling invariance, we normalize $N_*=1$ and define
\begin{equation*}
u_{lo}(t) =P_{\leq 1}u(t), \quad u_{hi}(t) = P_{>1}u(t).
\end{equation*}
Since we assume $1=N_*<c(\eta_2)N_{min}$, we have $1<c(\eta_2)N(t),
\forall t\in I_0$. From Corollary $\ref{frequencylocalization}$ and
the Sobolev embedding, we have the low frequency estimate
\begin{equation}\label{lowmedf}
\big\|u_{<\frac{1}{\eta_2}}\big\|_{L^{\infty}_t\dot{H}^1_x(I_0
\times \mathbb{R}^n)} +
\big\|u_{<\frac{1}{\eta_2}}\big\|_{L^{\infty}_tL^{\frac{2n}{n-2}}_x(I_0
\times \mathbb{R}^n)} \lesssim \eta_2,
\end{equation}
if $c(\eta_2)$ was chosen sufficiently small. In particular, this
implies that $u_{lo}$ has small energy
\begin{equation}\label{lowf}
\big\|u_{lo}\big\|_{L^{\infty}_t\dot{H}^1_x(I_0 \times
\mathbb{R}^n)} +
\big\|u_{lo}\big\|_{L^{\infty}_tL^{\frac{2n}{n-2}}_x(I_0 \times
\mathbb{R}^n)} \lesssim \eta_2.
\end{equation}

Using the Bernstein estimate and $(\ref{lowmedf})$, one also sees
that $u_{hi}$ has small mass
\begin{equation}\label{hif}
\aligned
 \big\|u_{hi}\big\|_{L^{\infty}_tL^2_x(I_0 \times
\mathbb{R}^n)} & \lesssim \sum_{1< N < \frac{1}{\eta_2}}
\big\|u_{N}\big\|_{L^{\infty}_tL^2_x} + \sum_{N\geq
\frac{1}{\eta_2}} \big\|u_{N}\big\|_{L^{\infty}_tL^2_x} \lesssim
\eta_2.
\endaligned
\end{equation}
Our goal is to prove $(\ref{flime})$, which is equivalent to
\begin{equation}\label{bootstrap1g}
\big\| |u_{hi} |^2 \big\|_{L^2_t \dot{H}^{-\frac{n-3}{2}}_x(I_0
\times \mathbb{R}^n)} \lesssim \eta^{\frac12}_1.
\end{equation}

By a standard continuity argument, it suffices to prove
$(\ref{bootstrap1g})$ under the bootstrap hypothesis
\begin{equation}\label{bootstrap1b}
\big\| |u_{hi} |^2 \big\|_{L^2_t \dot{H}^{-\frac{n-3}{2}}_x(I_0
\times \mathbb{R}^n)} \lesssim \big(C_0\eta_1\big)^{\frac12}.
\end{equation}
for a large constant $C_0$ depending on energy but not on any of the
$\eta$'s.

First, let us note that $(\ref{bootstrap1b})$ and Lemma $5.6$ in
\cite{Vi05} imply
\begin{equation}\label{bootstrap1bb}
\big\| |\nabla |^{-\frac{n-3}{4}} u_{hi} \big\|_{L^4_{t,x}(I_0
\times \mathbb{R}^n)} \lesssim \big(C_0\eta_1\big)^{\frac14}.
\end{equation}

We now use Proposition $\ref{IMG}$ to derive an interaction Morawetz
estimate for $\phi=u_{hi}$.

\begin{proposition}\label{imee}
With the notation and assumptions above, we have
\begin{eqnarray}
 &\displaystyle & \iiint_{I_0 \times \mathbb{R}^{n}\times \mathbb{R}^{n}}  \frac{\big| u_{hi}(t,y)
\big|^2\big| u_{hi}(t,x) \big|^2 }{|x-y|^3} -  \big| u_{hi}(t,y)
\big|^2 \frac{x-y }{|x-y|} \Big(\nabla
\big|\nabla\big|^{-(n-4)}|u_{hi}|^2 |u_{hi}|^2 \Big)(t,x) \
dxdydt \nonumber \\
 & \leq  &  \eta^3_2 \label{ine0} \\
&   &  + \eta_2\iint_{I_0\times \mathbb{R}^n} \big| u_{hi}\
P_{hi}(\big|\nabla\big|^{-(n-4)}|u|^2 u -
\big|\nabla\big|^{-(n-4)}|u_{lo}|^2 u_{lo}
 -\big|\nabla\big|^{-(n-4)}|u_{hi}|^2 u_{hi} )\big|(t,x)  dxdt \label{ine1} \\
&   &+\eta_2 \iint_{I_0\times \mathbb{R}^n}\ \big| u_{hi}\ P_{lo}(\big|\nabla\big|^{-(n-4)}|u_{hi}|^2 u_{hi})\big|(t,x) dxdt \label{ine2}\\
&  & +\eta_2 \iint_{I_0\times \mathbb{R}^n}\ \big| u_{hi}\ P_{hi}(\big|\nabla\big|^{-(n-4)}|u_{lo}|^2 u_{lo})\big|(t,x) dxdt \label{ine3}\\
&  &+ \eta^2_2 \iint_{I_0\times \mathbb{R}^n}\ \big|\nabla\big|^{-(n-4)}|u_{lo}|^2  |u_{hi}||\nabla u_{lo}|(t,x)  dxdt \label{ine4}\\
& &+  \eta^2_2  \iint_{I_0\times \mathbb{R}^n}\
\big|\nabla\big|^{-(n-4)}(|u_{lo}| \ |u_{hi}|)  \big(|u_{lo}|+ |u_{hi}|\big)|\nabla u_{lo}|(t,x)  dxdt \label{ine5}\\
& &+  \eta^2_2  \iint_{I_0\times \mathbb{R}^n}\
\big|\nabla\big|^{-(n-4)}|u_{hi}|^2 \big(|u_{lo}|+
|u_{hi}|\big)|\nabla u_{lo}|(t,x)  \  dxdt \label{ine6} \\
&  &+ \eta^2_2 \iint_{I_0\times \mathbb{R}^n}\ \Big| \nabla \big|\nabla\big|^{-(n-4)}|u_{lo}|^2 \big(  |u_{lo}|  |u_{hi}| + |u_{hi}|^2 \big)  (t,x) \Big| dxdt \label{ine16}\\
&& +  \eta^2_2 \iint_{I_0\times \mathbb{R}^n}\ \Big| \nabla \big|\nabla\big|^{-(n-4)}(|u_{lo}| \ |u_{hi}|)  \big(|u_{lo}|^2+ |u_{hi}|^2\big)(t,x) \Big|  dxdt \label{ine17}\\
& &+  \eta^2_2  \iint_{I_0\times \mathbb{R}^n}\ \Big| \nabla
\big|\nabla\big|^{-(n-4)}|u_{hi}|^2 \big(|u_{lo}|^2 + | u_{lo}||u_{hi}|\big)(t,x)  \Big|  dxdt \label{ine18} \\
&&  + \eta^2_2 \iint_{I_0\times \mathbb{R}^n}\big|\nabla P_{lo}(\big|\nabla\big|^{-(n-4)}|u|^2 u)\big|(t,x) \big|u_{hi}\big|(t,x) \ dxdt \label{ine7} \\
&  &+ \iiint_{I_0\times \mathbb{R}^n\times \mathbb{R}^n}\
\frac{\big| u_{hi}(t,y)\big|^2}{|x-y|}\ \big|\nabla\big|^{-(n-4)}|u_{lo}|^2 \big( |u_{lo}||u_{hi}|+|u_{hi}|^2\big)(t,x)  dxdydt \label{ine11}\\
& &+ \iiint_{I_0\times \mathbb{R}^n\times \mathbb{R}^n}\ \frac{\big|
u_{hi}(t,y)\big|^2}{|x-y|}\
\big|\nabla\big|^{-(n-4)}(|u_{lo}| \ |u_{hi}|)  \big(|u_{lo}|^2+|u_{lo}u_{hi}| + |u_{hi}|^2 \big)(t,x)  dxdydt \label{ine12}\\
& &+ \iiint_{I_0\times \mathbb{R}^n\times \mathbb{R}^n}\ \frac{\big|
 u_{hi}(t,y)\big|^2}{|x-y|}\
\big|\nabla\big|^{-(n-4)}|u_{hi}|^2 \big(|u_{lo}|^2+ |u_{lo}||u_{hi}|\big)(t,x) \  dxdydt \label{ine13}\\
& &+\iiint_{I_0\times \mathbb{R}^n\times \mathbb{R}^n}\ \big|
u_{hi}(t,y) \big|^2
\frac{\big|P_{lo}(\big|\nabla\big|^{-(n-4)}|u_{hi}|^2
u_{hi})\big|(t,x) \big|u_{hi}\big|(t,x) }{|x-y|} \ dxdydt
\label{ine14}
\end{eqnarray}
\end{proposition}

\begin{remark}
By symmetry, we have
\begin{equation*}
\aligned - \iiint_{I_0\times \mathbb{R}^n \times \mathbb{R}^n}\
\big| u_{hi}(t,y) \big|^2 \frac{x-y }{|x-y|} \Big(\nabla
\big|\nabla\big|^{-(n-4)}|u_{hi}|^2 |u_{hi}|^2 \Big)(t,x) \ dxdydt
\geq 0.
\endaligned
\end{equation*}
\end{remark}

{\bf Proof: } Applying Proposition \ref{IMG} with $\phi=u_{hi}$ and
$\mathcal{N}=P_{hi}(\big|\nabla\big|^{-(n-4)}|u|^2 u)$, we have
\begin{equation*}
\aligned& \iiint_{I_0\times \mathbb{R}^{n}\times \mathbb{R}^{n}}   (n-1)(n-3) \frac{\big| u_{hi}(t,y) \big|^2\big| u_{hi}(t,x) \big|^2 }{|x-y|^3} \ dxdydt \\
& \qquad +2\iiint_{I_0\times \mathbb{R}^{n}\times \mathbb{R}^{n}}   \ \big| u_{hi}(t,y) \big|^2 \frac{x-y }{|x-y|} \{P_{hi}(\big|\nabla\big|^{-(n-4)}|u|^2 u), u_{hi} \}_{p}(t,x) \ dxdydt \\
& \leq 4\big\|u_{hi} \big\|^3_{L^{\infty}L^2_x(I_0\times \mathbb{R}^n)} \big\|u_{hi} \big\|_{L^{\infty}_t\dot{H}^1_x(I_0\times \mathbb{R}^n)} \\
& \qquad + 4 \iiint_{I_0\times \mathbb{R}^{n}\times \mathbb{R}^{n}}
\big| \{P_{hi}(\big|\nabla\big|^{-(n-4)}|u|^2 u) , u_{hi} \}_{m}(t,
y) \big|  \big| \nabla u_{hi}(t,x) \big|   \big| u_{hi}(t,x) \big| \
dxdydt.
\endaligned
\end{equation*}

Observe that $(\ref{hif})$ plus the conservation of energy implies
\begin{equation*}
\big\|u_{hi} \big\|^3_{L^{\infty}L^2_x(I_0\times \mathbb{R}^n)}
\big\|u_{hi} \big\|_{L^{\infty}_t\dot{H}^1_x(I_0\times
\mathbb{R}^n)} \leq \eta^3_2,
\end{equation*}
which is the error terms $(\ref{ine0})$.

We consider the mass bracket term first. Exploiting cancelation, we
write
\begin{equation*}
\aligned   \{P_{hi}(\big|\nabla\big|^{-(n-4)}|u|^2 u) , u_{hi}
\}_{m} = & \big\{P_{hi}(\big|\nabla\big|^{-(n-4)}|u|^2 u -
\big|\nabla\big|^{-(n-4)}|u_{lo}|^2 u_{lo}
-\big|\nabla\big|^{-(n-4)}|u_{hi}|^2 u_{hi} ) , u_{hi} \big\}_{m} \\
&  -\{P_{lo}(\big|\nabla\big|^{-(n-4)}|u_{hi}|^2 u_{hi}) , u_{hi}
\}_{m}  +\{P_{hi}(\big|\nabla\big|^{-(n-4)}|u_{lo}|^2 u_{lo}) ,
u_{hi} \}_{m}.
\endaligned
\end{equation*}
Because
\begin{equation*}
 \int_{\mathbb{R}^n} \big| \nabla u_{hi}(t,x) \big| \ \big| u_{hi}(t,x) \big| \
dx \leq \big\|u_{hi} \big\|_{L^{\infty}L^2_x(I_0\times
\mathbb{R}^n)} \big\|u_{hi}
\big\|_{L^{\infty}_t\dot{H}^1_x(I_0\times \mathbb{R}^n)} \leq
\eta_2,
\end{equation*}
we can bound the contribution of the mass bracket term by the
following
\begin{equation*}
\aligned &  \iiint_{I_0\times \mathbb{R}^{n}\times \mathbb{R}^{n}}
\big| \{P_{hi}(\big|\nabla\big|^{-(n-4)}|u|^2 u) , u_{hi} \}_{m}(t,
y) \big|\ \big| \nabla u_{hi}(t,x) \big| \ \big| u_{hi}(t,x) \big| \
dxdydt\\
\leq & \eta_2   \iint_{I_0\times \mathbb{R}^{n}}  \big| u_{hi}(t,x)
\big| \  \big|P_{hi}(\big|\nabla\big|^{-(n-4)}|u|^2 u -
\big|\nabla\big|^{-(n-4)}|u_{lo}|^2 u_{lo}
-\big|\nabla\big|^{-(n-4)}|u_{hi}|^2 u_{hi} )(t,x) \big| dxdt \\
& +  \eta_2 \iint_{I_0\times \mathbb{R}^{n}} \big| u_{hi}(t,x) \big| \  \big(\big| P_{lo}(\big|\nabla\big|^{-(n-4)}|u_{hi}|^2 u_{hi})\big|+ \big| P_{hi}(\big|\nabla\big|^{-(n-4)}|u_{lo}|^2 u_{lo})\big| \big)(t,x)dxdt\\
\endaligned
\end{equation*}
which are the error terms $(\ref{ine1})$, $(\ref{ine2})$ and
$(\ref{ine3})$.

We turn now towards the momentum bracket term and write
\begin{equation*}
\aligned  \{P_{hi}(\big| &  \nabla\big|^{-(n-4)}|u|^2 u), u_{hi}
\}_{p}
\\ = & \quad \ \{\big|\nabla\big|^{-(n-4)}|u|^2 u, u \}_{p}
-\{\big|\nabla\big|^{-(n-4)}|u_{lo}|^2 u_{lo}, u_{lo} \}_{p}
\\
& - \{\big|\nabla\big|^{-(n-4)}|u|^2
u-\big|\nabla\big|^{-(n-4)}|u_{lo}|^2 u_{lo}, u_{lo} \}_{p}
 - \{P_{lo}(\big|\nabla\big|^{-(n-4)}|u|^2 u), u_{hi} \}_{p} \\
= & - \big(\nabla \big|\nabla\big|^{-(n-4)}|u|^2 |u|^2 - \nabla
\big|\nabla\big|^{-(n-4)}|u_{lo}|^2 |u_{lo}|^2\big)
\\
& + \nabla\emptyset[(\big|\nabla\big|^{-(n-4)}|u|^2
u-\big|\nabla\big|^{-(n-4)}|u_{lo}|^2 u_{lo}) u_{lo}]+
\emptyset[(\big|\nabla\big|^{-(n-4)}|u|^2
u-\big|\nabla\big|^{-(n-4)}|u_{lo}|^2 u_{lo})\nabla u_{lo}]
\\
& - \{P_{lo}(\big|\nabla\big|^{-(n-4)}|u|^2 u), u_{hi} \}_{p}
 =:
I +II+III
\endaligned
\end{equation*}
where $\emptyset(X)$ denotes an expression which is schematically of
the form $X$.

To estimate the contribution coming from $I$, we write
\begin{equation*}
\aligned &  \iiint_{I_0\times \mathbb{R}^{n}\times \mathbb{R}^{n}}\
\big| u_{hi}(t,y) \big|^2 \frac{x-y }{|x-y|} \Big(\nabla
\big|\nabla\big|^{-(n-4)}|u|^2 |u|^2 - \nabla
\big|\nabla\big|^{-(n-4)}|u_{lo}|^2 |u_{lo}|^2 \Big)(t,x) \
dxdydt \\
=&    \iiint_{I_0\times \mathbb{R}^{n}\times \mathbb{R}^{n}}\ \big|
u_{hi}(t,y) \big|^2 \frac{x-y }{|x-y|} \Big( \nabla
\big|\nabla\big|^{-(n-4)}|u_{hi}|^2 |u_{hi}|^2 \Big)(t,x) \
dxdydt \\
&  +  \iiint_{I_0\times \mathbb{R}^{n}\times \mathbb{R}^{n}}\ \big|
u_{hi}(t,y) \big|^2 \frac{x-y }{|x-y|} \nabla
\big|\nabla\big|^{-(n-4)}\Big(|u|^2 |u|^2 - |u_{lo}|^2 |u_{lo}|^2  -
|u_{hi}|^2 |u_{hi}|^2 \Big)(t,x) \ dxdydt.
\endaligned
\end{equation*}
The first term is the left-hand side term in Proposition
$\ref{imee}$. On the other hand, observing that
\begin{equation}\label{ine15}
\aligned \big\| \big|u_{hi}\big|^2\big\|_{L^{\infty}_t L^1_x(I_0
\times \mathbb{R}^n)} \lesssim \big\| u_{hi}\big\|^2_{L^{\infty}_t
L^2_x(I_0 \times \mathbb{R}^n)} \lesssim \eta^2_2,
\endaligned
\end{equation}
we take the absolute values inside the integrals and use
$(\ref{ine15})$ to obtain
\begin{equation*}
\aligned  \Big| & \iiint_{I_0\times \mathbb{R}^{n}\times
\mathbb{R}^{n}}\ \big| u_{hi}(t,y) \big|^2 \frac{x-y }{|x-y|}
\Big(\nabla \big|\nabla\big|^{-(n-4)}|u|^2 |u|^2 - \nabla
\big|\nabla\big|^{-(n-4)}|u_{lo}|^2 |u_{lo}|^2 \\
& \qquad \qquad \qquad \qquad \qquad \qquad \qquad \qquad \qquad
\qquad  \ \ - \nabla \big|\nabla\big|^{-(n-4)}|u_{hi}|^2 |u_{hi}|^2
\Big)(t,x) \ dxdydt
\Big|\\
\leq &\quad \  \eta^2_2 \iint_{I_0\times \mathbb{R}^n}\Big| \nabla \big|\nabla\big|^{-(n-4)}|u_{lo}|^2 \big(  |u_{lo}|  |u_{hi}| + |u_{hi}|^2 \big)  (t,x)  \Big| dxdt  \\
& +  \eta^2_2  \iint_{I_0\times \mathbb{R}^n}\Big| \nabla \big|\nabla\big|^{-(n-4)}(|u_{lo}| \ |u_{hi}|)  \big(|u_{lo}|^2+ |u_{hi}|^2\big)(t,x)  \Big| dxdt  \\
& +  \eta^2_2  \iint_{I_0\times \mathbb{R}^n}\Big| \nabla
\big|\nabla\big|^{-(n-4)}|u_{hi}|^2 \big(|u_{lo}|^2 + | u_{lo}||u_{hi}|\big)(t,x)  \Big| dxdt   \\
\endaligned
\end{equation*}
which are the error terms $(\ref{ine16})$, $(\ref{ine17})$ and
$(\ref{ine18})$.

To estimate the contribution of $II$. Integrating by parts for the
first term, we obtain
\begin{equation*}
\aligned  \iiint_{I_0\times \mathbb{R}^{n}\times \mathbb{R}^{n}}\
\big| u_{hi}(t,y) \big|^2 \frac{\Big| \big|\nabla\big|^{-(n-4)}|u|^2
u(t,x)-\big|\nabla\big|^{-(n-4)}|u_{lo}|^2 u_{lo}(t,x)\Big|
\big|u_{lo}(t,x)\big|}{|x-y|} \ dxdydt\\
+  \iiint_{I_0\times \mathbb{R}^{n}\times \mathbb{R}^{n}}\ \big|
u_{hi}(t,y) \big|^2 \Big| \big|\nabla\big|^{-(n-4)}|u|^2
u(t,x)-\big|\nabla\big|^{-(n-4)}|u_{lo}|^2 u_{lo}(t,x)\Big|
\big|\nabla u_{lo}(t,x)\big| \ dxdydt.
\endaligned
\end{equation*}

We estimate the error term coming from the first term
\begin{equation*}
\aligned & \quad \  \iiint_{I_0\times \mathbb{R}^{n}\times
\mathbb{R}^{n}}\ \big| u_{hi}(t,y) \big|^2 \frac{\Big|
\big|\nabla\big|^{-(n-4)}|u|^2
u(t,x)-\big|\nabla\big|^{-(n-4)}|u_{lo}|^2 u_{lo}(t,x)\Big|
\big|u_{lo}(t,x)\big|}{|x-y|} \ dxdydt\\
\leq & \quad \  \iiint_{I_0\times \mathbb{R}^{n}\times
\mathbb{R}^{n}}\
\frac{\big| u_{hi}(t,y)\big|^2}{|x-y|}\ \big|\nabla\big|^{-(n-4)}|u_{lo}|^2  |u_{lo}||u_{hi}|(t,x)  dxdydt\\
& +  \iiint_{I_0\times \mathbb{R}^{n}\times \mathbb{R}^{n}}\
\frac{\big| u_{hi}(t,y)\big|^2}{|x-y|}\
\big|\nabla\big|^{-(n-4)}(|u_{lo}| \ |u_{hi}|)  \big(|u_{lo}|^2+ |u_{lo}||u_{hi}|\big)(t,x)  dxdydt\\
& +  \iiint_{I_0\times \mathbb{R}^{n}\times \mathbb{R}^{n}}\
\frac{\big| u_{hi}(t,y)\big|^2}{|x-y|}\
\big|\nabla\big|^{-(n-4)}|u_{hi}|^2 \big(|u_{lo}|^2+ |u_{lo}||u_{hi}|\big)(t,x) \  dxdydt\\
\endaligned
\end{equation*}
which are controlled by the error terms of $(\ref{ine11})$,
$(\ref{ine12})$ and $(\ref{ine13})$.

We now turn to the contribution of the second term. We  take the
absolute values inside the integrals and use $(\ref{ine15})$ to
obtain
\begin{equation*}
\aligned \iiint_{I_0\times \mathbb{R}^{n}\times \mathbb{R}^{n}}&\
\big| u_{hi}(t,y) \big|^2 \Big| \big|\nabla\big|^{-(n-4)}|u|^2
u-\big|\nabla\big|^{-(n-4)}|u_{lo}|^2 u_{lo}\Big|
\big|\nabla u_{lo}\big| \ dxdydt\\
\leq & \quad \ \eta^2_2 \iint_{I_0\times \mathbb{R}^n}\  \Big|
\big|\nabla\big|^{-(n-4)}|u|^2
u(t,x)-\big|\nabla\big|^{-(n-4)}|u_{lo}|^2 u_{lo}(t,x)\Big|
\big|\nabla u_{lo}(t,x)\big| \ dxdt\\
\leq & \quad \ \eta^2_2 \iint_{I_0\times \mathbb{R}^n}\ \big|\nabla\big|^{-(n-4)}|u_{lo}|^2  |u_{hi}||\nabla u_{lo}|(t,x)  dxdt\\
& +  \eta^2_2  \iint_{I_0\times \mathbb{R}^n}\
\big|\nabla\big|^{-(n-4)}(|u_{lo}| \ |u_{hi}|)  \big(|u_{lo}|+ |u_{hi}|\big)|\nabla u_{lo}|(t,x)  dxdt\\
& +  \eta^2_2  \iint_{I_0\times \mathbb{R}^n}\
\big|\nabla\big|^{-(n-4)}|u_{hi}|^2 \big(|u_{lo}|+ |u_{hi}|\big)|\nabla u_{lo}|(t,x)  \  dxdt\\
\endaligned
\end{equation*}
which are the error terms of $(\ref{ine4})$, $(\ref{ine5})$ and
$(\ref{ine6})$.

We consider next the contribution of $III$ to the momentum bracket
term.

When the derivative falls on $P_{lo}(\big|\nabla\big|^{-(n-4)}|u|^2
u)$, we take the absolute values inside the integrals and use
$(\ref{ine15})$ to estimate this contribution by
\begin{equation*}
\aligned
   \iiint_{I_0\times \mathbb{R}^n \times \mathbb{R}^{n}}\ & \big| u_{hi}(t,y) \big|^2 \big|\nabla P_{lo}(\big|\nabla\big|^{-(n-4)}|u|^2 u)\big|(t,x) \big|u_{hi}\big|(t,x) \ dxdydt \\
\leq & \quad \eta^2_2\iint_{I_0\times \mathbb{R}^n}\ \big|\nabla P_{lo}(\big|\nabla\big|^{-(n-4)}|u|^2 u)\big|(t,x) \big|u_{hi}\big|(t,x) \ dxdt ,\\
 \endaligned
\end{equation*}
which is the error term $(\ref{ine7})$.

When the derivative falls on $u_{hi}$, we first integrate by parts
and then take the absolute values inside the integrals to obtain,
\begin{equation*}
\aligned & \quad  \iiint_{I_0\times \mathbb{R}^{n}\times \mathbb{R}^{n}}\ \big| u_{hi}(t,y) \big|^2 \big|\nabla P_{lo}(\big|\nabla\big|^{-(n-4)}|u|^2 u)\big|(t,x) \big|u_{hi}\big|(t,x) \ dxdydt \\
&   +   \iiint_{I_0\times \mathbb{R}^{n}\times \mathbb{R}^{n}}\ \big| u_{hi}(t,y) \big|^2 \frac{\big|P_{lo}(\big|\nabla\big|^{-(n-4)}|u|^2 u)\big|(t,x) \big|u_{hi}\big|(t,x) }{|x-y|} \ dxdydt \\
\leq & \quad   \iiint_{I_0\times \mathbb{R}^{n}\times \mathbb{R}^{n}}\ \big| u_{hi}(t,y) \big|^2 \big|\nabla P_{lo}(\big|\nabla\big|^{-(n-4)}|u|^2 u)\big|(t,x) \big|u_{hi}\big|(t,x) \ dxdydt \\
& +    \iiint_{I_0\times \mathbb{R}^{n}\times \mathbb{R}^{n}}\ \big| u_{hi}(t,y) \big|^2 \frac{\big|P_{lo}(\big|\nabla\big|^{-(n-4)}|u_{hi}|^2 u_{hi})\big|(t,x) \big|u_{hi}\big|(t,x) }{|x-y|} \ dxdydt \\
& +   \iiint_{I_0\times \mathbb{R}^{n}\times \mathbb{R}^{n}}\ \big|
u_{hi}(t,y) \big|^2 \frac{\big|P_{lo}(\big|\nabla\big|^{-(n-4)}|u|^2
u- \big|\nabla\big|^{-(n-4)}|u_{hi}|^2 u_{hi})\big|(t,x)
\big|u_{hi}\big|(t,x) }{|x-y|} \ dxdydt .
\endaligned
\end{equation*}
The first term on the right-hand side of the above inequality is
controlled by $(\ref{ine7})$. The second term is controlled by
$(\ref{ine14})$. The third term is estimated by
\begin{equation*}
\aligned
\quad \   \iiint_{I_0\times \mathbb{R}^{n}\times \mathbb{R}^{n}}\ &  \big| u_{hi}(t,y) \big|^2 \frac{\big|P_{lo}(\big|\nabla\big|^{-(n-4)}|u|^2 u- \big|\nabla\big|^{-(n-4)}|u_{hi}|^2 u_{hi})\big|(t,x) \big|u_{hi}\big|(t,x) }{|x-y|} \ dxdydt \\
\leq &  \quad \   \iiint_{I_0\times \mathbb{R}^{n}\times
\mathbb{R}^{n}}\
\frac{\big| u_{hi}(t,y)\big|^2}{|x-y|}\ \big|\nabla\big|^{-(n-4)}|u_{lo}|^2 \big( |u_{lo}||u_{hi}|+|u_{hi}|^2\big)(t,x)  dxdydt\\
& +   \iiint_{I_0\times \mathbb{R}^{n}\times \mathbb{R}^{n}}\
\frac{\big| u_{hi}(t,y)\big|^2}{|x-y|}\
\big|\nabla\big|^{-(n-4)}(|u_{lo}| \ |u_{hi}|)  \big(|u_{lo}||u_{hi}|+|u_{hi}|^2 \big)(t,x)  dxdydt\\
& +   \iiint_{I_0\times \mathbb{R}^{n}\times \mathbb{R}^{n}}\
\frac{\big| u_{hi}(t,y)\big|^2}{|x-y|}\
\big|\nabla\big|^{-(n-4)}|u_{hi}|^2 \big( |u_{lo}||u_{hi}|\big)(t,x) \  dxdydt\\
\endaligned
\end{equation*}
which are controlled by $(\ref{ine11})$, $(\ref{ine12})$ and
$(\ref{ine13})$.

\subsection{Interaction Morawetz: Strichartz control}
The purpose of this section is to obtain estimates on the low and
high-frequency parts of $u$, which we will use to bound the error
terms in Proposition $\ref{imee}$.

\begin{proposition}[Strichartz control on low and high frequencies]\label{imsc}
These exists a constant $C_1$ possibly depending on the energy, but
not on any of the $\eta$'s, such that
\begin{eqnarray}
\big\| u_{lo} \big\|_{\dot{S}^1(I_0 \times \mathbb{R}^n)} & \leq &
C_1
(C_0 \eta_1)^{\frac12}, \label{lose}\\
\big\| u_{hi}\big\|_{L^2_tL^{\frac{2n}{n-2}}_x(I_0 \times
\mathbb{R}^n)} & \leq & C_1 (C_0 \eta_1)^{\frac12}. \label{hise}
\end{eqnarray}
\end{proposition}

{\bf Proof: } To prove this Proposition, we will use a bootstrap
argument. Fix $t_0:=\inf I_0$ and let $\Omega_1$ be the set of all
times $ T\in I_0 $ such that $ (\ref{lose})$ and $ (\ref{hise})$
hold on $ [t_0, T]$.

Define also $\Omega_2$ to be the set of all times $T\in I_0$ such
that
\begin{eqnarray}
\big\| u_{lo} \big\|_{\dot{S}^1([t_0, T] \times \mathbb{R}^n)} &
\leq & 2C_1
(C_0 \eta_1)^{\frac12} \label{loseb}\\
\big\| u_{hi}\big\|_{L^2_tL^{\frac{2n}{n-2}}_x([t_0, T] \times
\mathbb{R}^n)} & \leq & 2C_1 (C_0 \eta_1)^{\frac12} \label{hiseb}
\end{eqnarray}
hold.

In order to run a bootstrap argument successfully, we need to check
four things:

\begin{enumerate}
\item[1.] First, we need to see that $t_0 \in \Omega_1$; this
follows immediately from the definition of $u_{lo}$  and $u_{hi}$ at
time $t=t_0$, provided $C_1$ is sufficiently large.
\item[2.] Next, we need $\Omega_1$ to be closed; this follows
from the definition of $\Omega_1$ and Fatou's lemma.
\item[3.] Third, we need to prove that if $T\in \Omega_1$, then there
exists a small neighborhood of $T$ contained in $\Omega_2$. This
property follows from the dominated convergence theorem and the fact
that $u_{lo}$ is not only in $\dot{S}^1([t_0, T]\times
\mathbb{R}^n)$, but also in $C^0_t\dot{H}^1_x([t_0, T]\times
\mathbb{R}^n)$ because of the smoothing effect of the free
propagator.
\item[4.] The last step one needs to check is that $\Omega_2 \subset \Omega_1$ and this is
what we will focus on for the rest of the proof of Proposition
$\ref{imsc}$. Now, we fix $T\in \Omega_2$. Throughout the rest of
the proof, all spacetime norms will be on $[t_0, T]\times
\mathbb{R}^n$.
\end{enumerate}

We first consider the low frequencies. By the Strichartz estimate,
we have
\begin{equation*}
\aligned  \big\| u_{lo}\big\|_{\dot{S}^1}   \lesssim & \quad \
\big\| \nabla u_{lo} (t_0)\big\|_{L^{\infty}_tL^2_x} \\
& + \big\|
\nabla P_{lo} \big( |\nabla|^{-(n-4)} |u_{lo}|^2u_{lo} \big)
\big\|_{L^2L^{\frac{2n}{n+2}}}  +  \big\| \nabla P_{lo} \big(
|\nabla|^{-(n-4)}
|u_{lo}|^2u_{hi} \big) \big\|_{L^2L^{\frac{2n}{n+2}}}\\
&      + \big\| \nabla P_{lo} \big( |\nabla|^{-(n-4)} \text{Re}\big(
\overline{u}_{hi}u_{lo}\big)u_{lo} \big)
\big\|_{L^2L^{\frac{2n}{n+2}}}  + \big\| \nabla P_{lo} \big(
|\nabla|^{-(n-4)} \text{Re}\big( \overline{u}_{hi}u_{lo}\big)u_{hi}
\big)
\big\|_{L^2L^{\frac{2n}{n+2}}} \\
 &    + \big\|
\nabla P_{lo} \big( |\nabla|^{-(n-4)} |u_{hi}|^2u_{lo} \big)
\big\|_{L^2L^{\frac{2n}{n+2}}}  + \big\| \nabla P_{lo} \big(
|\nabla|^{-(n-4)} |u_{hi}|^2u_{hi} \big)
\big\|_{L^{\frac43}L^{\frac{2n}{n+1}}}.
\endaligned
\end{equation*}
By $(\ref{lowf})$, we have
\begin{equation*}
\big\| \nabla u_{lo}(t_0) \big\|_{ L^2_x} \lesssim \eta_2.
\end{equation*}
By the H\"{o}lder inequaltiy, Hardy-Littlewood-Sobolev inequality
and $(\ref{lowf})$, we get
\begin{equation*}
\aligned \big\| \nabla P_{lo} \big( |\nabla|^{-(n-4)}
|u_{lo}|^2u_{lo} \big) \big\|_{L^2L^{\frac{2n}{n+2}}} & \lesssim
\big\|\nabla u_{lo}\big\|_{L^2L^{\frac{2n}{n-2}}} \big\|
u_{lo}\big\|^2_{L^{\infty}L^{\frac{2n}{n-2}}}  \lesssim
2C_1(C_0\eta_1)^{\frac12} \eta^2_2 .
\endaligned
\end{equation*}
Similarly, by the Bernstein estimate, H\"{o}lder inequaltiy,
Hardy-Littlewood-Sobolev inequality, $(\ref{lowf})$ and
$(\ref{hif})$, we get
\begin{equation*}
\aligned \big\| \nabla P_{lo} \big( |\nabla|^{-(n-4)}
|u_{lo}|^2u_{hi} \big) \big\|_{L^2L^{\frac{2n}{n+2}}} & \lesssim
\big\|P_{lo} \big( |\nabla|^{-(n-4)} |u_{lo}|^2u_{hi} \big)
\big\|_{L^2L^{\frac{2n}{n+2}}} \\
& \lesssim \big\| u_{lo}\big\|^2_{L^4L^{\frac{2n}{n-3}}} \big\|
u_{hi} \big\|_{L^{\infty}L^2},
 \lesssim (2C_1)^2C_0\eta_1 \eta_2, \\
\big\| \nabla P_{lo} \big( |\nabla|^{-(n-4)} \text{Re}\big(
\overline{u}_{hi}u_{lo}\big)u_{lo} \big)
\big\|_{L^2L^{\frac{2n}{n+2}}} & \lesssim (2C_1)^2C_0\eta_1 \eta_2 \\
\big\| \nabla P_{lo} \big( |\nabla|^{-(n-4)} |u_{hi}|^2u_{lo} \big)
\big\|_{L^2L^{\frac{2n}{n+2}}}& \lesssim \big\|   P_{lo} \big(
|\nabla|^{-(n-4)} |u_{hi}|^2u_{lo} \big)
\big\|_{L^2L^{\frac{2n}{n+2}}} \\
& \lesssim \big\|   P_{lo} \big( |\nabla|^{-(n-4)} |u_{hi}|^2u_{lo}
\big) \big\|_{L^2L^{\frac{2n}{n+4}}}\\
& \lesssim \big\| u_{lo}\big\|_{L^2L^{\frac{2n}{n-4}}} \big\| u_{hi}
\big\|^2_{L^{\infty}L^2} \lesssim 2C_1(C_0\eta_1)^{\frac12}\eta^2_2,\\
\big\| \nabla P_{lo} \big( |\nabla|^{-(n-4)} \text{Re}\big(
\overline{u}_{hi}u_{lo}\big)u_{hi} \big)
\big\|_{L^2L^{\frac{2n}{n+2}}} & \lesssim 2C_1(C_0\eta_1)^{\frac12}\eta^2_2, \\
\endaligned
\end{equation*}
and
\begin{equation*}
\aligned \big\| \nabla P_{lo}   \big( |\nabla & |^{-(n-4)}
|u_{hi}|^2u_{hi} \big) \big\|_{L^{\frac43}L^{\frac{2n}{n+1}}} \\
\lesssim & \big\| P_{lo} \big( |\nabla|^{-(n-4)} |u_{hi}|^2u_{hi}
\big) \big\|_{L^{\frac43}L^{\frac{2n}{n+1}}}  \lesssim \big\|
P_{lo} \big( |\nabla|^{-(n-4)} |u_{hi}|^2u_{hi}
\big) \big\|_{L^{\frac43}L^{\frac{2n}{n+5}}}\\
\lesssim &  \big\| u_{hi}\big\|^3_{L^4L^{\frac{2n}{n-1}}} \lesssim
\big\|   u_{hi} \big\|^{\frac32}_{L^2L^{\frac{2n}{n-2}}} \big\|
u_{hi} \big\|^{\frac32}_{L^{\infty}L^2} \lesssim
\big(2C_1(C_0\eta_1)^{\frac12}\big)^{\frac32} \eta^{\frac32}_2.
\endaligned
\end{equation*}
Combining the above estimates, we get
\begin{equation*}
\aligned \big\| u_{lo}\big\|_{\dot{S}^1}  &  \lesssim \eta_2+
2C_1(C_0\eta_1)^{\frac12} \eta^2_2 +   (2C_1)^2C_0\eta_1 \eta_2 +
\big(2C_1(C_0\eta_1)^{\frac12}\big)^{\frac32} \eta^{\frac32}_2
 \leq C_1(C_0\eta_1)^{\frac12},
\endaligned
\end{equation*}
provided we choose $\eta_2$ sufficiently small.

We turn now toward the high frequencies of $u$. By the Strichartz
estimate,
\begin{equation*}
\aligned  & \quad \ \big\|
u_{hi}\big\|_{L^2_tL^{\frac{2n}{n-2}}_x}\\
 \lesssim & \quad \ \big\|  u_{hi}(t_0)\big\|_{L^2_x} \\
& + \big\| P_{hi}( |\nabla|^{-(n-4)}|u_{lo}|^2u_{lo})
\big\|_{L^2_tL^{\frac{2n}{n+2}}_x}  + \big\| P_{hi}(
|\nabla|^{-(n-4)}|u_{lo}|^2u_{hi})
\big\|_{L^2_tL^{\frac{2n}{n+2}}_x}\\
&   +\big\| P_{hi}( |\nabla|^{-(n-4)}
\text{Re}\big(\overline{u}_{hi}u_{lo}\big)u_{lo})
\big\|_{L^2_tL^{\frac{2n}{n+2}}_x} +\big\| P_{hi}( |\nabla|^{-(n-4)}
\text{Re}\big(\overline{u}_{hi}u_{lo}\big)u_{hi})
\big\|_{L^2_tL^{\frac{2n}{n+2}}_x}\\
& +\big\|P_{hi}( |\nabla|^{-(n-4)}|u_{hi}|^2u_{lo})
\big\|_{L^2_tL^{\frac{2n}{n+2}}_x}  + \big\| P_{hi}(
|\nabla|^{-(n-4)}|u_{hi}|^2u_{hi})
\big\|_{L^2_tL^{\frac{2n}{n+2}}_x}.
\endaligned
\end{equation*}
By $(\ref{hif})$, we have
\begin{equation*}
\aligned  \big\|   u_{hi}(t_0)\big\|_{L^2_x} \lesssim \eta_2.
\endaligned
\end{equation*}
By the Bernstein estimate, H\"{o}lder inequaltiy,
Hardy-Littlewood-Sobolev inequality, $(\ref{lowf})$, $(\ref{hif})$
and $(\ref{bootstrap1bb})$, we get

\begin{equation*}
\aligned \big\| P_{hi}( |\nabla|^{-(n-4)}|u_{lo}|^2u_{lo})
\big\|_{L^2_tL^{\frac{2n}{n+2}}_x} & \lesssim \big\| \nabla P_{hi}
\big( |\nabla|^{-(n-4)}
|u_{lo}|^2u_{lo} \big) \big\|_{L^2L^{\frac{2n}{n+2}}} \\
& \lesssim
\big\|\nabla u_{lo}\big\|_{L^2L^{\frac{2n}{n-2}}} \big\|
u_{lo}\big\|^2_{L^{\infty}L^{\frac{2n}{n-2}}}   \lesssim 2C_1(C_0\eta_1)^{\frac12} \eta^2_2, \\
\big\| P_{hi}( |\nabla|^{-(n-4)}|u_{lo}|^2u_{hi})
\big\|_{L^2_tL^{\frac{2n}{n+2}}_x} & \lesssim \big\|
u_{lo}\big\|^2_{L^4L^{\frac{2n}{n-3}}} \big\| u_{hi}
\big\|_{L^{\infty}L^2}
 \lesssim (2C_1)^2C_0\eta_1 \eta_2, \\
 \big\| P_{hi}( |\nabla|^{-(n-4)}
\text{Re}\big(\overline{u}_{hi}u_{lo}\big)u_{lo})
\big\|_{L^2_tL^{\frac{2n}{n+2}}_x} & \lesssim (2C_1)^2C_0\eta_1
\eta_2, \\
\big\| P_{hi}( |\nabla|^{-(n-4)}|u_{hi}|^2u_{lo})
\big\|_{L^2_tL^{\frac{2n}{n+2}}_x}& \lesssim  \big\| P_{hi}(
|\nabla|^{-(n-4)}|u_{hi}|^2u_{lo})
\big\|_{L^2_tL^{\frac{2n}{n+3}}_x}\\
& \lesssim  \big\| u_{hi}\big\|^2_{L^4L^{\frac{2n}{n-1}}} \big\| u_{lo} \big\|_{L^{\infty}L^{\frac{2n}{n-3}}}\\
& \lesssim  \big\| u_{hi}\big\|^2_{L^4L^{\frac{2n}{n-1}}} \big\| u_{lo} \big\|_{L^{\infty}L^{\frac{2n}{n-2}}}\lesssim  2C_1(C_0\eta_1)^{\frac12}  \eta^{\frac12}_2 \eta_2,\\
  \big\| P_{hi}( |\nabla|^{-(n-4)}
\text{Re}\big(\overline{u}_{hi}u_{lo}\big)u_{hi})
\big\|_{L^2_tL^{\frac{2n}{n+2}}_x} & \lesssim
2C_1(C_0\eta_1)^{\frac12}  \eta^{\frac12}_2 \eta_2,
\endaligned
\end{equation*}
and
\begin{equation*}
\aligned \big\| P_{hi}( |\nabla|^{-(n-4)}|u_{hi}|^2u_{hi})
\big\|_{L^2_tL^{\frac{2n}{n+2}}_x} & \lesssim \big\| P_{hi}(
|\nabla|^{-(n-4)}|u_{hi}|^2u_{hi})
\big\|_{L^2_tL^{\frac{2n}{n+3}}_x}  \lesssim \big\| u_{hi}\big\|^3_{L^6L^{\frac{6n}{3n-5}}} \\
& \lesssim \big\| |\nabla |^{-\frac{n-3}{4}} u_{hi}
\big\|^{\frac23}_{L^4_{t,x}}
 \big\| \nabla
u_{hi}\big\|^{\frac43}_{L^{\infty}L^2} \big\|
u_{hi} \big\|_{L^3L^{\frac{6n}{3n-4}}} \\
&   \lesssim \big(C_0\eta_1\big)^{\frac14\cdot \frac23} \big( 2C_1
(C_0 \eta_1)^{\frac12}\big)^{\frac23}= (2C_1)^{\frac23} \big( C_0
\eta_1\big)^{\frac12}.
\endaligned
\end{equation*}
Combining the above estimates, we obtain
\begin{equation*}
\aligned \big\| u_{hi}\big\|_{L^2_tL^{\frac{2n}{n-2}}_x}
 & \lesssim \eta_2 +  2C_1(C_0\eta_1)^{\frac12} \eta^2_2 + (2C_1)^2C_0\eta_1
 \eta_2 +2C_1(C_0\eta_1)^{\frac12}  \eta^{\frac12}_2 \eta_2 + (2C_1)^{\frac23} \big( C_0
\eta_1\big)^{\frac12}\\
& \leq C_1(C_0\eta_1)^{\frac12},
 \endaligned
\end{equation*}
provided $C_1$ sufficiently large and $\eta_2$ sufficiently small.

\begin{remark}\label{rem}
Interpolating between $(\ref{lowf})$ and $(\ref{lose})$, for any
Schr\"{o}dinger sharp admissible pair $(q,r)$ we obtain
\begin{equation}\label{lowfreq2}
\big\| \nabla u_{lo}\big\|_{L^qL^r(I_0 \times R^n)}\lesssim
C^{\frac{2}{q}}_1(C_0\eta_1)^{\frac1q} \eta^{1-\frac2q}_2\lesssim
(C_0\eta_1)^{\frac1q}.
\end{equation}
Similarly, interpolating between $(\ref{hif})$, $(\ref{hise})$ and
the boundness of the energy, for any Schr\"{o}dinger sharp
admissible pair $(q,r)$, we get
\begin{equation}\label{hifreq2}
\big\| u_{hi}\big\|_{L^qL^r(I_0 \times R^n)}\lesssim
C^{\frac{2}{q}}_1(C_0\eta_1)^{\frac1q} \eta^{1-\frac2q}_2\lesssim
(C_0\eta_1)^{\frac1q}.
\end{equation}
\end{remark}

\subsection{Interaction Morawetz: Error estimates}
In this section, we use the control on $u_{lo}$ and $u_{hi}$ in
Proposition $\ref{imsc}$ to bound the error terms on the right-hand
side of Proposition $\ref{imee}$. Throughout the rest of the section
all spacetime norms will be on $I_0\times R^n$.

Consider $(\ref{ine1})$, by the Bernstein estimate, H\"{o}lder
inequality, Sobolev embedding, Proposition $\ref{imsc}$, and Remark
$\ref{rem}$, we have
\begin{equation*}
\aligned
(\ref{ine1}) & \lesssim \eta_2 \big(\big\| u_{hi} |\nabla|^{-(n-4)}|u_{hi} |^2 u_{lo} \big\|_{L^1L^1} +
\big\|u_{hi} |\nabla|^{-(n-4)}|u_{lo}|^2 u_{hi} \big\|_{L^1L^1}  \\
& \quad \qquad + \big\|u_{hi}
|\nabla|^{-(n-4)}\text{Re}\big(\overline{u}_{hi} u_{lo}\big) u_{hi}
\big\|_{L^1L^1}  + \big\| u_{hi}
|\nabla|^{-(n-4)}\text{Re}\big(\overline{u}_{hi} u_{lo}\big) u_{lo}
\big\|_{L^1L^1}  \big)  \\
 & \lesssim \eta_2 \big( \big\| u_{hi} \big\|^2_{L^4L^{\frac{2n}{n-1}}} \big\| u_{hi} \big\|_{L^{\infty}L^{\frac{2n}{n-2}}} \big\| u_{lo}\big\|_{L^2L^{\frac{2n}{n-4}}}+ \big\| u_{hi} \big\|^2_{L^4L^{\frac{2n}{n-1}}} \big\| u_{lo} \big\|^2_{L^{4}L^{\frac{2n}{n-3}}}   \big)\\
 & \lesssim \eta_2 \big(  (C_0 \eta_1)^{\frac12} (C_0 \eta_1)^{\frac12}  + (C_0 \eta_1)^{\frac12} (C_0 \eta_1)^{\frac12}  \big) \lesssim \eta_2 C_0 \eta_1. \endaligned
\end{equation*}

We now move on to $(\ref{ine2})$. Using the Bernstein estimate,
Proposition $\ref{imsc}$, and Remark $\ref{rem}$, we get
\begin{equation*}
\aligned
(\ref{ine2}) & \lesssim \eta_2 \big\| u_{hi} P_{lo}\big( |\nabla|^{-(n-4)}|u_{hi} |^2 u_{hi}\big) \big\|_{L^1L^1}\\
& \lesssim \eta_2  \big\| u_{hi}\big\|_{L^3L^{\frac{6n}{3n-4}}} \big\| P_{lo}\big( |\nabla|^{-(n-4)}|u_{hi} |^2 u_{hi}\big) \big\|_{L^{\frac32} L^{\frac{6n}{3n+4}}}\\
& \lesssim \eta_2  \big\| u_{hi}\big\|_{L^3L^{\frac{6n}{3n-4}}} \big\| P_{lo}\big( |\nabla|^{-(n-4)}|u_{hi} |^2 u_{hi}\big) \big\|_{L^{\frac32} L^{\frac{6n}{3n+10}}}\\
& \lesssim \eta_2  \big\| u_{hi}\big\|^3_{L^3L^{\frac{6n}{3n-4}}}
\big\| u_{hi} \big\|_{L^{\infty}L^{\frac{2n}{n-2}}} \lesssim \eta_2
C_0 \eta_1.
\endaligned
\end{equation*}

We next estimate $(\ref{ine3}) $. By the Bernstein estimate, Sobolev
embedding, Proposition $\ref{imsc}$, and Remark $\ref{rem}$, we
estimate
\begin{equation*}
\aligned
(\ref{ine3}) & \lesssim \eta_2 \big\| u_{hi} P_{hi}\big( |\nabla|^{-(n-4)}|u_{lo} |^2 u_{lo}\big) \big\|_{L^1L^1}\\
& \lesssim \eta_2 \big\| u_{hi} \big\|_{L^4L^{\frac{2n}{n-1}}} \big\|P_{hi}\big( |\nabla|^{-(n-4)}|u_{lo} |^2 u_{lo}\big) \big\|_{L^{\frac43}L^{\frac{2n}{n+1}}}\\
& \lesssim \eta_2 \big\| u_{hi} \big\|_{L^4L^{\frac{2n}{n-1}}} \big\|\nabla P_{hi}\big( |\nabla|^{-(n-4)}|u_{lo} |^2 u_{lo}\big) \big\|_{L^{\frac43}L^{\frac{2n}{n+1}}}\\
& \lesssim \eta_2 \big\| u_{hi} \big\|_{L^4L^{\frac{2n}{n-1}}}
\big\| \nabla u_{lo} \big\|_{L^4L^{\frac{2n}{n-1}}} \big\| u_{lo}
\big\|^2_{L^4L^{\frac{2n}{n-3}}}  \lesssim \eta_2 C_0 \eta_1 .
\endaligned
\end{equation*}

We now turn toward $(\ref{ine4})- (\ref{ine6})$, and use the
H\"{o}lder, Sobolev embedding,  Proposition $\ref{imsc}$, and Remark
$\ref{rem}$ to obtain
\begin{equation*}
\aligned
(\ref{ine4}) & \lesssim \eta^2_2 \big\| |\nabla |^{-(n-4)}|u_{lo}|^2 |u_{hi}||\nabla u_{lo}| \big\|_{L^1L^1}\\
& \lesssim \eta^2_2 \big\| u_{hi} \big\|_{L^4L^{\frac{2n}{n-1}}}
\big\| \nabla u_{lo} \big\|_{L^4L^{\frac{2n}{n-1}}} \big\| u_{lo}
\big\|^2_{L^4L^{\frac{2n}{n-3}}} \lesssim \eta^2_2 C_0 \eta_1,
\\
(\ref{ine5}) & \lesssim \eta^2_2 \big\| |\nabla |^{-(n-4)}|u_{lo}u_{hi}|\big(|u_{lo}|+ |u_{hi}|\big)|\nabla u_{lo}| \big\|_{L^1L^1}\\
& \lesssim \eta^2_2 \big\| u_{hi} \big\|_{L^4L^{\frac{2n}{n-1}}}  \big\| \nabla u_{lo} \big\|_{L^4L^{\frac{2n}{n-1}}} \big\| u_{lo} \big\|^2_{L^4L^{\frac{2n}{n-3}}}\\
& \quad  + \eta^2_2 \big\| u_{hi}\big\|_{L^3L^{\frac{6n}{3n-4}}}
\big\|u_{hi}\big\|_{L^{\infty}L^{\frac{2n}{n-2}}} \big\|u_{lo}
\big\|_{L^3L^{\frac{6n}{3n-10}}} \big\|\nabla u_{lo}
\big\|_{L^3L^{\frac{6n}{3n-4}}} \lesssim \eta^2_2 C_0 \eta_1,
\\
(\ref{ine6}) & \lesssim \eta^2_2 \big\| |\nabla |^{-(n-4)} |u_{hi}|^2 \big( |u_{lo}| + |u_{hi}| \big) |\nabla u_{lo}| \big\|_{L^1L^1}\\
& \lesssim \eta^2_2 \big\| u_{hi}\big\|_{L^3L^{\frac{6n}{3n-4}}} \big\|u_{hi}\big\|_{L^{\infty}L^{\frac{2n}{n-2}}} \big\|u_{lo} \big\|_{L^3L^{\frac{6n}{3n-10}}} \big\|\nabla u_{lo} \big\|_{L^3L^{\frac{6n}{3n-4}}}\\
& \quad  + \eta^2_2 \big\|u_{hi}\big\|^2_{L^3L^{\frac{6n}{3n-4}}}
\big\| u_{hi} \big\|_{L^{\infty}L^{\frac{2n}{n-2}}} \big\|\nabla
u_{lo}\big\|_{L^3L^{\frac{6n}{3n-10}}} \lesssim \eta^2_2 C_0 \eta_1,
\endaligned
\end{equation*}
where in the last inequality we use the fact that
\begin{equation*}
\aligned
 \big\|\nabla u_{lo}\big\|_{L^3L^{\frac{6n}{3n-10}}}  \lesssim \big\|  u_{lo}\big\|_{L^3L^{\frac{6n}{3n-10}}} \lesssim  \big( C_0
 \eta_1\big)^{\frac13}.
\endaligned
\end{equation*}

 We next consider $(\ref{ine16})- (\ref{ine18})$, and use the
H\"{o}lder, Sobolev embedding,  Proposition $\ref{imsc}$, and Remark
$\ref{rem}$ to obtain
\begin{equation*}
\aligned
(\ref{ine16}) & \lesssim \eta^2_2 \big\|   \nabla \big|\nabla\big|^{-(n-4)}|u_{lo}|^2 \big(  |u_{lo}|  |u_{hi}| + |u_{hi}|^2 \big)    \big\|_{L^1L^1}\\
& \lesssim \eta^2_2 \big\| u_{hi} \big\|_{L^4L^{\frac{2n}{n-1}}}   \big\| u_{lo} \big\|^3_{L^4L^{\frac{2n}{n-3}}} + \eta^2_2 \big\| u_{hi}\big\|_{L^3L^{\frac{6n}{3n-4}}} \big\|u_{hi}\big\|_{L^{\infty}L^{\frac{2n}{n-2}}} \big\|u_{lo} \big\|^2_{L^3L^{\frac{6n}{3n-10}}}   \lesssim \eta^2_2 C_0 \eta_1, \\
(\ref{ine17}) & \lesssim \eta^2_2 \big\|   \nabla \big|\nabla\big|^{-(n-4)}(|u_{lo}| \ |u_{hi}|)  \big(|u_{lo}|^2+ |u_{hi}|^2\big)    \big\|_{L^1L^1}\\
& \lesssim \eta^2_2 \big\| u_{hi} \big\|_{L^4L^{\frac{2n}{n-1}}}   \big\| u_{lo} \big\|^3_{L^4L^{\frac{2n}{n-3}}}+ \eta^2_2 \big\|u_{hi}\big\|^2_{L^3L^{\frac{6n}{3n-4}}} \big\| u_{hi} \big\|_{L^{\infty}L^{\frac{2n}{n-2}}} \big\| u_{lo}\big\|_{L^3L^{\frac{6n}{3n-10}}}\lesssim \eta^2_2 C_0 \eta_1,\\
(\ref{ine18}) & \lesssim \eta^2_2 \big\|  \nabla
\big|\nabla\big|^{-(n-4)}|u_{hi}|^2 \big(|u_{lo}|^2 + | u_{lo}||u_{hi}|\big)    \big\|_{L^1L^1}\\
& \lesssim \eta^2_2 \big\| u_{hi}\big\|_{L^3L^{\frac{6n}{3n-4}}}
\big\|u_{hi}\big\|_{L^{\infty}L^{\frac{2n}{n-2}}} \big\|u_{lo}
\big\|^2_{L^3L^{\frac{6n}{3n-10}}} + \eta^2_2 \big\|u_{hi}\big\|^{\frac32}_{L^3L^{\frac{6n}{3n-4}}} \big\| u_{hi} \big\|^{\frac32}_{L^{\infty}L^{\frac{2n}{n-2}}} \big\|  u_{lo}\big\|_{L^2L^{\nu}}\\
& \lesssim \eta^2_2 C_0 \eta_1,
\endaligned
\end{equation*}
where $\nu = \infty$ for $n=5$ and $\nu = \frac{2n}{n-5}$ for $n\geq
6$ and  we use the fact that
\begin{equation*}
\aligned
  \big\|  u_{lo}\big\|_{L^2L^{\nu}} \lesssim   \big\|  u_{lo}\big\|_{L^2L^{\frac{2n}{n-4}} } \lesssim \big( C_0
 \eta_1\big)^{\frac12}
\endaligned
\end{equation*}
in the last inequality.

Now we turn toward $(\ref{ine7}) $. By the triangle inequality and
the similar estimates as $(\ref{ine1})$, $(\ref{ine2})$ and
$(\ref{ine3})$, we get
\begin{equation*}
\aligned
(\ref{ine7})  & \lesssim \eta^2_2 \big\| u_{hi} \nabla P_{lo} \big( |\nabla|^{-(n-4)} |u|^2 u \big) \big\|_{L^1L^1} \\
 & \lesssim \eta^2_2 \big\| u_{hi} \nabla P_{lo} \big( |\nabla|^{-(n-4)} (|u|^2 u -|u_{lo}|^2u_{lo}  -|u_{hi}|^2u_{hi}) \big) \big\|_{L^1L^1} \\
 & \quad + \eta^2_2 \big\| u_{hi} \nabla P_{lo} \big( |\nabla|^{-(n-4)} |u_{lo}|^2 u_{lo} \big) \big\|_{L^1L^1} \\
  & \quad + \eta^2_2 \big\| u_{hi} \nabla P_{lo} \big( |\nabla|^{-(n-4)} |u_{hi}|^2 u_{hi} \big) \big\|_{L^1L^1} \\
& \lesssim \eta^2_2 \big\| u_{hi} \big\|_{L^4L^{\frac{2n}{n-1}}} \big\| |\nabla|^{-(n-4)} (|u|^2 u -|u_{lo}|^2u_{lo}  -|u_{hi}|^2u_{hi})  \big\|_{L^{\frac43}L^{\frac{2n}{n+1}}} \\
& \quad + \eta^2_2 \big\| u_{hi} \big\|_{L^4L^{\frac{2n}{n-1}}} \big\|\nabla P_{lo}\big( |\nabla|^{-(n-4)}|u_{lo} |^2  u_{lo}\big) \big\|_{L^{\frac43}L^{\frac{2n}{n+1}}}\\
  & \quad + \eta^2_2  \big\| u_{hi}\big\|_{L^3L^{\frac{6n}{3n-4}}} \big\| P_{lo}\big( |\nabla|^{-(n-4)}|u_{hi} |^2 u_{hi}\big) \big\|_{L^{\frac32} L^{\frac{6n}{3n+4}}}\\
  & \lesssim \eta^2_2 C_0 \eta_1 .
\endaligned
\end{equation*}

We turn now to the error terms $(\ref{ine11})$ through
$(\ref{ine14})$ . We notice that  they are of the form $\langle
|u_{hi}|^2* \frac{1}{|x|}, f \rangle$ where
\begin{equation*}
\aligned f= \left\{\begin{array}{rl}
\big|\nabla\big|^{-(n-4)}|u_{lo}|^2 \big(
|u_{lo}||u_{hi}|+|u_{hi}|^2\big) \quad \text{in} \ (\ref{ine11}),\\
\big|\nabla\big|^{-(n-4)}(|u_{lo}| \ |u_{hi}|)  \big(|u_{lo}|^2+
|u_{lo}||u_{hi}| + |u_{hi}|^2 \big) \quad \text{in} \ (\ref{ine12}) ,\\
\big|\nabla\big|^{-(n-4)}|u_{hi}|^2 \big(|u_{lo}|^2+
|u_{lo}||u_{hi}|\big) \quad \text{in} \ (\ref{ine13}) ,\\
\big|u_{hi} P_{lo}(\big|\nabla\big|^{-(n-4)}|u_{hi}|^2 u_{hi})| \quad \text{in} \ (\ref{ine14}). \\
\end{array}
 \right.
\endaligned
\end{equation*}

Let us first note that as $u_{hi} \in L^3L^{\frac{6n}{3n-4}}$ and
$u_{hi} \in L^{\infty}L^{2}$,  we also have $|u_{hi}|^2 \in
L^3L^{\frac{3n}{3n-2}}$. Therefore, by the Hardy-Littlewood-Sobolev
inequality, we have $|u_{hi}|^2* \frac{1}{|x|} \in L^3L^{3n}$ and
\begin{equation*}
\aligned \langle |u_{hi}|^2* \frac{1}{|x|}, f \rangle & \lesssim
\big\| |u_{hi}|^2* \frac{1}{|x|} \big\|_{L^3L^{3n}} \big\|f
\big\|_{L^{\frac32}L^{\frac{3n}{3n-1}}} \\
& \lesssim \big\| u_{hi} \big\|_{L^3L^{\frac{6n}{3n-4}}} \big\|
u_{hi} \big\|_{L^{\infty}L^{2}} \big\|f
\big\|_{L^{\frac32}L^{\frac{3n}{3n-1}}}   \lesssim
(C_0\eta_1)^{\frac13} \eta_2 \big\|f
\big\|_{L^{\frac32}L^{\frac{3n}{3n-1}}}.
\endaligned
\end{equation*}

Consider the case of $(\ref{ine11})$, that is
$f=\big|\nabla\big|^{-(n-4)}|u_{lo}|^2 \big(
|u_{lo}||u_{hi}|+|u_{hi}|^2\big)$. By the H\"{o}lder inequality,
Hardy-Littlewood-Sobolev inequality, Proposition $\ref{imsc}$, and
Remark $\ref{rem}$, we estimate
\begin{equation*}
\aligned (\ref{ine11}) & \lesssim (C_0\eta_1)^{\frac13} \eta_2
\big\|\big|\nabla\big|^{-(n-4)}|u_{lo}|^2 \big(
|u_{lo}||u_{hi}|+|u_{hi}|^2\big)
\big\|_{L^{\frac32}L^{\frac{3n}{3n-1}}} \\
& \lesssim (C_0\eta_1)^{\frac13} \eta_2
\big\|u_{lo}\big\|^2_{L^6L^{\frac{6n}{3n-8}}}
\big\|u_{hi}\big\|_{L^3L^{\frac{6n}{3n-4}}} \big(
\big\|u_{lo}\big\|_{L^{\infty}L^{\frac{2n}{n-2}}}
+ \big\|u_{hi}\big\|_{L^{\infty}L^{\frac{2n}{n-2}}}  \big) \\
& \lesssim (C_0\eta_1)^{\frac13} \eta_2 (C_0\eta_1)^{\frac13}
(C_0\eta_1)^{\frac13}   \lesssim C_0\eta_1 \eta_2.\endaligned
\end{equation*}

Consider next the error term $(\ref{ine12})$, that is,
$f=\big|\nabla\big|^{-(n-4)}(|u_{lo}| \ |u_{hi}|) \big(|u_{lo}|^2+
|u_{lo}||u_{hi}| + |u_{hi}|^2 \big)$. By the H\"{o}lder inequality,
Hardy-Littlewood-Sobolev inequality, Proposition $\ref{imsc}$, and
Remark $\ref{rem}$, we have
\begin{equation*}
\aligned
(\ref{ine12}) & \lesssim (C_0\eta_1)^{\frac13} \eta_2
\big\|\big|\nabla\big|^{-(n-4)}(|u_{lo}| \ |u_{hi}|)
\big(|u_{lo}|^2+ |u_{lo}||u_{hi}| + |u_{hi}|^2 \big)
\big\|_{L^{\frac32}L^{\frac{3n}{3n-1}}} \\
& \lesssim (C_0\eta_1)^{\frac13} \eta_2 \big\| u_{lo} \big\|_{L^3L^{\frac{6n}{3n-10}}} \big\|u_{lo}\big\|^2_{L^{\infty}L^{\frac{2n}{n-2}}} \big\| u_{hi} \big\|_{L^3L^{\frac{6n}{3n-4}}}\\
& \quad + (C_0\eta_1)^{\frac13} \eta_2 \big\| u_{lo} \big\|_{L^3L^{\frac{6n}{3n-10}}} \big\|u_{lo}\big\|_{L^{\infty}L^{\frac{2n}{n-2}}} \big\|u_{hi}\big\|_{L^{\infty}L^{\frac{2n}{n-2}}} \big\| u_{hi} \big\|_{L^3L^{\frac{6n}{3n-4}}}\\
 & \quad + (C_0\eta_1)^{\frac13} \eta_2 \big\| u_{lo}
\big\|_{L^3L^{\frac{6n}{3n-10}}} \big\|u_{hi}
\big\|^2_{L^{\infty}L^{\frac{2n}{n-2}}} \big\| u_{hi}
\big\|_{L^3L^{\frac{6n}{3n-4}}}  \lesssim C_0\eta_1 \eta_2.
\endaligned
 \end{equation*}

Similarly, we can estimate the error term $(\ref{ine13})$, that is,
\begin{equation*} \aligned
(\ref{ine13})
& \lesssim (C_0\eta_1)^{\frac13} \eta_2
\big\|\big|\nabla\big|^{-(n-4)}|u_{hi}|^2 \big(|u_{lo}|^2+
|u_{lo}||u_{hi}|\big)
\big\|_{L^{\frac32}L^{\frac{3n}{3n-1}}} \\
& \lesssim (C_0\eta_1)^{\frac13} \eta_2 \big\| u_{lo} \big\|_{L^3L^{\frac{6n}{3n-10}}} \big\|u_{lo}\big\|_{L^{\infty}L^{\frac{2n}{n-2}}} \big\|u_{hi}\big\|_{L^{\infty}L^{\frac{2n}{n-2}}} \big\| u_{hi} \big\|_{L^3L^{\frac{6n}{3n-4}}}\\
 & \quad + (C_0\eta_1)^{\frac13} \eta_2 \big\| u_{lo}
\big\|_{L^3L^{\frac{6n}{3n-10}}} \big\|u_{hi}
\big\|^2_{L^{\infty}L^{\frac{2n}{n-2}}} \big\| u_{hi}
\big\|_{L^3L^{\frac{6n}{3n-4}}}  \lesssim C_0\eta_1 \eta_2.
\endaligned
 \end{equation*}

The last error term left to estimate is $(\ref{ine14})$. Using the
Bernstein estimate, Proposition $\ref{imsc}$, and Remark
$\ref{rem}$, we obtain
 \begin{equation*}
\aligned (\ref{ine14}) & \lesssim (C_0\eta_1)^{\frac13} \eta_2
\big\|\big|u_{hi} P_{lo}(\big|\nabla\big|^{-(n-4)}|u_{hi}|^2
u_{hi})|
\big\|_{L^{\frac32}L^{\frac{3n}{3n-1}}} \\
& \lesssim (C_0\eta_1)^{\frac13} \eta_2 \big\| u_{hi}
\big\|_{L^3L^{\frac{6n}{3n-4}}}
\big\|P_{lo}(\big|\nabla\big|^{-(n-4)}|u_{hi}|^2 u_{hi})
\big\|_{L^3L^{\frac{6n}{3n+2}}} \\
& \lesssim (C_0\eta_1)^{\frac13} \eta_2 \big\| u_{hi}
\big\|_{L^3L^{\frac{6n}{3n-4}}}
\big\|P_{lo}(\big|\nabla\big|^{-(n-4)}|u_{hi}|^2 u_{hi})
\big\|_{L^3L^{\frac{6n}{3n+8}}} \\
& \lesssim (C_0\eta_1)^{\frac13} \eta_2 \big\| u_{hi}
\big\|^2_{L^3L^{\frac{6n}{3n-4}}} \big\|
u_{hi}\big\|^2_{L^{\infty}L^{\frac{2n}{n-2}}}  \lesssim C_0\eta_1
\eta_2.
\endaligned
 \end{equation*}

 Collecting all the above estimates, we obtain that all the error
 terms on the right-hand side of Proposition $\ref{imee}$ are
 controlled by $\eta_1$. Upon rescaling, this concludes the proof of
 Proposition $\ref{flime}$ in all dimensions $n\geq 5$.

\section{Preventing energy evacuation
}\label{preventenergyevacuation}
\setcounter{section}{7}\setcounter{equation}{0} We now prove
Proposition $\ref{ecneflf}$ with the aid of almost conservation law
of frequency localized mass just as in \cite{CKSTT07}, \cite{RyV05}
and \cite{Vi05}. By the scaling $(\ref{scaling})$, we may take
$N_{min}=1$.

\subsection{The setup and contradiction argument}
Since $N(t) \in 2^{\mathbb{Z}}$, there exists $t_{min}\in I_0$ such
that $N(t_{min})=N_{min}=1$.

At time $t=t_{min}$, we have a considerable amount of mass at medium
frequencies:
\begin{equation}\label{massconcen}
\big\| P_{c(\eta_0)< \cdot < C(\eta_0)}u(t_{min})
\big\|_{L^2_x}\gtrsim c(\eta_0).
\end{equation}
However, by the Bernstein estimate, there is not much mass at
frequencies higher than $C(\eta_0)$
\begin{equation*}
\big\| P_{> C(\eta_0)}u(t_{min}) \big\|_{L^2_x}\lesssim c(\eta_0).
\end{equation*}

Let's assume for a contradiction that there exists  $t_{evac} \in
I_0$ such that $N(t_{evac}) \gg C(\eta_4)$. By time reversal
symmetry, we may assume $t_{min}< t_{evac}$. If $C(\eta_4)$ is
sufficiently large, we then see from Corollary
$\ref{frequencylocalization}$ that energy has been almost entirely
evacuated from low and medium frequencies at time $t_{evac}$:
\begin{equation}\label{lowenergycontra}
\big\| P_{< \frac{1}{\eta_4}} u(t_{evac})\big\|_{\dot{H}^1_x} \leq
\eta_4.
\end{equation}

We define
\begin{equation*}
u_{lo}(t)=P_{<\eta^{10n}_3}u(t), \quad  u_{hi}(t)=P_{\geq
\eta^{10n}_3}u(t).
\end{equation*}
Then by $(\ref{massconcen})$,
\begin{equation}\label{himassc}
\big\|u_{hi}(t_{min}) \big\|_{L^2_x}\geq \eta_1.
\end{equation}
Suppose we could show that a big portion of the mass sticks around
until time $t_{evac},$ i.e.,
\begin{equation}\label{malmostc}
\big\|u_{hi}(t_{evac})  \big\|_{L^2_x}\geq \frac12\eta_1.
\end{equation}
Since we have by the Bernstein estimate
\begin{equation*}
\big\|P_{>C(\eta_1)}u_{hi}(t_{evac})  \big\|_{L^2_x}\leq c(\eta_1),
\end{equation*}
then the triangle inequality would imply
\begin{equation*}
\big\|P_{\leq C(\eta_1)}u_{hi}(t_{evac})  \big\|_{L^2_x} \geq
\frac14\eta_1.
\end{equation*}
Another application of the Bernstein estimate would give
\begin{equation*}
\big\|P_{\leq C(\eta_1)}u(t_{evac})  \big\|_{\dot{H}^1_x} \gtrsim
c(\eta_1,\eta_3),
\end{equation*}
which would contradict $(\ref{lowenergycontra})$ if $\eta_4$ were
chosen sufficiently small.

It therefore remains to show $(\ref{malmostc})$. In order to prove
it we assume that there exists a time $t_*$ such that $t_{min}\leq
t_* \leq t_{evac}$ and
\begin{equation}\label{bootstrap3}
\aligned \inf_{t_{min}\leq t \leq t_*} \big\|u_{hi}(t)
\big\|_{L^2_x}\geq \frac12\eta_1.
\endaligned
\end{equation}
We will show that this can be bootstrapped to
\begin{equation}\label{bootstrap3g}
\aligned \inf_{t_{min}\leq t \leq t_*} \big\|u_{hi}(t)
\big\|_{L^2_x}\geq \frac34\eta_1.
\endaligned
\end{equation}
Hence, $\{t_*\in [t_{min}, t_{evac}]: (\ref{bootstrap3}) \
\text{holds} \}$ is both open and closed in $[t_{min}, t_{evac}]$
and $(\ref{malmostc})$ holds.

In order to show that $(\ref{bootstrap3})$ implies
$(\ref{bootstrap3g})$, we will treat the $L^2_x$-norm of $u_{hi}$ as
an almost conserved quantity. Define
\begin{equation*}
L(t)=\int_{\mathbb{R}^n} \big| u_{hi}(t,x)\big|^2 dx.
\end{equation*}
By $(\ref{himassc})$, we have $L(t_{min}) \geq \eta^2_1$. Hence, by
the Foundamental Theorem of Calculus it suffices to show that
\begin{equation*}
\int^{t_*}_{t_{min}} \big| \partial_t L(t) \big| dt \leq
\frac{1}{100}\eta^2_1.
\end{equation*}
Since
\begin{equation*}
\aligned
\partial_t L(t) & = 2\int_{\mathbb{R}^n} \big\{P_{hi}(|\nabla|^{-(n-4)}|u|^2u), u_{hi} \big\}_{m}
dx  = 2\int_{\mathbb{R}^n}
\big\{P_{hi}(|\nabla|^{-(n-4)}|u|^2u)-|\nabla|^{-(n-4)}|u_{hi}|^2u_{hi}
, u_{hi} \big\}_{m} dx,
\endaligned
\end{equation*}
we need to show
\begin{equation}\label{macl}
\int^{t_*}_{t_{min}} \big| \int_{\mathbb{R}^n}
\big\{P_{hi}(|\nabla|^{-(n-4)}|u|^2u)-|\nabla|^{-(n-4)}|u_{hi}|^2u_{hi}
, u_{hi} \big\}_{m} dx \big| dt \leq \frac{1}{100}\eta^2_1.
\end{equation}

In order to prove $(\ref{macl})$, we need to control the various
interactions between low, medium, and high frequencies. In the next
section we will make some preliminary estimates that will make this
goal possible.

\subsection{Spacetime estimates for high, medium, and low frequencies}
Remember that the frequency-localized interaction Morawetz
inequality implies that for $N < c(\eta_2)N_{min}$,
\begin{equation*}
\aligned
   \int^{t_{evac}}_{t_{min}}\int_{\mathbb{R}^n}\int_{\mathbb{R}^n} \frac{\big| P_{\geq N}u(t,y) \big|^2\big| P_{\geq N}u(t,x) \big|^2 }{|x-y|^3} \ dxdydt
   \lesssim \eta_1 N^{-3}.
\endaligned
\end{equation*}

This estimate is useful for medium and high frequencies; however it
is extremely bad for low frequencies since $N^{-3}$ gets
increasingly larger as $N \rightarrow 0$. We therefore need to
develop better estimates in this case. Since $u_{\leq \eta_3}$ has
extremely small energy at $t=t_{evac}$ (see
$(\ref{lowenergycontra})$), we expect it to have small energy at all
time in $[t_{min}, t_{evac}]$. Of course, there is energy leaking
from the high frequencies to the low frequencies, but the
interaction Morawetz estimate limits this leakage. Indeed, we have

\begin{proposition}\label{lowcon}
 under the assumptions above,
\begin{equation*}
\big\| P_{\leq N} u\big\|_{\dot{S}^1([t_{min}, t_{evac}]\times
\mathbb{R}^n)} \lesssim \eta_4 + \eta_0 \eta_3^{-2}N^{2},
\end{equation*}
for all $N\leq \eta_3$.
\end{proposition}

{\bf Proof: } Consider the set
\begin{equation*}
\Omega = \big\{ t\in [t_{min}, t_{evac}): \big\| P_{\leq N}
u\big\|_{\dot{S}^1([t , t_{evac}]\times \mathbb{R}^n)} \leq
C_0\eta_4 + \eta_0 \eta_3^{-2}N^{ 2}  \big\}
\end{equation*}
where $C_0$ is a large constant to be chosen later and not depending
on any of the $\eta$'.

Our goal is to show that $t_{min} \in \Omega$. First we can show
that $t\in \Omega$ for $t$ close to $t_{evac}$.

Now suppose that $t\in \Omega$. We will show that
\begin{equation}\label{goal1}
\big\| P_{\leq N} u\big\|_{\dot{S}^1([t , t_{evac}]\times
\mathbb{R}^n)} \leq \frac12 C_0\eta_4 + \frac12 \eta_0 \eta_3^{-
2}N^{ 2}
\end{equation}
holds for any $ N \leq \eta_3$. thus, $\Omega$ is both open and
closed in $[t_{min}, t_{evac}]$ and we have $t_{min} \in \Omega$ as
desired.

Fixing $N \leq \eta_3$, the Strichartz estimate implies
\begin{equation*}
\big\| P_{\leq N} u\big\|_{\dot{S}^1([t , t_{evac}]\times
\mathbb{R}^n)} \lesssim \big\| P_{\leq
N}u(t_{evac})\big\|_{\dot{H}^1_x} +  \big\| \nabla P_{\leq N}\big(
\big|\nabla \big|^{-(n-4)}|u|^2u
\big)\big\|_{L^{\frac32}L^{\frac{6n}{3n+4}}}.
\end{equation*}

By $(\ref{lowenergycontra})$, we have
\begin{equation}\label{lowfrele}
\big\| P_{\leq N}u(t_{evac})\big\|_{\dot{H}^1_x} \lesssim \eta_4
\end{equation}
which is acceptable for $(\ref{goal1})$ if $C_0$ is chosen
sufficiently large.

To handle the nonlinearity, we decompose $u = u_{<\eta_4}+u_{\eta_4
\leq \cdot \leq \eta_3}+ u_{> \eta_3}$ and use the triangle
inequality to estimate
\begin{eqnarray}
  \big\| \nabla  P_{\leq N}\big( \big|\nabla
\big|^{-(n-4)}|u|^2u \big)\big\|_{L^{\frac32}L^{\frac{6n}{3n+4}}}
 &\lesssim &   \big\| \nabla P'_{\leq N}\big( \big|\nabla
\big|^{-(n-4)}|u_{<\eta_4}|^2|u_{< \eta_4}| \big)\big\|_{L^{\frac32}L^{\frac{6n}{3n+4}}} \label{lowfrenle1}\\
& + & \big\| \nabla P'_{\leq N}\big( \big|\nabla
\big|^{-(n-4)}|u_{<\eta_4}|^2|u_{\eta_4 \leq \cdot \leq \eta_3}| \big)\big\|_{L^{\frac32}L^{\frac{6n}{3n+4}}}  \label{lowfrenle2}\\
& + & \big\| \nabla P'_{\leq N}\big( \big|\nabla
\big|^{-(n-4)}|u_{<\eta_4}|^2|u_{> \eta_3}| \big)\big\|_{L^{\frac32}L^{\frac{6n}{3n+4}}} \label{lowfrenle3}\\
& + & \big\| \nabla P'_{\leq N}\big( \big|\nabla
\big|^{-(n-4)}|u_{\eta_4 \leq \cdot \leq \eta_3}|^2|u_{<\eta_4}| \big)\big\|_{L^{\frac32}L^{\frac{6n}{3n+4}}}  \label{lowfrenle4}\\
& + & \big\| \nabla P'_{\leq N}\big( \big|\nabla
\big|^{-(n-4)}|u_{\eta_4 \leq \cdot \leq \eta_3}|^2|u_{\eta_4 \leq \cdot \leq \eta_3}| \big)\big\|_{L^{\frac32}L^{\frac{6n}{3n+4}}}  \label{lowfrenle5}\\
& + & \big\| \nabla P'_{\leq N}\big( \big|\nabla
\big|^{-(n-4)}|u_{\eta_4 \leq \cdot \leq \eta_3}|^2|u_{>\eta_3}| \big)\big\|_{L^{\frac32}L^{\frac{6n}{3n+4}}}  \label{lowfrenle6}\\
& + & \big\| \nabla  P'_{\leq N}\big( \big|\nabla
\big|^{-(n-4)}|u_{>\eta_3}|^2|u |
\big)\big\|_{L^{\frac32}L^{\frac{6n}{3n+4}}}.  \label{lowfrenle7}
\end{eqnarray}

Using the bootstrap hypothesis $t\in \Omega$, we have
\begin{eqnarray}\label{mede}
  \big\| u_{ \eta_4\leq \cdot \leq \eta_3}\big\|_{\dot{S}^0}
&\lesssim & \sum_{\eta_4 \leq M \leq \eta_3} \big\|
P_Mu\big\|_{\dot{S}^0}   \lesssim   \sum_{\eta_4 \leq
M \leq \eta_3} M^{-1} \big\|\nabla P_Mu\big\|_{\dot{S}^0} \nonumber\\
& \lesssim & \sum_{\eta_4 \leq M \leq \eta_3} M^{-1} \big\|
P_Mu\big\|_{\dot{S}^1}  \lesssim  \sum_{\eta_4 \leq M \leq \eta_3}
M^{-1} \big( C_0\eta_4 + \eta_0 \eta^{- 2}_3 M^{ 2} \big)   \lesssim
\eta_0 \eta^{-1}_3. \label{mede}
 \end{eqnarray}

We now turn to estimate the nonlinearity. Using the bootstrap
hypothesis, $t\in \Omega$, we estimate
\begin{equation*}
\aligned (\ref{lowfrenle1}) & \lesssim \big\| \nabla P'_{\leq
N}\big( \big|\nabla \big|^{-(n-4)}|u_{<\eta_4}|^2|u_{< \eta_4}|
\big)\big\|_{L^{\frac32}L^{\frac{6n}{3n+4}}} \\
& \lesssim \big\| \nabla u_{<\eta_4}
\big\|_{L^3L^{\frac{6n}{3n-4}}}\big\|u_{<\eta_4}\big\|^2_{L^6L^{\frac{6n}{3n-8}}}
 \lesssim \big\|u_{<\eta_4}\big\|^3_{\dot{S}^1} \lesssim \big(
C_0\eta_4 + \eta_0 \eta_3^{- 2}\eta^{2}_4 \big)^3 \lesssim \eta_4,
\endaligned
\end{equation*}
which again is acceptable for $(\ref{goal1})$ provided $C_0$ is
sufficiently large.

By the Bernstein estimate, $(\ref{mede})$, Corollary $\ref{himora}$
and $t\in \Omega$, we obtain
\begin{equation*}
\aligned (\ref{lowfrenle2}) & \lesssim  N \big\| \big|\nabla
\big|^{-(n-4)}|u_{<\eta_4}|^2|u_{ \eta_4\leq \cdot \leq \eta_3}|
\big\|_{L^{\frac32}L^{\frac{6n}{3n+4}}} \\
&  \lesssim  N \big\| u_{<\eta_4}
 \big\|^2_{L^6L^{\frac{6n}{3n-8}}}\big\| u_{ \eta_4\leq \cdot \leq
\eta_3}\big\|_{L^3L^{\frac{6n}{3n-4}}}\\
&  \lesssim  N \big( C_0 \eta_4 + \eta_0 \eta^{- 2}_3 \eta^{ 2}_4
\big)^2 \eta_0
 \eta^{-1}_3  \lesssim \eta_4,
\endaligned
\end{equation*}
and
\begin{equation*}
\aligned
  (\ref{lowfrenle3})
& \lesssim N \big\| \big|\nabla \big|^{-(n-4)}|u_{<\eta_4}|^2|u_{>
\eta_3}| \big\|_{L^{\frac32}L^{\frac{6n}{3n+4}}}  \\
& \lesssim  N \big\| u_{<\eta_4}
 \big\|^2_{L^6L^{\frac{6n}{3n-8}}}\big\| u_{ >
\eta_3}\big\|_{L^3L^{\frac{6n}{3n-4}}}\\
& \lesssim  N \big( C_0 \eta_4 + \eta_0 \eta^{- 2}_3 \eta^{ 2}_4
 \big)^2 \eta^{\frac13}_1 \eta^{-1}_3
\lesssim \eta_4,
\endaligned
\end{equation*}
which again is acceptable for $(\ref{goal1})$ provided $C_0$ is
sufficiently large.

By the Bernstein estimate, $(\ref{mede})$, Corollary $\ref{himora}$
and $t\in \Omega$, we have
\begin{equation*}
\aligned (\ref{lowfrenle4}) & \lesssim  N^2 \big\| \big|\nabla
\big|^{-(n-4)}|u_{ \eta_4\leq \cdot \leq \eta_3}|^2 |u_{<\eta_4}|
\big\|_{L^{\frac32}L^{\frac{6n}{3n+10}}} \\
&  \lesssim  N^2 \big\| u_{<\eta_4}
 \big\|_{L^{\infty}L^{\frac{2n}{n-2}}} \big\| u_{ \eta_4\leq \cdot \leq
\eta_3}\big\|^2_{L^3L^{\frac{6n}{3n-4}}}\\
&  \lesssim  N^2 \big( C_0 \eta_4 + \eta_0 \eta^{- 2}_3 \eta^{ 2}_4
\big) \big( \eta_0
 \eta^{-1}_3 \big)^2 \lesssim \eta_4,\\
(\ref{lowfrenle5}) & \lesssim  N^2 \big\| \big|\nabla
\big|^{-(n-4)}|u_{ \eta_4\leq \cdot \leq \eta_3}|^2 u_{\eta_4\leq
\cdot \leq \eta_3}
\big\|_{L^{\frac32}L^{\frac{6n}{3n+10}}} \\
&  \lesssim  N^2 \big\| u_{\eta_4\leq \cdot \leq \eta_3 }
 \big\|_{L^{\infty}L^{\frac{2n}{n-2}}} \big\| u_{ \eta_4\leq \cdot \leq
\eta_3}\big\|^2_{L^3L^{\frac{6n}{3n-4}}}\\
&  \lesssim  N^2   \big( \eta_0
 \eta^{-1}_3 \big)^2 = \eta^2_0 \eta^{-2}_3 N^2,
 \endaligned
\end{equation*}
and
 \begin{equation*}
\aligned  (\ref{lowfrenle6})& \lesssim N^2 \big\| \big|\nabla
\big|^{-(n-4)}|u_{ \eta_4\leq \cdot \leq \eta_3}|^2u_{>
\eta_3} \big\|_{L^{\frac32}L^{\frac{6n}{3n+10}}}  \\
& \lesssim  N^2 \big\| u_{\eta_4\leq \cdot \leq \eta_3 }
 \big\|_{L^6L^{\frac{6n}{3n-8}}} \big\| u_{\eta_4\leq \cdot \leq \eta_3 }
 \big\|_{L^6L^{\frac{6n}{3n-2}}} \big\| u_{ >
\eta_3}\big\|_{L^3L^{\frac{6n}{3n-4}}}\\
& \lesssim  N^2 \eta_0 \eta^{-1}_3\eta^{\frac13}_1 \eta^{-1}_3 =
\eta_0 \eta^{\frac13}_1 \eta^{-2}_3 N^2,
\endaligned
\end{equation*}
which again is acceptable for $(\ref{goal1})$ provided $C_0$ is
sufficiently large.

By the Bernstein estimate, $(\ref{mede})$, Corollary $\ref{himora}$
and $t\in \Omega$, we have
\begin{equation*}
\aligned  (\ref{lowfrenle7})& \lesssim  N^2 \big\|
 \big|\nabla \big|^{-(n-4)}|u_{\geq \eta_3}|^2u
 \big\|_{L^{\frac32}L^{\frac{6n}{3n+10}}} \\
 & \lesssim  N^2 \big\| u \big\|_{L^{\infty}L^{\frac{2n}{n-2}}}\big\| u_{>
  \eta_3}\big\|^2_{L^3L^{\frac{6n}{3n-4}}}  \lesssim  N^2  (\eta^{\frac13}_1 \eta^{-1}_3)^2 =
\eta^{\frac23}_1\eta^{-2}_3N^2
\endaligned
\end{equation*}
which again is acceptable for $(\ref{goal1})$ provided $C_0$ is
sufficiently large.

The proposition is complete.

\subsection{Controlling the localized $L^2$ mass increment}
We now have good enough control over low, medium, and high
frequencies to prove $(\ref{macl})$. Writing
\begin{equation*}
\aligned
P_{hi}\big(|\nabla|^{-(n-4)}|u|^2u\big)-|\nabla|^{-(n-4)}|u_{hi}|^2u_{hi}
=
&   P_{hi}\big(|\nabla|^{-(n-4)}|u|^2u-|\nabla|^{-(n-4)}|u_{hi}|^2u_{hi} - |\nabla|^{-(n-4)}|u_{lo}|^2u_{lo}\big) \\
& - P_{lo}\big(|\nabla|^{-(n-4)}|u_{hi}|^2u_{hi} \big)  +
P_{hi}\big( |\nabla|^{-(n-4)}|u_{lo}|^2u_{lo} \big).
\endaligned
\end{equation*}
Clearly, we only have to consider the following terms
\begin{eqnarray}
\int^{t_*}_{t_{min}} \Big| \int_{\mathbb{R}^n} \overline{u}_{hi}
P_{hi}(|\nabla|^{-(n-4)}|u|^2u-|\nabla|^{-(n-4)}|u_{hi}|^2u_{hi} -
|\nabla|^{-(n-4)}|u_{lo}|^2u_{lo})dx\Big| dt. \label{mine1}\\
\int^{t_*}_{t_{min}} \Big| \int_{\mathbb{R}^n} \overline{u}_{hi}
P_{lo}\big(|\nabla|^{-(n-4)}|u_{hi}|^2u_{hi} \big)dx\Big| dt. \label{mine2}\\
\int^{t_*}_{t_{min}} \Big| \int_{\mathbb{R}^n} \overline{u}_{hi}
P_{hi}\big( |\nabla|^{-(n-4)}|u_{lo}|^2u_{lo}  \big) dx\Big| dt.
\label{mine3}
\end{eqnarray}

{\bf $\bullet$ Case 1. Estimation of $(\ref{mine1})$.}

We move the self-adjoint operator $P_{hi}$ onto $\overline{u}_{hi}$,
and obtain
\begin{equation*}
\aligned & \
|\nabla|^{-(n-4)}|u|^2u-|\nabla|^{-(n-4)}|u_{hi}|^2u_{hi} -
|\nabla|^{-(n-4)}|u_{lo}|^2u_{lo}  \\
= & \quad \ |\nabla|^{-(n-4)}|u_{lo}|^2u_{hi}   +
2|\nabla|^{-(n-4)}\text{Re}(u_{lo}\overline{u}_{hi})\big(u_{lo}+u_{hi}\big)
+ |\nabla|^{-(n-4)}|u_{hi}|^2u_{lo}.
\endaligned
\end{equation*}

We first consider the contribution of
$|\nabla|^{-(n-4)}|u_{lo}|^2u_{hi}$. By Corollary $\ref{himora}$ and
Proposition $\ref{lowcon}$,  we have
\begin{equation*}
\aligned
 \int^{t_*}_{t_{min}} \int_{\mathbb{R}^n} \big|
P_{hi}u_{hi}\big|
 |\nabla|^{-(n-4)}|u_{lo}|^2 \big|u_{hi}\big|dx dt & \lesssim \big\| P_{hi}u_{hi}
 \big\|_{L^3L^{\frac{6n}{3n-4}}} \big\| u_{hi}
 \big\|_{L^3L^{\frac{6n}{3n-4}}} \big\|
 u_{lo}\big\|^2_{L^6L^{\frac{6n}{3n-8}}} \\
 & \lesssim \big( \eta_1 (\eta^{10n}_3)^{-3} \big)^{\frac23} \big(
 \eta_4+  (\eta^{-1}_3\eta^{10n}_3)^{ 2}\big)^2 \\
 & \lesssim \eta^{\frac23}_1 \eta^{20n-4}_3\ll \eta^2_1.
 \endaligned
 \end{equation*}

 We now turn towards the contribution of
 $|\nabla|^{-(n-4)}|u_{hi}|^2u_{lo}$. We decompose $u_{hi}=u_{\eta^{10n}_3\leq \cdot \leq \eta_3} +
 u_{>\eta_3}$ and obtain by the Berstein estimate, Corollary $\ref{himora}$ and
Proposition $\ref{lowcon}$,
\begin{equation*}
\aligned \big\| u_{>\eta_3}\big\|_{L^3L^{\frac{6n}{3n-4}}} &
\lesssim
\eta^{\frac13}_1 \eta^{-1}_3,\\
\big\| u_{ \eta^{10n}_3\leq \cdot \leq \eta_3
}\big\|_{L^3L^{\frac{6n}{3n-4}}} & \lesssim \sum_{\eta^{10n}_3\leq N
\leq \eta_3} \big\|u_{N}\big\|_{L^3L^{\frac{6n}{3n-4}}}  \lesssim
\sum_{\eta^{10n}_3\leq N \leq \eta_3} N^{-1} \big\|\nabla
u_{N}\big\|_{L^3L^{\frac{6n}{3n-4}}} \\
&   \lesssim \sum_{\eta^{10n}_3\leq N \leq \eta_3} N^{-1} \big(
\eta_4+ \eta^{- 2}_3 N^{ 2} \big)    \lesssim
\eta^{- 2}_3\eta_3=\eta^{-1}_3,\\
\big\| u_{lo}\big\|_{L^{\infty}L^{\frac{2n}{n-4}}} & \lesssim
\eta^{10n}_3 \big\| u_{lo}\big\|_{L^{\infty}L^{\frac{2n}{n-2}}}
\lesssim \eta^{10n}_3  \big( \eta_4 + \eta^{- 2}_3 (\eta^{10n}_3)^{
2}\big)= \eta^{30n- 2}_3,
\endaligned
\end{equation*}
then
\begin{equation*}
\aligned
 \int^{t_*}_{t_{min}} \int_{\mathbb{R}^n} \big|
P_{hi}u_{hi}\big|
 |\nabla|^{-(n-4)}|u_{hi}|^2\big| u_{lo}\big| dxdt & \lesssim
 \sum^3_{j=0} \big\| u_{ \eta^{10n}_3\leq \cdot \leq \eta_3
}\big\|^{3-j}_{L^3L^{\frac{6n}{3n-4}}} \big\| u_{>\eta_3}\big\|^{
j}_{L^3L^{\frac{6n}{3n-4}}} \big\|
u_{lo}\big\|_{L^{\infty}L^{\frac{2n}{n-4}}} \\
& \lesssim
 \sum^3_{j=0} \big( \eta^{\frac13}_1 \eta^{-1}_3\big)^{j} \big( \eta^{-1}_3\big)^{3-j} \eta^{30n- 2}_3   \lesssim \eta^{30n-5}_3 \ll  \eta^2_1.
 \endaligned
 \end{equation*}

Now we consider the contribution of
$2|\nabla|^{-(n-4)}\text{Re}(u_{lo}\overline{u}_{hi})\big(u_{lo}+u_{hi}\big)$.
Similarly, we have
 \begin{equation*}
\aligned
 \int^{t_*}_{t_{min}} \int_{\mathbb{R}^n} \big|
P_{hi}u_{hi}\big|\big|
  |\nabla|^{-(n-4)}(u_{lo}\overline{u}_{hi})\big(u_{lo}+u_{hi}\big)\big| dt & \lesssim \big\| P_{hi}u_{hi}
 \big\|_{L^3L^{\frac{6n}{3n-4}}} \big\| u_{hi}
 \big\|_{L^3L^{\frac{6n}{3n-4}}} \big\|
 u_{lo}\big\|^2_{L^6L^{\frac{6n}{3n-8}}}
 \\
  + \sum^3_{j=0} & \big\| u_{ \eta^{10n}_3\leq \cdot \leq \eta_3
}\big\|^{j}_{L^3L^{\frac{6n}{3n-4}}} \big\|
u_{>\eta_3}\big\|^{3-j}_{L^3L^{\frac{6n}{3n-4}}} \big\|
u_{lo}\big\|_{L^{\infty}L^{\frac{2n}{n-4}}} \\
& \lesssim \eta^{\frac23}_1 \eta^{20n-4}_3 + \eta^{30n-5}_3 \ll
\eta^2_1.
 \endaligned
 \end{equation*}

 Therefore
 \begin{equation*}
(\ref{mine1})  \ll  \eta^2_1.
 \end{equation*}

{\bf $\bullet$ Case 2. Estimation of $(\ref{mine2})$.}

Moving the projection $P_{lo}$ onto $\overline{u}_{hi}$ and writing
$P_{lo}u_{hi}= P_{hi}u_{lo}$, and get
\begin{equation*}
\aligned (\ref{mine2}) & = \int^{t_*}_{t_{min}} \Big|
\int_{\mathbb{R}^n} P_{hi}\overline{u}_{lo}
\big(|\nabla|^{-(n-4)}|u_{hi}|^2u_{hi} \big)dx\Big| dt \\
& \lesssim
 \sum^3_{j=0} \big\| u_{ \eta^{10n}_3\leq \cdot \leq \eta_3
}\big\|^{j}_{L^3L^{\frac{6n}{3n-4}}} \big\|
u_{>\eta_3}\big\|^{3-j}_{L^3L^{\frac{6n}{3n-4}}} \big\|P_{hi}
u_{lo}\big\|_{L^{\infty}L^{\frac{2n}{n-4}}} \\
& \lesssim
 \sum^3_{j=0} \big\| u_{ \eta^{10n}_3\leq \cdot \leq \eta_3
}\big\|^{j}_{L^3L^{\frac{6n}{3n-4}}} \big\|
u_{>\eta_3}\big\|^{3-j}_{L^3L^{\frac{6n}{3n-4}}} \big\|
u_{lo}\big\|_{L^{\infty}L^{\frac{2n}{n-4}}}   \lesssim
\eta^{30n-5}_3 \ll \eta^2_1.
\endaligned
\end{equation*}

{\bf $\bullet$ Case 3. Estimation of $(\ref{mine3})$.}

By Bernstein estimate, Corollary $\ref{himora}$ and Proposition
$\ref{lowcon}$, we have
\begin{equation*}
\aligned (\ref{mine3}) & \lesssim \big\|
u_{hi}\big\|_{L^3L^{\frac{6n}{3n-4}}} \big\| P_{hi}\big(
|\nabla|^{-(n-4)}|u_{lo}|^2u_{lo} \big) \big\|_{L^{\frac32}L^{\frac{6n}{3n+4}}}\\
& \lesssim \big\| u_{hi}\big\|_{L^3L^{\frac{6n}{3n-4}}} \big\|
\nabla P_{hi}\big(
|\nabla|^{-(n-4)}|u_{lo}|^2u_{lo} \big) \big\|_{L^{\frac32}L^{\frac{6n}{3n+4}}}\\
& \lesssim \big\| u_{hi}\big\|_{L^3L^{\frac{6n}{3n-4}}} \big\|
\nabla u_{lo}\big\|_{L^3 L^{\frac{6n}{3n-4}}} \big\|
u_{lo}\big\|^2_{L^6L^{\frac{6n}{3n-8}}} \\
& \lesssim  \eta^{\frac13}_1 (\eta^{10n}_3)^{-1} \big(
 \eta_4+  (\eta^{-1}_3\eta^{10n}_3)^{ 2}\big)^3   =\eta^{\frac13}_1
 \eta^{50n-6}_3 \ll  \eta^2_1.
\endaligned
\end{equation*}

\textbf{Acknowledgements:} The authors thank  thank  Professor T.
Tao for his helpful suggestion. C. Miao and G. Xu are partly
supported by the NSF of China (No. 10725102, No. 10726053), and L.
Zhao is supported by China postdoctoral science foundation
project.

\begin{center}

\end{center}
\end{document}